\documentclass[aos, preprint]{imsart}

\RequirePackage[OT1]{fontenc}
\usepackage{graphicx}
\RequirePackage{amsthm,amsmath,amssymb, float}
\RequirePackage[numbers]{natbib}

\usepackage{bm, bbm}
\usepackage{color}
\usepackage{enumerate}
\usepackage{mathabx}

\usepackage{adjustbox}
\usepackage{chngpage}
\usepackage{calc}

\usepackage{multirow}
\usepackage{multicol}
\usepackage{etoolbox}
\patchcmd{\footnotemark}{\stepcounter{footnote}}{\refstepcounter{footnote}}{}{}
\RequirePackage[colorlinks,citecolor=blue,urlcolor=blue,hyperfootnotes=true]{hyperref}

\newcommand{\bigcell}[2]{\begin{tabular}{@{}#1@{}}#2\end{tabular}}

\usepackage{lastpage}

\startlocaldefs
\numberwithin{equation}{section}
\theoremstyle{plain}
\newtheorem{thm}{Theorem}[section]
\newtheorem{lem}{Lemma}[section]

\newtheorem{defn}{Definition}[section]

\newtheorem{example}{Example}[section]
\newtheorem{remark}{Remark}[section]
\endlocaldefs

\usepackage{accents}

\renewcommand{\hat}{\widehat}
\renewcommand{\tilde}{\widetilde}
\renewcommand{\check}{\widecheck}

\renewcommand{\phi}{\varphi}

\newcommand{\M}{{\cal M}}

\newcommand{\ep}{\mathbb{E}}
\newcommand{\pr}{\mathbb{P}}
\newcommand{\var}{\mathrm{Var}}
\newcommand{\cov}{\mathrm{Cov}}

\newcommand{\tblue}[1]{\textcolor{blue}{#1}}
\newcommand{\tred}[1]{\textcolor{red}{#1}}

\newcommand{\ii}{{\mathbbm{i}}}
\newcommand{\td}{{\mathrm{d}}}

\newcommand{\tildehighorderbound}{{\tOp\left({\cal M}(\rho_n,n; R)\right)}}
\newcommand{\Ohighorderbound}{{O\left({\cal M}(\rho_n,n; R)\right)}}

\newcommand{\tOp}{{\tilde{O}_p}}
\renewcommand{\top}{{\tilde{o}_p}}

\usepackage{xcolor}
	\definecolor{yuancolor}{rgb}{0.01, 0.75, 0.24}
    \definecolor{yuancolor2}{rgb}{0,0,0}
    \definecolor{yuancolor3}{rgb}{0,0,0}

\begin{document}

\begin{frontmatter}
\title{Edgeworth expansions for network moments}

\runtitle{Network Edgeworth expansion}

\begin{aug}
\author{\fnms{Yuan} \snm{Zhang}\ead[label=e1]{yzhanghf@stat.osu.edu}}
\and
\author{\fnms{Dong} \snm{Xia}\ead[label=e2]{madxia@ust.hk}}

\runauthor{Zhang and Xia}


%

\end{aug}

\begin{abstract}
	
	\label{sec::abstract}Network method of moments \citep{bickel2011method} is an important tool for nonparametric network inference.  However, there has been little investigation on accurate descriptions of the sampling distributions of network moment statistics.  In this paper, we present the first higher-order accurate approximation to the sampling CDF of a studentized network moment by Edgeworth expansion.  In sharp contrast to classical literature on \emph{noiseless} U-statistics, we {show} that the Edgeworth expansion of a network moment statistic as a \emph{noisy} U-statistic can achieve higher-order accuracy without non-lattice or smoothness assumptions but just requiring weak regularity conditions.  Behind this result is our surprising discovery that the two typically-hated factors in network analysis, namely, sparsity and edge-wise observational errors, jointly play a blessing role, contributing a crucial \emph{self-smoothing} effect in the network moment statistic and making it analytically tractable.  Our assumptions match the minimum requirements in related literature.  
{For sparse networks, our theory shows that a simple normal approximation achieves a gradually depreciating Berry-Esseen type bound as the network becomes sparser.  This result also refines the best previous theoretical result.}
	
For practitioners, our empirical Edgeworth expansion is highly accurate and computationally efficient.  It is also easy to implement.  {We demonstrate the clear advantage of our method by} comprehensive simulation studies.

We showcase three applications of our results in network inference.  We {prove}, to our knowledge, {the first theoretical guarantee of higher-order accuracy for} some network bootstrap {schemes, and moreover, the first} theoretical guidance for {selecting the sub-sample size for} network sub-sampling.  We also derive a one-sample test and {the} Cornish-Fisher confidence interval for a given moment with {higher-order accurate controls of confidence level and type I error, respectively}.

\end{abstract}



\end{frontmatter}

\section{Introduction}
\label{sec::introduction}

\subsection{Overview}

\emph{Network moments} are the frequencies of particular patterns, called \emph{motifs}, that repeatedly occur in networks \citep{milo2002network, alon2007network, rubinov2010complex}.   {Examples include} triangles, stars and wheels.  They provide {succinct and} informative sketches of potentially very high-dimensional network population distributions.  Pioneered by \cite{bickel2011method, lovasz2009very}, the \emph{method of moments} for network data has become a powerful tool for frequentist nonparametric network inferences \citep{ambroise2012new, maugis2017statistical, wegner2018identifying, ali2016comparison, lunde2019subsampling, matsushita2020jackknife}.  Compared to model-based network inference methods \citep{lei2016goodness, bickel2016hypothesis, wang2017likelihood, li2018two}, moment method enjoy{s} several unique advantages.

First, {network moments play important roles in network modeling.}
They are the building blocks of the well-known exponential random graph models (ERGM) \citep{hunter2008ergm, yan2016asymptotics}.  {More generally, under an exchangeable network assumption,} the deep theory by \cite{bickel2011method} (Theorem 3) and \cite{borgs2010moments} (Theorem 2.1) show that knowing all population moments can uniquely determine the network model up to weak isomorphism, despite no explicit inversion formula is yet available.  
{From an inference perspective,} the evaluation of network moments is completely model-free,
making them objective evidences for the specification{, validation} and comparison of network models \citep{broido2019scale, Seshadhri201911030, tsiotas2019detecting, ouadah2019degree}. 
Second, {network moments can be very efficiently computed, easily allowing parallel computing. 
This is a crucial advantage} in a big data era, where business and industry networks could contain $10^5\sim 10^7$ or even more nodes \citep{clauset2004finding, snapnets} and computation efficiency becomes a substantive practicality concern.
Model-fitting based network inference methods might face challenges in handling huge networks, while moment method equipped with proper sampling techniques \citep{rohe2019critical, crane2018probabilistic} {will scale more comfortably (also see our comment in Section \ref{sec::discussion})}.
Third, many network moments {and their derived functionals} are {important structural features of great practical} interest.  Examples include clustering coefficient \citep{holland1971transitivity, watts1998collective}, degree distribution \citep{prvzulj2007biological, stephen2009explaining}, transitivity \citep{rohe2013blessing}, and {more listed in Table A.1 in \cite{rubinov2010complex}.}

Despite the {importance and raising interest in} network moment method, the answer to the following core question remains under-explored:
\begin{center}
	\smallskip
	\emph{What is the sampling distribution of a network moment?}
	\smallskip
\end{center}
For a given network motif ${\cal R}$, let $\hat{U}_n$ denote its sample relative frequency {(see \eqref{def::h(A)} for a formal definition)} with expectation $\mu_n:=\ep[\hat{U}_n]$.  Let $\hat{S}_n^2$ be an estimator of $\var(\hat{U}_n)$ to be specified later.  We are mainly interested in finding the distribution of the studentized form $\hat{T}_n := (\hat{U}_n - \mu_n)/\hat{S}_n$.  It is well-known that under the widely-studied \emph{exchangeable network} model, $\hat{T}_n\stackrel{d}{\to} N(0,1)$ uniformly \citep{bickel2011method, bhattacharyya2015subsampling, green2017bootstrapping}, but usually, $N(0,1)$ only provides a {rough characterization of $F_{\hat T_n}$}, and one naturally yearns for a finer approximation.  To this end, several network bootstrap methods have been recently proposed \citep{bickel2011method, bhattacharyya2015subsampling, green2017bootstrapping, levin2019bootstrapping, lunde2019subsampling} in an attempt to address this question, and they quickly inspired many follow-up works \cite{thompson2016using, tang2017semiparametric, gel2017bootstrap, chen2019bootstrap}, which clearly reflects the need of {an accurate approximation by data analysts}.  However, compared to their empirical effectiveness, the theoretical {foundation} of network bootstraps remains weak.  Almost all existing justifications of network bootstraps critically depend on the following type of results
{
\begin{align*}
    |\hat{U}_n^*-\hat{U}_n|=o_p(n^{-1/2}),
    &\quad \textrm{ and }\quad
    |\hat{U}_n-U_n|=o_p(n^{-1/2});\\
    \textrm{ or similarly, }\left|\hat{T}_n^* - \hat{T}_n\right| = o_p(1),
    &\quad \textrm{ and }\quad
    |\hat{T}_n-T_n|=o_p(1);
\end{align*}}
where $\hat{U}_n^*$ or $\hat{T}_n^*$ are bootstrapped statistics.  {Then the validity of network bootstraps is implied by the {well-known} asymptotic normality of $ U_n$ or $T_n$} \citep{bhattacharyya2015subsampling, green2017bootstrapping, lunde2019subsampling}.
However, this approach cannot show whether network bootstraps have any accuracy advantage over a simple normal approximation, especially considering the much higher computational costs of bootstraps.

In this paper, we propose the first provable \emph{higher-order} approximation to the sampling distribution of a given studentized network moment.  
{To our knowledge, we are the first to realize that in fact the noisy $\hat U_n$ and $\hat T_n$ are usually more analytically tractable than the noiseless versions $U_n$ and $T_n$.  This sharply contrasts the existing literature on network bootstraps that attempt to prove the asymptotic properties of $\hat U_n$ by reducing it to $U_n$.}
We briefly summarize our main results into an informal theorem {here.  It turns out that the error bound depends on the shape of the motif $R$.}
{
\begin{defn}[Acyclic and cyclic motifs, see also \cite{bickel2011method,bhattacharyya2015subsampling, green2017bootstrapping, levin2019bootstrapping}]
    \label{definition::acyclic-cyclic-R}
	A motif $R$ is called \emph{acyclic}, if its edge set is a subset of an $r$-tree.  The motif is called \emph{cyclic}, if it is \emph{connected} and contains at least one cycle.  In other words, a \emph{cyclic} motif is connected but not a tree.
\end{defn}
\begin{defn}
	\label{definition::M}
	To simplify the statements of our method's error bound under different motif shapes, especially in Table \ref{table::comparison-of-results} and proof steps, define the following shorthand
	\begin{equation}
		{\cal M}(\rho_n,n; R):=
		\begin{cases}
			 \left(\rho_n\cdot n \right)^{-1}{\cdot \log^{1/2} n+n^{-1}\cdot \log^{3/2}n}, & \textrm{ For acyclic }R\\
			 \rho_n^{-r/2}\cdot n^{-1}{\cdot \log^{1/2} n+n^{-1}\cdot \log^{3/2}n}, & \textrm{For cyclic }R
		\end{cases}
		\label{M}
	\end{equation}
\end{defn}
}
Now we are ready to present the informal statement of our main results.
\begin{thm}[Informal statement of main results]
	\label{theorem::introduction-main}
	Assume the network is generated from an exchangeable network model.  Define the Edgeworth expansion for a given network moment ${\cal R}$ with $r$ nodes $s$ edges as follows:
	\begin{align*}
		G_n(x) &:= \Phi(x) + \frac{\phi(x)}{\sqrt{n}\cdot \xi_1^3} \cdot \Bigg\{ \frac{ 2 x^2 + 1}{6} \cdot\ep[g_1^3(X_1)]\notag\\
		& + \frac{r-1}2\cdot \left(  x^2 + 1 \right)\ep[g_1(X_1)g_1(X_2)g_2(X_1,X_2)]\Bigg\},
	\end{align*}
	where $\Phi, \phi$ are the CDF and PDF of $N(0,1)$, and the {estimable coefficients components $\xi_1$, $\ep[g_1^3(X_1)]$ and $\ep[g_1(X_1)g_1(X_2)g_2(X_1,X_2)]$,} to be defined in Section \ref{section::moments-and-Edgeworth-expansion}, {only depend} on the graphon $f$ and the motif $R$.  Let $\rho_n$ denote the network sparsity parameter.  {In the dense regime, where we assume:}
	\begin{enumerate}
		\item $\rho_n^{-2s}\cdot \var(g_1(X_1))\geq\textrm{constant}>0$;
		\item $\rho_n = \omega(n^{-1/2})$ for acyclic $R$, or $\rho_n = \omega(n^{-1/r})$ for cyclic $R$;
		\item Either $\rho_n\preceq (\log n)^{-1}$, or $\limsup_{t\to\infty}\left|\ep\left[e^{\ii t g_1(X_1)/\xi_1}\right]\right|<1$;
	\end{enumerate}
	we have
	\begin{equation}
		\left\| F_{\hat{T}_n}(u) - G_n(u) \right\|_\infty = \Ohighorderbound,
		\label{thmeqn::introduction-main}
	\end{equation}
	where $\|H(u)\|_\infty := \sup_{u\in\mathbb{R}}\left| H(u) \right|$, and ${\cal M}(\rho_n,n;R)$, defined in \eqref{M}, satisfies ${\cal M}(\rho_n,n;R)\ll n^{-1/2}$.  Under the same conditions, the empirical Edgeworth expansion $\hat{G}_n$ with estimated coefficients (see \eqref{Edgeworth::empirical}) satisfies
	{
	\begin{equation}
	\pr\left(\left\| F_{\hat{T}_n}(u) - \hat{G}_n(u) \right\|_\infty > C\cdot {\cal M}(\rho_n,n;R)\right) = O(n^{-1}).
	\label{thmeqn::introduction-main-empirical}
	\end{equation}
	for a large enough absolute constant $C$.
	
	In the sparse regime, where we replace condition 2 by: 
	\begin{enumerate}
	    \item[2'.] $n^{-1}\prec \rho_n\preceq n^{-1/2}$ for acyclic $R$, or $n^{-2/r}\prec \rho_n\preceq n^{-1/r}$ for cyclic $R$,
	\end{enumerate}
	a simple $N(0,1)$ approximation achieves the following Berry-Esseen bound:
	\begin{align}
	    \label{thmeqn::introduction-main-sparse}
        \left\|
            F_{\hat T_n}(u) - G_n(u)
        \right\|_\infty
        &
        \asymp
        \left\|
            F_{\hat T_n}(u) - \Phi(u)
        \right\|_\infty =
        \Ohighorderbound \bigwedge o(1).\notag
    \end{align}
    Moreover, we have
    $$
        \left\|
            F_{\hat T_n}(u) - \hat G_n(u)
        \right\|_\infty
        =
        \begin{cases}
            \tildehighorderbound, &\textrm{ if }\rho_n=\omega(n^{-1}\log^{1/2} n),\\
            o_p(1), &\textrm{ if }\rho_n = \omega(n^{-1})
        \end{cases}
    $$
    That is, in the sparse regime, the empirical Edgeworth expansion has the same order of approximation error as $N(0,1)$.}
\end{thm}

\subsection{Our contributions}

Our contributions are three-fold.
First, we establish the first accurate distribution approximations for network moments \eqref{thmeqn::introduction-main}.  {The results} originated from our {discovery of the surprisingly blessing} roles that network noise and sparsity jointly play in this setting.  {Our work reveals a new dimension to the understanding of these two components in network analysis.}  
Second, we propose a provably highly accurate and computationally efficient empirical Edgeworth approximation \eqref{thmeqn::introduction-main-empirical} for practical use.
Third, our results {enable} accurate and fast nonparametric network inference procedures.

To understand the strength of our main results \eqref{thmeqn::introduction-main} and \eqref{thmeqn::introduction-main-empirical}, notice that for mildly sparse networks, we achieve \emph{higher-order accuracy} in distribution approximation \emph{without non-lattice or smoothness assumption}.  To our best knowledge, the non-lattice assumption is universally required to achieve higher-order accuracy in all literature {for similar settings}.  However, this assumption is violated by some popular network models such as stochastic block model, arguably {one of the most widely-used} network models.  Waiving the graphon smoothness assumption makes our approach a powerful tool for model-free exploratory network analysis and for analyzing networks with {high complexity and} irregularities.

{
In the sparse regime, our modified Berry-Esseen bound \eqref{thmeqn::introduction-main-sparse} significantly improves over the previous best known bound $o(1)$ in existing literature \citep{bickel2011method, bhattacharyya2015subsampling, green2017bootstrapping, levin2019bootstrapping} and fills a large missing part in the big picture.  As the network sparsity $\rho_n$ declines from $n^{-1/2}$ towards $n^{-1}$ for acyclic $R$ and from $n^{-1/r}$ towards $n^{-2/r}$ for cyclic $R$, our result reveals a gradually depreciating uniform error bound.  When $\rho_n$ hits the minimum assumption boundary, our result matches the uniform consistency result in classical literature.
}

Our key insight is {to view} the sample network moment $\hat{U}_n$ as a \emph{noisy U-statistic}, where ``noise'' refers to edge-wise observational errors in $A$.  Our analysis reveals the connection and differences between the noisy and the conventional \emph{noiseless} U-statistic settings.  We discover the surprisingly blessing roles that the two typically-hated factors, namely, \emph{edge-wise observational errors} and \emph{network sparsity} jointly play in this setting:
\begin{enumerate}
	\item The edge-wise errors behave like a smoother that tames potential distribution discontinuity due to a lattice or discrete network population\footnote{More precisely speaking, such irregularity is jointly induced by both the network population distribution and the shape of the motif, but the former is usually the determining factor.};
	\item Network sparsity boosts the smoothing effect of the error term to a sufficient level such that $F_{\hat{T}_n}$ becomes analytically tractable.
\end{enumerate}
{At first sight, the smoothing effect of edge-wise errors is rather counter-intuitive, since generating $A$ from $W$ \emph{discretizes} edge probabilities from numbers originally in a continuum $[0,1]$ to binary entries.  How could this eventually yield a smoothing effect?  In Section \ref{subsec::U-hat}, we present two simple examples to illustrate the intuitive reason.}  In our proofs, we present original analysis to carefully quantify the impact of such smoothing effect.  Our analysis techniques are very different from those in network bootstrap papers \cite{bhattacharyya2015subsampling, green2017bootstrapping, levin2019bootstrapping, lunde2019subsampling}.  It seems unlikely that our assumptions can be substantially relaxed since they match the well-known minimum conditions in related settings in \cite{lahiri1993bootstrapping}.

Our empirical Edgeworth expansion \eqref{thmeqn::introduction-main-empirical} is very fast, much more scalable than network bootstraps, and easily permits parallel computing.

We showcase three applications of our theory.  We present the first proof of the higher-order accuracy of some mainstream network bootstrap techniques under certain conditions, which their original proposing papers did not prove.  Our results also enable rich future works on accurate and {computationally} very efficient network inferences.  We present two immediate applications to testing and Cornish-Fisher type confidence interval for network moments with explicit accuracy guarantees.

\subsection{Paper organization}

The rest of this paper is organized as follows.  In Section \ref{sec::problem-set-up-and-literature-review}, we formally set up the problem and provide a detailed literature review.  In Section \ref{section::moments-and-Edgeworth-expansion}, we present our core ideas, derive the Edgeworth expansions and establish their uniform approximation error bounds.  We discuss different versions of the studentization form.  {We also present our modified Berry-Esseen theorem for the sparse regime.}  In Section \ref{sec::applications}, we present three applications of our results:  bootstrap accuracy, one-sample test, and one-sample Cornish-Fisher confidence interval.  In Section \ref{sec::simulations}, {we conduct three simulations to evaluate the performance of our method from various aspects}.  Section \ref{sec::discussion} discusses interesting implications of our results and future work.  

{
\subsection{Big-O and small-o notation system}
\label{subsec::tOp}
In this paper, we will make frequent references to the big-O and small-o notation system.  We use the same definitions of $O(\cdot)$, $o(\cdot)$, $\Omega(\cdot)$ and $\omega(\cdot)$ as that in standard mathematical analysis, and the same $O_p(\cdot)$ and $o_p(\cdot)$ as that in probability theory. 
For a random variable $Z$ and a deterministic sequence $\{\alpha_n\}$, define $\tOp(\cdot)$ as follows
\begin{align}
    Z := \tOp(\alpha_n), & \textrm{ if } \pr(|Z|\geq C\alpha_n) = O(n^{-1}) \textrm{ for some constant $C>0$.}
\end{align}
This is similar to ``o$_p$'' in \cite{maesono1997edgeworth} (see the remark beneath Lemma 2) and Assumption (A1) in \cite{lai1993edgeworth}.  For technical reasons, in this paper, we do not need to define a $\top(\cdot)$ sign.

}

\section{Problem set up and literature review}
\label{sec::problem-set-up-and-literature-review}
\subsection{Exchangeable networks and graphon model}
\label{subsec::problem-setup-graphon-model}

The base model of this paper is exchangeable network model \citep{diaconis2007graph, bickel2009nonparametric}.  Exchangeability describes the unlabeled nature of many networks in social, knowledge and biological contexts, where node indices do not carry meaningful information.  It is a very rich family that contains many popular models as special cases, including the stochastic block model and its variants \footnote{{Here we follow the convention of \cite{bickel2009nonparametric} and view community memberships as randomly sampled from a multinomial distribution.}} \citep{holland1983stochastic, zhao2012consistency, zhang2016minimax, airoldi2008mixed, karrer2011stochastic, zhang2014detecting, jing2020community}, the configuration model \citep{chung2002connected, newman2009random}, latent space models \citep{hoff2002latent, grover2016node2vec} and general smooth graphon models \citep{choi2012stochastic, gao2015rate, zhang2017estimating}\footnote{{Smooth graphon: we can simply think that a graphon is called ``smooth'' if $f(\cdot,\cdot)$ is a smooth function.  In the rigorous definition, $f$ is smooth if $f(\psi(\cdot),\psi(\cdot))$ is smooth under some measure-preserving map $\psi:[0,1]\to[0,1]$, see \cite{bickel2009nonparametric, gao2015rate, zhang2017estimating}.}}. 
{In this paper, we base our study on the following exchangeable network model called \emph{graphon model}.  The framework is closely related to the Aldous-Hoover representation  for infinite matrices \citep{aldous1981representations, hoover1979relations}.  Under a graphon model,}  the $n$ nodes correspond to latent space positions $X_1,\ldots,X_n\stackrel{\textrm{i.i.d.}}{\sim}$~Uniform$[0,1]$.  Network generation is governed by a measurable latent graphon function $f(\cdot,\cdot): [0,1]^2\to [0,1]$, $f(x,y)=f(y,x)$ that encodes all structures.  The edge probability between nodes $(i,j)$ is
\begin{equation}
	W_{ij} = W_{ji} := \rho_n\cdot f(X_i,X_j); \quad 1\leq i<j\leq n,
	\label{graphon-model::W}
\end{equation}
where the sparsity parameter $\rho_n\in(0,1)$ absorbs the constant factor, and we fix {$\int_{[0,1]^2}f(u,v)\td u\td v=$~constant}.  We only observe the adjacency matrix $A$:
\begin{equation}
	A_{ij}=A_{ji}|W \sim \textrm{Bernoulli}(W_{ij}), \forall 1\leq i<j\leq n.
	\label{graphon-model::A|W}
\end{equation}
The model defined by \eqref{graphon-model::W} and \eqref{graphon-model::A|W} has a well-known issue that both $f$ and $\{X_1,\ldots,X_n\}$ are only identifiable up to equivalence classes \citep{chan2014consistent}.  This may pose significant challenges for model-based network inference, {especially those based on parameter estimations}.  On the other hand, network moments are permutation-invariant and thus clearly immune to this identification issue{, making them attractive candidates for inference objective}.

\subsection{Network moment statistics}
\label{subsec::problem-setup-network-moments}


To formalize network moments, it is more convenient to first define the sample version and then the population version.  Each network moment is indexed by the corresponding motif ${\cal R}$.  For simplicity, we follow the convention to focus on connected motifs.  Let $R$ represent the adjacency matrix of ${\cal R}$ with $r$ nodes and $s$ edges.  For any $r$-node sub-network $A_{i_1,\ldots,i_r}$ of $A$, define
\begin{equation}
	h(A_{i_1,\ldots,i_r}) := \mathbbm{1}_{[A_{i_1,\ldots,i_r}{\sqsupseteq} R]}\footnote{Since we consider an arbitrary but fixed $R$ throughout this paper, without causing confusion, we drop the dependency on $R$ in symbols such as $h$ to simplify notation.},\quad \textrm{ for all }1\leq i_1<\cdots<i_r\leq n,
	\label{def::h(A)}
\end{equation}
Here, ``$A_{i_1,\ldots,i_r}{{\sqsupseteq}} R$'' means there exists a permutation map $\pi:\{1,\ldots,r\}\to\{1,\ldots,r\}$, such that $A_{i_1,\ldots,i_r} {\geq} R_\pi$, {where the ``$\geq$'' is entry-wise} and $R_\pi$ is defined as $(R_\pi)_{ij}:= R_{\pi(i)\pi(j)}$.
{
Our definition of $h(A_{i_1,\ldots,i_r})$ here is similar to the ``$Q(R)$'' defined in \cite{bickel2011method}.  One can similarly define
\begin{equation}
    \tilde h(A_{i_1,\ldots,i_r}) := \mathbbm{1}_{[A_{i_1,\ldots,i_r}\cong R]},\quad \textrm{ for all }1\leq i_1<\cdots<i_r\leq n,
	\label{def::htilde(A)}
\end{equation}
where ``$A_{i_1,\ldots,i_r}\cong R$'' means there exists a permutation map $\pi:\{1,\ldots,r\}\to\{1,\ldots,r\}$, such that $A_{i_1,\ldots,i_r} = R_\pi$.
The definition of $\tilde h$ corresponds to the ``$P(R)$'' studied in \cite{bickel2011method, bhattacharyya2015subsampling}, and \cite{green2017bootstrapping}.  As noted by \cite{bickel2011method}, each $h$ can be explicitly expressed as a linear combination of $\tilde h$ terms, and vice versa.  Therefore, they are usually treated with conceptual equivalence in literature, and most existing papers would choose one of them to study.}
{
For technical cleanness, in this paper we focus on $h$.  We believe our method is also applicable to analyzing $\tilde h$, but the analysis is much more complicated and we leave it to future work.
}
Define the \emph{sample network moment} as
\begin{equation}
\hat{U}_n := \dfrac1{\binom{n}r}\sum_{1\leq i_1<\cdots<i_r\leq n} h(A_{i_1,\ldots,i_r}),
\label{def::sample-moments}
\end{equation}
Then we define the \emph{sample-population version} and \emph{population version} of $\hat{U}_n$ to be $U_n := \ep[ \hat{U}_n|W ]$ and $\mu_n := \ep[U_n] = \ep[\hat{U}_n ]$, respectively.  We refer to $\hat{U}_n$ as the \emph{noisy} U-statistic, and call $U_n := \binom{n}r^{-1}\sum_{1\leq i_1<\cdots<i_r\leq n} h(W_{i_1,\ldots,i_r}) = \binom{n}r^{-1}\sum_{1\leq i_1<\cdots<i_r\leq n} h(X_{i_1},\ldots,X_{i_r})$\footnotemark\label{footnote::def::h} the (conventional) \emph{noiseless} U-statistic.  Similar to the {insight that studentization is key to achieve higher-order accurate approximations} in the i.i.d. setting (Section 3.5 of \cite{wasserman2006all}), we study
\begin{equation}
\hat{T}_n := \frac{\hat{U}_n - \mu_n}{\hat{S}_n},\notag
\end{equation}
where $\hat{S}_n$ will be specified later {in \eqref{def::jackknife-variance-estimator} and 
\eqref{def::our-variance-estimator}}.  We can similarly {standardize} or studentize the noiseless U-statistic $U_n$ by
$\check{T}_n:=(U_n-\mu_n)/\sigma_n$ and
$T_n:=(U_n-\mu_n)/S_n$, respectively, where $\sigma_n^2:=\var(U_n)$ and $S_n^2$ is a {$\sqrt{n}$-consistent estimator\footnotemark\label{footnote::def::sqrt-n-consistency} for $\sigma_n^2$, which will be specified later}.

\footnotetext[\getrefnumber{footnote::def::h}]{{Here, without causing confusion, we slightly abused the notation of $h(\cdot)$, letting it take either $W$ or $X$ as its argument, noticing that $W$ is determined by $X_1,\ldots,X_n$.}}

\footnotetext[\getrefnumber{footnote::def::sqrt-n-consistency}]{{$\sqrt{n}$-consistency of $S_n^2$ means that $\sqrt{n}(S_n^2-\sigma_n^2)=o_p(1)$, see \cite{bhattacharyya2015subsampling, levin2019bootstrapping} for definition.}}

\subsection{Edgeworth expansions for i.i.d. data and noiseless U-statistics}

Edgeworth expansion \citep{edgeworth1905law, wallace1958asymptotic} refines the central limit theorem.  It is the supporting pillar in the justification of bootstrap's higher-order accuracy.  In this subsection, we review the literature on Edgeworth expansions for i.i.d. data and {conventional noiseless} U-statistics, due to their close connection.  Under mild conditions, the one-term Edgeworth expansion for {the sample mean of $n$} i.i.d. mean-zero and unit-variance $X_1,\ldots,X_n$ {reads} $F_{n^{{1/2}}(\bar{X} - \ep[X_1])/\sigma_X}(u) = \Phi(u) -  n^{-1/2}\cdot \ep[X_1^3](u^2-1)\phi(u)/6 + O(n^{-1})$, where $\Phi$ and $\phi$ are the CDF and PDF of $N(0,1)$, respectively.  Edgeworth terms of even higher orders can be derived \citep{hall2013bootstrap} but are not meaningful in practice unless we know a few true population moments.  The minimax rate for estimating $\ep[X_1^3]$ is $O_p(n^{-1/2})$, so $O(n^{-1})$ is the best practical remainder bound for an Edgeworth expansion.  For further references, see \cite{bickel1974edgeworth, sakov2000edgeworth, bhattacharya1978validity, hall1987edgeworth, hall1993edgeworth, babu1984one} and textbooks \cite{hall2013bootstrap, davison1997bootstrap, wasserman2006all}.

The literature on Edgeworth expansions for U-statistics concentrates on the noiseless version.  In early 1980's, \cite{callaert1978berry, janssen1981rate, callaert1981order} established the asymptotic normality of {the standarized} and {the studentized U-statistics, respectively, both} with $O(n^{-1/2})$ Berry-Esseen type bounds.
Then \cite{callaert1980edgeworth, bickel1986edgeworth,lai1993edgeworth} approximated degree-two (i.e. $r=2$) standardized U-statistics with an $o(n^{-1})$ remainder {with known population moments}, and \cite{bentkus1997edgeworth} established an $O(n^{-1})$ bound under relaxed conditions for more general symmetric statistics.
Later, \cite{helmers1991edgeworth, putter1998empirical} studied empirical Edgeworth expansions (EEE) {with estimated coefficients} and established $o(n^{-1/2})$ bounds.  For finite populations, \cite{babu1985edgeworth, bloznelis2001orthogonal, bloznelis2002edgeworth, bloznelis2003edgeworth} established the earliest results, and we will use some of their results in our analysis of network bootstraps.  An incomplete list of other notable works {on Edgeworth expansions for noiseless U-statistics with various finite moment assumptions} includes \cite{bentkus1994lower, hall1995uniform, jing2003edgeworth, maesono1997edgeworth, bentkus2009normal, jing2010unified}.

\subsection{The non-lattice condition and lattice Edgeworth expansions in the i.i.d. setting}
\label{subsec::Edgeworth-lattice-iid}

A major assumption called the \emph{non-lattice condition}  is critical for achieving $o(n^{-1/2})$ accuracy in Edgeworth expansions {and is needed by} all results in the i.i.d. setting {without} oracle moment knowledge and all results for noiseless U-statistics, but this condition is clearly not required {for} an $O(n^{-1/2})$ accuracy bound\footnote{Simply use a Berry-Esseen theorem.}.  A random variable $X_1$ is called \emph{lattice}, if it is supported on $\{a + bk: k\in\mathbb{Z}\}$ for some $a,b\in\mathbb{R}$ where $b\neq 0$.  General discrete distributions are ``nearly lattice''
\footnote{``A discrete distribution is nearly-lattice'': a discrete distribution, if not already lattice, can be viewed as a lattice distribution with diminishing periodicity.}.  
A distribution is essentially \emph{non-lattice} if it contains a continuous component.  In many works, the non-lattice condition is replaced by the stronger Cramer's condition \citep{cramer1928composition}:
\begin{equation}
\limsup_{t\to\infty}\left| \ep\left[ e^{\ii t X_1} \right] \right| <1.\notag
\end{equation}
For U-statistics, this condition is imposed on $g_1(X_1):=\ep[h(X_1,\ldots,X_r)|X_1] - \mu_n$.  Cramer's condition can be relaxed \citep{angst2017weak, mattner2017optimal, song2016ordering, song2018uniform} towards a non-lattice condition, but all {existing} relaxations come at the price of essentially depreciated error bounds
\footnote{{To our knowledge,} existing {results} assuming only non-lattice{ness} achieve no better than $o(n^{-1/2})$ approximation errors.  {For example,} \cite{bentkus1997edgeworth} replaces {the RHS} ``1'' in Cramer's condition by $1-q$ {and assumes} it holds for $t\preceq n^{1/2}$.  {They} obtain an error bound proportional to $q^{-2}$.  Another example is \cite{bloznelis2002edgeworth}.  It replaces \cite{bentkus1997edgeworth}'s $t$ range by $t\preceq \pi$ {(their $\pi$ is a variable)} and obtains an error bound proportional to $q^{-2}\pi^{-2}$.  Also see the comment {beneath} equation (4.7) of \cite{putter1998empirical}.}.  
Therefore, for simplicity, in Theorems \ref{thm::main-theorem} and \ref{thm::bootstrap-accuracy}, we use Cramer's condition to represent the non-lattice setting.

However, in network analysis, Cramer's condition {may be a strong assumption, for the following reasons.  First, it} is violated by stochastic block model, {a very popular and important} network model.  In a block model, $g_1(X_1)$ only depends on node $1$'s community membership, thus is discrete.  
{Second, this} condition is difficult to check in practice.  
{Third, }some smooth models may even {induce a lattice $g_1(X_1)$ under certain} motifs {and a non-lattice $g_1(X_1)$ under a different motif}.  For example, under the graphon model $f(x,y):= 0.3+0.1\cdot \mathbbm{1}_{[x>1/2; y>1/2]} + 0.1 \sin\left( 2\pi(x+y) \right)$,
{$g_1(X_1)$ is lattice} when $R$ is an edge, but it is non-lattice when $R$ is a triangle.

Next, we {review existing treatments of Edgeworth expansion in the lattice case that will spark} the key inspiration to our work.  
{In current literature, in the lattice case, we could approximate the CDF of an i.i.d. sample mean at higher-order accuracy, where the lattice Edgeworth expansion would contain an order $n^{-1/2}$ jump function; whereas to our best knowledge, no analogous result exists for U-statistics.}
{Available} approaches can be categorized into two mainstreams: (1) adding an artificial error term to the sample mean to smooth out lattice-induced discontinuity \citep{singh1981asymptotic, lahiri1993bootstrapping}; and (2) formulating the lattice version Edgeworth expansion with a jump function \citep{singh1981asymptotic}.  The seminal work \cite{singh1981asymptotic} adds a uniform error of bandwidth $n^{-1/2}$, and by {inverting} its impact {on} the smoothed distribution function, it {explicitly} formulates the lattice Edgeworth expansion with an $O(n^{-1})$ remainder.  Another classical work \cite{lahiri1993bootstrapping} uses a normal {artificial} error instead of uniform and shows that the Gaussian bandwidth must be $\omega((\log n/ n)^{1/2})$ and $o(1)$ to provide sufficient smoothing effect without {causing} an $\omega(n^{-1/2})$ distribution distortion.  Other notable works include \cite{woodroofe1988singh, kolassa1990edgeworth, babu1989edgeworth},
{in which, \cite{woodroofe1988singh} and \cite{kolassa1990edgeworth} also formulate lattice Edgeworth expansions in the i.i.d. univariate setting, and \cite{babu1989edgeworth} studies Edgeworth expansions for the sample mean of i.i.d. random vectors, where some dimensions are lattice and the others are non-lattice.}

{Despite the significant achievements of these treatments, latticeness remains an obstacle in practice.  The difficulties are two-fold.}
{On one hand, if we introduce an artificial error to smooth the distribution, it will unavoidably} bring an $\Omega(n^{-1/2})$ distortion to the original distribution\footnote{To see this, simply notice that the original distribution contains $n^{-1/2}$ jumps, but the smoothed distribution does not, {so an $o(n^{-1/2})$ approximation error is impossible} \citep{bickel1986edgeworth}.}.  
{On the other hand, the exact formulation of a lattice Edgeworth expansion contains} an $n^{-1/2}$ jump term.  In many examples such as bootstrap, the jump locations depend on the true population variance, laying an uncrossable $\Omega(n^{-1/2})$ barrier for practical CDF approximation.  {For more details, see page 91 of \cite{hall2013bootstrap}.}

\section{Edgeworth expansions for network moments}
\label{section::moments-and-Edgeworth-expansion}

\subsection{Outline and core ideas to analyze $\hat{T}_n$}
\label{subsec::Edgeworth-outline}

Our key discovery is that the studentized noisy U-statistic $\hat{T}_n$ can be decomposed as follows:
\begin{equation}
\hat{T}_n = \tilde{T}_n + {\check{\Delta}_n + \tildehighorderbound},
\label{eqn::T-hat-main-decomposition}
\end{equation}
where $\tilde{T}_n$, {to be formally defined in \eqref{def::T-tilde},} can be roughly understood as a studentized noiseless U-statistic {who could be approximately decomposed into an $O_p(1)$ linear part and an $O_p(n^{-1/2})$ quadratic part}, 
${\check{\Delta}_n} \approx N(0,\sigma^2\asymp (\rho_n\cdot n)^{-1})$, and {recall the symbol $\tOp$ from Section \ref{subsec::tOp}}.
{Here the remainder term would be ignorable if the network is mildly dense.  In the sparse regime, the remainder will dominate both ${\check\Delta_n}$ and the quadratic part of $\tilde T_n$.}

Our decomposition \eqref{eqn::T-hat-main-decomposition} is a renaissance of the spirits of \cite{singh1981asymptotic} and \cite{lahiri1993bootstrapping}, but with the following crucial {conceptual} differences.  First and most important, the error term ${\check\Delta_n}$ in our formula is \emph{not} artificial, but naturally a constituting component of $\hat{T}_n$.  Therefore, the smoother does \emph{not} distort the objective distribution, that is, $\hat{T}_n$ is \emph{self-smoothed}.  The second difference lies in the bandwidth of the smoothing error term. 
{Since the smoothing error terms in \cite{singh1981asymptotic} and \cite{lahiri1993bootstrapping} are artificial, the user is at the freedom to choose these bandwidths.  In our setting, the bandwidth of the smoothing term $(\rho_n\cdot n)^{-1/2}$ is not managed by the user, but governed by the network sparsity.  Therefore, when Cramer's condition fails, we make the very mild sparsity assumption that $\rho_n=O((\log n)^{-1})$ to ensure enough smoothing effect.
}

This echoes the lower bound on Gausssian bandwidth in \cite{lahiri1993bootstrapping}.  We also need $\rho_n$ to be lower bounded {to effectively bound the remainder term}, see Lemma \ref{lemma::term-approx}{-\eqref{lemma::term-approx-Delta-hat-conditional-normal}}.  Third, our error term ${\check\Delta_n}$ is \emph{dependent} on $\tilde{T}_n$ through $W$. 
{Last, the proof technique of \cite{singh1981asymptotic} is inapplicable to our setting due to the quadratic part in $\tilde T_n$; and \cite{lahiri1993bootstrapping} obtains an $o(n^{-1/2})$ error bound\footnote{{The $o(n^{-1/2})$ error bound in \cite{lahiri1993bootstrapping} holds on some ${\frak B}\subset \mathbb R$ with ``diminishing boundary'', while our error bounds hold on the entire $\mathbb{R}$.}}, while we aim at stronger results under a more complicated U-statistic setting with degree-two terms.}
In our proofs, we carefully handle {these challenges} with original analysis.


\subsection{Decomposition of the stochastic variations of $\hat{U}_n$}
\label{subsec::U-hat}

To simplify narration, here we focus on analyzing $\hat{U}_n$, and the analysis of $\hat{T}_n$ is conceptually similar.  The stochastic variations in $\hat{U}_n=U_n+(\hat{U}_n-U_n)$ {come} from two sources: the randomness in $U_n$ due to $W$ and ultimately $X_1,\ldots,X_n$, and the randomness in $\hat{U}_n-U_n$ due to $A|W$, the edge-wise observational errors.  The stochastic variations in $U_n$ as a conventional noiseless U-statistic is well-understood due to Hoeffding's decomposition \citep{hoeffding1948class}:
\begin{align}
U_n -\mu_n
&=\frac{r}{n} \sum_{i=1}^n g_1(X_i) + \frac{r(r-1)}{n(n-1)}\sum_{1\leq i<j\leq n} g_2(X_i,X_j) + {\tOp(\rho_n^s\cdot n^{-3/2}\log^{3/2}n)}
\label{eqn::Hoeffding's decomposition}
\end{align}
where $g_1,\ldots,g_r$ are defined as follows.  To avoid complicated subscripts, without confusion we define $g_k$'s for special indexes $(i_1,\ldots,i_r)=(1,\ldots,r)$.  For indexes $1$, $k=\{2,\ldots,r-1\}$ (only when $r\geq 3$) and $r$, define $g_1(x_1) := \ep[h(X_1,\ldots,X_r)|X_1=x_1] - \mu_n$, $g_k(x_1,\ldots,x_k) := \ep[h(X_1,\ldots,X_r)|X_1=x_1,\ldots,X_{{k}}=x_{{k}}] - \mu_n - \sum_{k'=1}^{k-1} \sum_{1\leq i_1<\ldots <i_{k'}\leq r} g_{k'}(x_{i_1},\ldots,x_{i_{k'}})$ for $2\leq k\leq r-1$ and $g_r(x_1,\ldots,x_r) := h(x_1,\ldots,x_r) - \mu_n$.  From classical literature, we know that $\ep[g_k(X_{i_1},\ldots,X_{i_k})|\{X_i:i\in{\cal I}_k\subset\{i_1,\ldots,i_k\}\}] = 0$, where the strict subset ${\cal I}_k$ could be $\emptyset$, and $\cov\left( g_k(X_{i_1},\ldots,X_{i_k}), g_\ell(X_{j_1},\ldots,X_{j_\ell}) \right)=0$ unless $k=\ell$ and $\{i_1,\ldots,i_k\} = \{j_1,\ldots,j_\ell\}$.  Consequently, the linear part in the Hoeffding's decomposition is dominant.  Define
\begin{equation}
\xi_1^2 := \var(g_1(X_1)).
\label{def::xi_1}
\end{equation}
We focus on discussing the stochastic variations in $\hat{U}_n-U_n$.  The typical treatment in network bootstrap literature is to simply bound and ignore this component, such as Lemma 7 in \cite{green2017bootstrapping}.  {In sharp contrast, by carefully quantifying its impact,} we shall reveal its key smoothing effect by a refined analysis.  To better understand the impact of $\hat{U}_n-U_n$, let us inspect two simple {and illustrative} examples.  {In these examples, we inspect $\hat T_n$ and use the fact that $\hat S_n\asymp \sigma_n\asymp \rho_n^s\cdot n^{-1/2}$ (by Lemma~\ref{lemma::term-approx})}.
\begin{example}
	\label{example::motif-edge}
	Let $R$ be an edge with $r=2$ and $s=1$, and $\hat{U}_n$ is simply the sample edge density.  By definition, all $h(A_{i_1,i_2})-h(W_{i_1,i_2})$ terms are mutually independent given $W$.  Then {the asymptotic behavior of the self-smoother term is}
	$$
	\frac{\hat{U}_n - U_n}{{\hat S_n}} 
	\stackrel{d}{\to}
	N\left(0,  
	\sigma_{\frac{\hat{U}_n-U_n}{{\hat S_n}}\big|W}^2
	{\asymp(\rho_n\cdot n)^{-1}}
	\right)
	$$ 
	{at} a uniform $O(\rho_n^{-1/2}\cdot n^{-1})$ Berry-Esseen CDF approximation error {rate}.
\end{example}
The next example shows that the insight of Example \ref{example::motif-edge} generalizes.
\begin{example}
	\label{example::motif-triangle}
	Let $R$ be a triangular motif with $r=3, s=3$, and $\hat{U}_n$ is the empirical triangle frequency.  We can decompose $\hat{U}_n - U_n$ as follows:
	\begin{align}
	&\frac{\hat{U}_n - U_n}{{\hat S_n}}
	= 
	\frac1{\binom{n}3} \sum_{1\leq i_1<i_2<i_3\leq n}
	\frac{\left\{h(A_{i_1,i_2,i_3})-h(W_{i_1,i_2,i_3})\right\}}{{\hat S_n}} \notag\\
	=&
	\frac1{\binom{n}3} \sum_{1\leq i_1<i_2<i_3\leq n}
	\frac{  \left( W_{i_1i_2} + \eta_{i_1i_2} \right)\left( W_{i_1i_3} + \eta_{i_1i_3} \right)\left( W_{i_2i_3} + \eta_{i_2i_3} \right) - W_{i_1i_2}W_{i_1i_3}W_{i_2i_3}  }{{\hat S_n}}\notag\\
	=&
	{\frac1{\binom{n}3}\Bigg\{ \sum_{\substack{1\leq i_1<i_2\leq n\\1\leq i_3\leq n\\i_3\neq i_1,i_2}} \frac{ W_{i_1i_3}W_{i_2i_3} \eta_{i_1i_2} +   W_{i_1i_2}\eta_{i_1i_3}\eta_{i_2i_3}}{\hat S_n} + \sum_{1\leq i_1<i_2<i_3\leq n} \frac{\eta_{i_1i_2}\eta_{i_1i_3}\eta_{i_2i_3}}{\hat S_n} \Bigg\}}\notag\\
	=&
	{ 
	\underbrace{\frac1{\binom{n}2}\sum_{1\leq i<j\leq n}\left(\frac{3\sum_{\substack{1\leq k\leq n\\k\neq i,j}}W_{ik}W_{jk}}{(n-2)\hat S_n}\right)\eta_{ij} }_{\textrm{Linear part}}
	}
	+
	{
	\underbrace{
	\frac1{\binom{n}3}\sum_{\substack{1\leq i<j\leq n\\1\leq k\leq n\\k\neq i,j}} \frac{W_{ij}}{\hat S_n} \eta_{ik}\eta_{jk}
	}_{\textrm{Quadratic part}}
	}\notag\\
	&+
	{
	\underbrace{
	\frac1{\binom{n}3}\sum_{1\leq i<j<k\leq n} \frac1{\hat S_n}\eta_{ij}\eta_{ik}\eta_{jk}
	}_{\textrm{Cubic part}}
	}\notag
	\end{align}
	where $\eta_{ij}:= A_{ij}-W_{ij}$.
    {
    The linear part is $\asymp \rho_n^{-1/2}\cdot n^{-1/2}$, the quadratic part is $\tOp(\rho_n^{-1}\cdot n^{-1}\log^{1/2}n)$ and the cubic part is $\tOp(\rho_n^{-3/2}\cdot n^{-1}\log^{1/2}n)$.  We make two observations.  First, the linear part in this example has the same asymptotic order as the linear part in Example \ref{example::motif-edge}.  This is not a coincidence and will be formalized by Lemma \ref{lemma::term-approx}-\eqref{lemma::term-approx-Delta-hat-conditional-normal}.  In other words, regardless of the shape of $R$, the linear part in such decomposition always provides smoothing effect at the same magnitude.  Second, different from Example \ref{example::motif-edge}, we now have higher-degree remainder terms.  The linear part nicely always dominates the quadratic part; but it only dominates the cubic part when $\rho_n=\omega(n^{-1/2}\log^{1/2}n)$.
    }
\end{example}

The insights of the two examples are generalized in Lemma \eqref{lemma::term-approx}-\eqref{lemma::term-approx-Delta-hat-conditional-normal}.  When the network is moderately dense, the linear part in $\hat{U}_n-U_n$ dominates.  Consequently, the overall contribution of the stochastic variations in $A|W$ approximates Gaussian {at} an $O(\rho_n^{-1/2}\cdot n^{-1})$ Berry-Esseen {error rate}.

\subsection{Studentization form}
\label{subsec::our-method-studentization-form}

The understanding of $\hat{U}_n$ in Section \ref{subsec::U-hat} prepares us to fully specify $\hat{T}_n = (\hat{U}_n - \mu_n)/\hat{S}_n$.  {Recall that we still need to} design $\hat{S}_n$.  In $\var(\hat{U}_n) = \ep[\var(\hat{U}_n|W)] + \var(\ep[\hat{U}_n|W])$, we observe $\var(\hat{U}_n|W)\asymp \rho_n^{2s-1}\cdot n^{-2}$ and $\var(\ep[\hat{U}_n|W]) = \var(U_n) \asymp \rho_n^{2s} \cdot n^{-1}$.  We shall assume $\rho_n \cdot n\to \infty$, so $\sigma_n^2=\var(U_n) = \var(\ep[\hat{U}_n|W])$ dominates.  There are two main choices of $\hat{S}_n$.  The conventional choice for studentizing noiseless U-statistics \citep{callaert1981order, helmers1991edgeworth, putter1998empirical} {uses} the jackknife estimator
\begin{equation}
n\cdot \hat{S}_{n;\textrm{jackknife}}^2 := (n-1)\sum_{i=1}^n\left( \hat{U}_n^{(-i)} - \hat{U}_n \right)^2,
\label{def::jackknife-variance-estimator}
\end{equation}
where $\hat{U}_n^{(-i)}$ is $\hat{U}_n$ calculated on the induced sub-network of $A$ with node $i$ removed.  Despite conceptual straightforwardness, the jackknife estimator unnecessarily complicates analysis.  Therefore, we use an estimator with a simpler formulation.  In $\var(\hat{U}_n) = \sigma_n^2 + O(\rho_n^{2s-1}n^{-2}) =  r^2\xi_1^2/n + O(\rho_n^{2s-1}n^{-2})$, replace $\xi_1$ by its moment estimator.  {Specifically, recall that $\xi_1^2 = \var(g_1(X_1)) = \ep[(\ep[h(X_1,\ldots,X_n)|X_1]-\mu_n)^2]$.  Replacing $\ep[h(X_1,\ldots,X_n)|X_1]$ and $\mu_n$ by their estimators based on observable data, we can}  design $\hat{S}_n$ as follows
\begin{equation}
n\cdot \hat{S}_n^2 := \frac{r^2}n\sum_{i=1}^n
\underbrace{\Bigg\{\frac1{\binom{n-1}{r-1}}\sum_{\substack{1\leq i_1<\cdots<i_{r-1}\leq n\\i_1,\ldots,i_{r-1}\neq i}} h(A_{i,i_1,\ldots,i_{r-1}}) - \hat{U}_n  \Bigg\}^2}_{ \textrm{Estimates }\xi_1^2 = \var(g_1(X_1))}.
\label{def::our-variance-estimator}
\end{equation}
We will show in Theorem \ref{thm::jackknife} that the {$|\hat S_n^2-\hat S_{n;\textrm{jackknife}}^2|$ is ignorable, but our estimator $\hat S_n$ is computationally more efficient than the jackknife estimator.}  Next, we expand $\hat{T}_n$.  For simplicity, define the following shorthand
\begin{align}
{ U_n^\#} := \frac{1}{\sqrt{n}\cdot \xi_1}\sum_{i=1}^n g_1(X_i), \quad&\Delta_n := \frac{r-1}{\sqrt{n}(n-1)\xi_1} \sum_{1\leq i<j\leq n} g_2(X_i,X_j),
\label{def::delta}\\
\hat{\Delta}_n:=(\hat{U}_n-U_n)/\sigma_n,\quad\delta_n := (\hat{\sigma}_n^2&-\sigma_n^2)/\sigma_n^2,
\quad\textrm{ and }\quad \hat{\delta}_n := (\hat{S}_n^2-\hat{\sigma}_n^2)/\sigma_n^2,\notag
\end{align}
where in \eqref{def::delta}, the technical intermediate term $\hat{\sigma}_n^{ 2}$ is defined as
\begin{align}
n\cdot \hat{\sigma}_n^2 &:=  \frac{r^2}n\sum_{i=1}^n\Bigg\{\frac1{\binom{n-1}{r-1}}\sum_{\substack{1\leq i_1<\cdots<i_{r-1}\leq n\\i_1,\ldots,i_{r-1}\neq i}} h(W_{i,i_1,\ldots,i_{r-1}}) - U_n  \Bigg\}^2.\notag
\end{align}
We now show that $\hat{T}_n$ can be expanded as follows.
\begin{align}
\hat{T}_n
&= \left( { U_n^\#} + \Delta_n + \hat{\Delta}_n + 
{\tOp(n^{-1}\log^{3/2}n)} 
\right)\cdot \left( 1+\hat{\delta}_n + \delta_n \right)^{-1/2}\notag\\
&= \tilde{T}_n + {\check{\Delta}_n} + \textrm{Remainder},
\label{main::hat-T_n::expansion}
\end{align}
where {$\check\Delta_n$ encodes the ``linear part'' of $\hat\Delta_n$ (see Lemma \ref{lemma::term-approx}-\eqref{lemma::term-approx-hat-delta}), recall the definition of $\tOp$ from Section \ref{subsec::tOp}, and define}
\begin{align}
\tilde{T}_n := { U_n^\#} + \Delta_n - \frac12 { U_n^\#}\cdot \delta_n.
\label{def::T-tilde}
\end{align}
{ The remainder in \eqref{main::hat-T_n::expansion} consists of the remainder terms from both $U_n-\mu_n$ and $\hat U_n-U_n$ approximations.  We will show that it is $\tOp({\cal M}(\rho_n,n;R))$.}
The form \eqref{main::hat-T_n::expansion} is partially justified by the Taylor expansion $(1+x)^{-1/2}= 1-x/2 + O(x^2)$, with $x:= (\hat S_n^2-\sigma_n^2)/\sigma_n^2 = O_p(n^{-1/2})$ \citep{maesono1997edgeworth}; and a complete justification comes from our main lemma, i.e. Lemma \ref{lemma::term-approx}.
{Now recalling the definition of acyclic and cyclic $R$ shapes from Definition \ref{definition::acyclic-cyclic-R}, the definition of $\M(\rho_n,n;R)$ from definition \ref{definition::M} in Section \ref{sec::introduction}, and the definition of $\tOp$, we are ready to state our main lemma. }

\begin{lem}  \label{lemma::term-approx}
	Assume the following conditions hold:
	\begin{enumerate}[(i).]
		\item $\rho_n^{-s}\cdot \xi_1 >C>0$,  where $C>0$ is a universal constant,
		\label{condition::xi_1-bounded-away-from-zero}
		\item 
        {
        $\rho_n = \omega(n^{-1})$ for acyclic $R$, or $\rho_n = \omega(n^{-2/r})$ for cyclic $R$,
        }
		\label{condition::rho_n}
	\end{enumerate}
	We have the following results:
	\begin{enumerate}[(a)]
		\item\label{lemma::term-approx-Tn-Unstar-Deltan} $\dfrac{U_n-{\mu_n}}{\sigma_n} = { U_n^\#} + \Delta_n + {\tOp(n^{-1}\cdot \log^{3/2} n)} $,
		\item\label{lemma::term-approx-Delta-hat-conditional-normal} We have 
		$$
		    \hat{\Delta}_n
		    =\frac{(\hat{U}_n-U_n)}{\sigma_n} 
		    = {\check{\Delta}_n + \check{R}_n},
		$$
		{
		where $\check\Delta_n$ and $\check{R}_n$ satisfy
		\begin{align}
		    \check{R}_n 
		    &=
		    \tildehighorderbound
		    \label{eqn::R-check-bound}
		    \\
		    \left\| F_{{\check{\Delta}_n}|W}(u) - \right.
		    &\left.F_{N(0,(\rho_n\cdot n)^{-1}\sigma_w^2)}(u) \right\|_\infty
		    =
		    {\tOp}\left(\rho_n^{-1/2}\cdot n^{-1}\right).
			\label{eqn::Gaussian-term-convergence}
		\end{align}
		where the order control in \eqref{eqn::Gaussian-term-convergence} is $\tOp(\cdot)$ rather than $O(\cdot)$ due to the randomness in $W$.  }The definition of $\sigma_w$ is lengthy and formally stated in Section \ref{appendix::def::sigma_w} in Supplemental Material.  As $n\to\infty$, we have $\sigma_w\stackrel{p}\asymp 1$.
		\item\label{lemma::term-approx-hat-delta} $\hat{\delta}_n =
        {\tildehighorderbound}
		$,
		\item\label{lemma::term-approx-delta} We have
		\begin{equation}
			\delta_n = \frac1n\sum_{i=1}^n\frac{g_1^2(X_i)-\xi_1^2}{\xi_1^2} + \dfrac{2(r-1)}{n(n-1)}\sum_{\substack{1\leq \{i,j\}\leq n\\i\neq j}}\dfrac{g_1(X_i)g_2(X_i,X_j)}{\xi_1^2}  +  {\tOp}(n^{-1}{\cdot \log n}).\notag
		\end{equation}
	\end{enumerate}
\end{lem}
\begin{remark}
    \label{remark::3-1-degenerate}
	Assumption \eqref{condition::xi_1-bounded-away-from-zero} is a standard non-degeneration assumption in literature.  It 
	{is different from a smoothness assumption on graphon $f$}\footnote{{Smooth graphon: $f$ is called \emph{smooth}, if there exists a measure-preserving map $\varrho: [0,1]\to[0,1]$ such that $f(\varrho(\cdot), \varrho(\cdot))$ is a smooth function.  See \cite{gao2015rate, zhang2017estimating} for more details.}}.  A globally smooth Erdos-Renyi graphon leads to a degenerate $g_1(X_1)$ {that $\xi_1^2=\var(g_1(X_1))=0$}.  In the degenerate setting, both the standardization/studentization and the analysis would be very different.  Asymptotic results for $r=2,3$ motifs under an Erdos-Renyi graphon {have been} established by \cite{gao2017testinga, gao2017testingb}.  Degenerate U-statistics are outside the scope of this paper.
\end{remark}
\begin{remark}
    {
    \label{remark::rho-n-lower-bound}
    We note that Lemma \ref{lemma::term-approx} only requires the weak assumption on $\rho_n$ (see Assumption \ref{condition::rho_n}).  This assumption matches the classical sparsity assumptions in network bootstrap literature \citep{bickel2011method, bhattacharyya2015subsampling, green2017bootstrapping}.  Using Lemma \ref{lemma::term-approx}, we prove a higher-order error bound of the Edgeworth expansion in Theorem \ref{thm::main-theorem} with a stronger density assumption; while in Theorem \ref{thm::sparse-main-theorem}, we prove a novel modified Berry-Esseen bound for the normal approximation.  Both downstream theorems significantly improve over existing best results.
    }
\end{remark}
\begin{remark}
	Lemma \ref{lemma::term-approx}-\eqref{lemma::term-approx-Tn-Unstar-Deltan} and \eqref{lemma::term-approx-delta} are similar to results in classical literature on Edgeworth expansion for noiseless U-statistics, but here we account for $\rho_n$.  Parts \eqref{lemma::term-approx-Delta-hat-conditional-normal} and \eqref{lemma::term-approx-hat-delta} are {new results} unique to the network setting.  Especially in the proof of part \eqref{lemma::term-approx-Delta-hat-conditional-normal}, we significantly refine the analysis of the randomness in $A|W$ in \cite{bhattacharyya2015subsampling} and \cite{green2017bootstrapping}.
\end{remark}
\begin{remark}
	\label{remark::knowing-rho_n}
	Our result \eqref{eqn::Gaussian-term-convergence} in Lemma \ref{lemma::term-approx}-\eqref{lemma::term-approx-Delta-hat-conditional-normal} {is different from} Theorem 1 of \cite{bickel2011method}.  {Here we distinguish} three quantities: the true $\rho_n$, $\tilde{\rho}_n = \mathrm{Mean}(W_{ij})$ and $\hat{\rho}_n = \mathrm{Mean}(A_{ij})$.  The convergence rate of $\hat{\rho}_n\to \tilde{\rho}_n$ is much faster than $\tilde{\rho}_n\to\rho_n$.    {The convergence rate here would contribute to the eventual CDF approximation error bound, therefore it is important that }our result \eqref{eqn::Gaussian-term-convergence} corresponds to $\hat{\rho}_n\to \tilde{\rho}_n$, thus avoids the bottleneck.  In contrast, \cite{bickel2011method} and later \cite{lunde2019subsampling} focus on $\hat{\rho}_n\to \rho_n$.
\end{remark}

\begin{table}[h!]
    \centering
    \vspace{-1em}
    \caption{Summary of the main components in $\hat T_n$}
    \begin{tabular}{c|c|c|c}\hline
        Component & Order std. dev. & \bigcell{c}{Impacts\\ Edgeworth formula} & \bigcell{c}{Smoothing\\effect}\\\hline
        $U_n^{\#}$ & $1$ & Yes & No\\
        $\Delta_n-\frac12 U_n^{\#}\cdot \delta_n$ & $n^{-1/2}$ & Yes & No\\
        \rule{0pt}{3ex}$\check\Delta_n$ & $(\rho_n\cdot n)^{-1/2}$ & No & Yes\\
        Remainder & $\tildehighorderbound$ & No & No\\\hline
    \end{tabular}
    \label{tab::Tn-hat-summary}
\end{table}

Overall, Lemma \ref{lemma::term-approx} clarifies the asymptotic orders of the leading terms in the expansion of $\hat{T}_n$.  In fact, Lemma \ref{lemma::term-approx} {has a parallel version for a jackknife $\hat{S}_{n}$} in view of Theorem \ref{thm::jackknife}, but we do not present it due to page limit.
{
We conclude this subsection by a summary table of the main components in $\hat T_n$.
Notice that despite smoother $\check\Delta_n$ is $\Omega(n^{-1/2})$, it does \emph{not} distort the $n^{-1/2}$ term in the Edgeworth expansion formula.  Similar phenomenon is observed in the i.i.d. setting, see \cite{singh1981asymptotic} (equation (2.8)) and \cite{lahiri1993bootstrapping} (Section 2.2).
}

\subsection{Population and   empirical Edgeworth expansions for network moments}

In this subsection, we present our main theorems.

\begin{thm}[Population network Edgeworth expansion]
	\label{thm::main-theorem}
	
	Define
	\begin{align}
		G_n(x) &:= \Phi(x) + \frac{\phi(x)}{\sqrt{n}\cdot \xi_1^3} \cdot \Bigg\{ \frac{ 2 x^2 + 1}{6} \cdot\ep[g_1^3(X_1)]\notag\\
		& + \frac{r-1}2\cdot \left(  x^2 + 1 \right)\ep[g_1(X_1)g_1(X_2)g_2(X_1,X_2)]\Bigg\},
		\label{eqn::main-theorem-1}
	\end{align}
	where $\Phi(x)$ and $\phi(x)$ are the CDF and PDF of $N(0,1)$. 
	{
	Assume condition (\ref{condition::xi_1-bounded-away-from-zero}) of Lemma \ref{lemma::term-approx} hold, and replace condition \eqref{condition::rho_n} by a stronger assumption that either $R$ is acyclic and $\rho_n = \omega(n^{-1/2})$, or $R$ is cyclic and $\rho_n = \omega(n^{-1/r})$.
	}
	Additionally, assume either $\rho_n = O((\log n)^{-1})$ or Cramer's condition $\limsup_{t\to\infty}\left| \ep\left[ e^{\ii t g_1(X_1)\cdot \xi_1^{-1}} \right] \right|<1$ holds.
	We have
	\begin{equation}
		\left\| F_{\hat{T}_n}(x) - G_n(x) \right\|_\infty = \Ohighorderbound.\notag
	\end{equation}
\end{thm}

\begin{remark}
	The assumed $\rho_n$'s upper bound in absence of Cramer's condition serves to sufficiently boost the smoothing power of $\check{\Delta}_n$, quantified in Lemma \ref{lemma::term-approx}-\eqref{eqn::Gaussian-term-convergence}.  This assumption is unlikely improvable, since its required Gaussian variance $(\rho_n\cdot n)^{-1}= \Omega(\log n\cdot n^{-1})$ matches the minimum Gaussian standard deviation requirement $\Omega((\log n)^{1/2}\cdot n^{-1/2})$ in Remark 2.4 in \cite{lahiri1993bootstrapping} for the i.i.d. setting.
\end{remark}

In \eqref{eqn::main-theorem-1}, the Edgeworth coefficients depend on true population moments.  In practice, they need to be estimated from data.  Define
\begin{align}
\hat{g}_1(X_i) &:= \frac1{\binom{n-1}{r-1}}\sum_{\substack{1\leq i_1<\ldots<i_{r-1}\leq n\\i_1,\ldots,i_{r-1}\neq i}}h(A_{i,i_1,\ldots,i_{r-1}}) - \hat{U}_n,\notag\\
\hat{g}_2(X_i,X_j) &:= \frac1{\binom{n-2}{r-2}}\sum_{\substack{1\leq i_1<\ldots<i_{r-2}\leq n\\i_1,\ldots,i_{r-2}\neq i,j}}h(A_{i,j,i_1,\ldots,i_{r-2}}) - \hat{U}_n  - \hat{g}_1(X_i) - \hat{g}_1(X_j),\notag
\end{align}
where we write ``$\hat{g}_1(X_i)$'' rather than ``$\hat{g_1(X_i)}$'' for cleanness.  We stress that the evaluation of $\hat{g}_1(X_i)$ and $\hat{g}_2(X_i,X_j)$ does \emph{not} require knowing the latent $X_i,X_j$.  The Edgeworth coefficients can be estimated by
\begin{align}
\hat{\xi}_1^2 := \frac{n\cdot \hat{S}_n^2}{r^2} = \frac1n\sum_{i=1}^n\hat{g}_1^2(X_i),&\quad\textrm{ and }\quad \hat{\ep}\left[ g_1^3(X_1) \right] := \frac1n\sum_{i=1}^n \hat{g}_1^3(X_i),\notag\\
\hat{\ep}\left[ g_1(X_1)g_1(X_2)g_2(X_1,X_2) \right] &:= \frac1{\binom{n}2} \sum_{1\leq i<j\leq n}\hat{g}_1(X_i)\hat{g}_1(X_j)\hat{g}_2(X_i,X_j).\notag
\end{align}

\begin{thm}[Empirical network Edgeworth expansion]
	\label{thm::main-empirical}
	Define the empirical Edgeworth expansion as follows:
	\begin{align}
	\hat{G}_n(x) &:= \Phi(x) + \frac{\phi(x)}{\sqrt{n}\cdot \hat{\xi}_1^3} \cdot \Bigg\{  \frac{ 2 x^2 + 1}{6} \cdot \hat{\ep}[g_1^3(X_1)]\notag\\
	& + \frac{r-1}2\cdot \left(  x^2 + 1 \right)\hat{\ep}[g_1(X_1)g_1(X_2)g_2(X_1,X_2)]\Bigg\},
	\label{Edgeworth::empirical}
	\end{align}
	Under the conditions of Theorem \ref{thm::main-theorem}, we have
	\begin{equation}
		\left\| F_{\hat{T}_n}(x) - \hat{G}_n(x) \right\|_\infty = {\tOp}({\cal M}(\rho_n,n;R)).\notag
	\end{equation}
\end{thm}
\begin{remark}
	{Another approach to estimate the unknown coefficients in Edgeworth expansion is by bootstrap.}  The concentration of $\hat{G}_n\to G_n$ should not be confused with the concentration $\hat{G}_n^*\to \hat{G}_n$, where $\hat{G}_n^*$ is the expansion with bootstrap-estimated coefficients.  See literature regarding the i.i.d. setting \citep{helmers1991edgeworth,maesono1997edgeworth}.  In $\hat{G}_n^*\to \hat{G}_n$, the convergence rate is not a concern, because, without constraining computation cost, one can let the number of bootstrap samples grow arbitrarily fast, so the proof of bootstrap concentration only requires consistency, but our proof regarding $\hat{G}_n\to G_n$ requires careful rate calculations.
\end{remark}
Next, we show that different choices of the variance estimators for studentization represent no essential discrepancy.
\begin{thm}[Studentizing by a jackknife variance estimator \eqref{def::jackknife-variance-estimator}]
	\label{thm::jackknife}
	Define
	$$
	\hat{T}_{n;\mathrm{jackknife}} := \frac{\hat{U}_n-\mu_n}{\hat{S}_{n;\mathrm{jackknife}}}.
	$$
	Under the assumptions of Theorem \ref{thm::main-theorem}, we have
	\begin{align}
	| \hat{S}_n - \hat{S}_{n;\mathrm{jackknife}} | &= O(\hat S_n\cdot n^{-1}),
	\label{proofeqn::jackknife-temp}\\
	\left\| F_{\hat{T}_{n;\mathrm{jackknife}}}(x) - G_n(x) \right\|_\infty &= \Ohighorderbound,\notag\\
	\left\| F_{\hat{T}_{n;\mathrm{jackknife}}}(x) - \hat{G}_n(x) \right\|_\infty &= {\tOp}({\cal M}(\rho_n,n;R)).\notag
	\end{align}
\end{thm}
Theorem \ref{thm::jackknife} states that on statistical properties, one does not need to differentiate between $\hat{T}_n$ and $\hat{T}_{n;\mathrm{jackknife}}$.  {The evaluation of $\hat S_{n;\textrm{jackknife}}$ costs $O(n^{r+1})$ time because each individual $\hat U_n^{(-i)}$ costs $O(n^r)$; whereas our estimator $\hat S_n$ costs $O(n^r)$.  Our estimator also has a more convenient form for theoretical analysis.}

\subsection{Remarks on non-smooth graphons and a comparison table of our results with literature}


Our results do not assume graphon smoothness or low-rankness.  This aligns with the literature on noiseless U-statistics but sharply contrasts network inferences based on model parameter estimation such as \cite{hoff2002latent,lei2016goodness} and network bootstraps based on model estimation \citep{green2017bootstrapping, levin2019bootstrapping}.  Notice that the concept ``non-smoothness'' usually emphasizes ``not assuming smoothness'' rather than explicitly describing irregularity.  It is a very useful tool for modeling networks with high structural complexity or unbalanced observations, examples include: (1) a small group of \emph{outlier} nodes that behave differently from the main network patterns \cite{cai2015robust}; (2) in networks exhibiting ``core-periphery'' structures \citep{della2013profiling, zhang2015identification}, we may wish to relax structural assumptions on periphery nodes due to scarcity of observations; and (3) networks generated from a mixture model \citep{newman2007mixture} with many small-probability mixing components may appear non-smooth in these parts.  Unfortunately, existing research on practical methods for non-smooth graphons is rather limited due to the obvious technical difficulty, but exceptions include \citep{choi2017co}.

Our results send the surprising message that under mild conditions, the sampling distribution of a network moment is still \emph{smooth} and can be \emph{accurately} approximated, even if the graphon is non-smooth.

{
\subsection{Sparse networks}
\label{subsec::sparse-networks}
We have been focusing on discussing mildly sparse networks, but many networks tend to be sparse \citep{guedon2016community}.  In this section, we investigate the following sparsity regime
\begin{equation}
    \rho_n:\ 
    \begin{cases}
        n^{-1}\prec \rho_n \preceq n^{-1/2}, &\textrm{ for acyclic }R\\
        n^{-2/r}\prec \rho_n \preceq n^{-1/r}, &\textrm{ for cyclic }R\\
    \end{cases}
    \label{rho_n-sparse-regime}
\end{equation}
It turns out that the Berry-Esseen bound would deteriorate and be worse than the conventional $n^{-1/2}$ in the i.i.d. and the noiseless U-statistic settings.  The exact reason is technical and will be better seen in the proof of Theorem \ref{thm::sparse-main-theorem}, but the intuitive explanation is that if $\rho_n$ is too small, the higher degree ($\geq 2$) random errors in $\hat U_n-U_n$, caused by the randomness in $A|W$, diminishes too slowly compared to the scale of the demoninator of $\hat T_n$.  If the network sparsity $\rho_n$ falls below the typically assumed lower bounds: $n^{-1}$ for acyclic $R$ and $n^{-2/r}$ for cyclic $R$ \citep{bickel2011method, bhattacharyya2015subsampling, green2017bootstrapping}, then no known consistency guarantee exists.
\begin{thm}\label{thm::sparse-main-theorem}
    Under the conditions of Lemma \ref{lemma::term-approx}, except replacing Condition \eqref{condition::rho_n} by \eqref{rho_n-sparse-regime}, we have the following modified Berry-Esseen bound
    \begin{align}
        \left\|
            F_{\hat T_n}(u) - G_n(u)
        \right\|_\infty
        &\asymp
        \left\|
            F_{\hat T_n}(u) - \Phi(u)
        \right\|_\infty=
        \Ohighorderbound \bigwedge o(1),\notag
    \end{align}
    where recall that $\Phi(\cdot)$ is the CDF of $N(0,1)$.  Moreover, 
    $$
        \left\|
            F_{\hat T_n}(u) - \hat G_n(u)
        \right\|_\infty
        =
        \begin{cases}
            \tildehighorderbound, &\textrm{ if }\rho_n=\omega(n^{-1}\log^{1/2} n),\\
            o_p(1), &\textrm{ if }\rho_n = \omega(n^{-1})
        \end{cases}
    $$
    %
\end{thm}
In the sparse regime, we could not control the uniform CDF approximation error bound below $n^{-1/2}$.  Consequently, using Edgeworth expansion would not bring asymptotic merit compared to a simple $N(0,1)$ approximation.  On the other hand, the conclusion of Theorem \ref{thm::sparse-main-theorem} connects the error bound results for dense and sparse regimes.  Interestingly, as the order of $\rho_n$ decreases from $n^{-1/2}$ to $n^{-1}$ for acyclic $R$, or from $n^{-1/r}$ to $n^{-2/r}$ for cyclic $R$, we see a gradual depreciation in the uniform CDF approximation error from the order of $n^{-1/2}$ to merely uniform consistency.  The classical literature only studied the boundary cases ($\rho_n=\omega(n^{-1})$ or $\rho_n=\omega(n^{-2/r})$, depending on $R$), and our result here reveals the complete picture.



}

We conclude this section by comparing our results to some representative works in classical and very recent literature.
\begin{table}[ht!]
	\centering
	\vspace{-1em}
	\caption{Comparison of CDF approximation methods for noisy/noiseless studentized U-statistics}
	\label{table::comparison-of-results}
	\begin{adjustbox}{center}
		\begin{tabular}{|c|c|c|c|c|c|c|}\hline
			Method & \bigcell{c}{U-stat.\\type} & \bigcell{c}{Popul.\\momt.\footnotemark\label{footnote::knownpopmoments}} & \bigcell{c}{Smooth\\graphon} & \bigcell{c}{Non lat.\\/Cramer} & \bigcell{c}{Network sparsity\\assumption on $\rho_n$\footnotemark\label{footnote::sparsity-assumption}} & \bigcell{c}{CDF approx.\\error rate}\\\hline
			
			\multirow{3}{*}{\bigcell{c}{Our method\\(empirical Edgeworth)}} & \multirow{3}{*}{Noisy} & \multirow{3}{*}{No} & \multirow{3}{*}{No} & If yes & $\omega(n^{{-2/r}})$(C); $\omega(n^{{-1}})$(Ac)\footnotemark\label{footnote::motifshape} & 
			\rule{0pt}{3ex}{$\tildehighorderbound \wedge o_p(1)$} {\bf (H)}\footnotemark\label{footnote::higher-order-accuracy}
			\rule[-1.2ex]{0pt}{0pt}\\\cline{5-7}
			& &  & & If no & \bigcell{c}{$\omega(n^{{-2/r}})$(C); $\omega(n^{{-1}})$(Ac)\\and $O\left((\log n)^{-1}\right)$(C, Ac)} & {$\tildehighorderbound \wedge o_p(1)$} {\bf (H)}\\\hline
			
			\bigcell{c}{Node re-/sub- sampling\\justified by our theory} & Noisy & No & No & Yes & $\omega(n^{-1/r})$(C); $\omega(n^{-1/2})$(Ac) & $o_p(n^{-1/2})$ {\bf (H)}\\\hline
			
			\citet{bickel2011method} & Noisy & No\footnotemark\label{footnote::knowing-true-rho} & No & No & $\omega(n^{-2/r})$(C); $\omega(n^{-1})$(Ac) & Consistency  \\\hline
			
			\citet{bhattacharyya2015subsampling} & Noisy & No & No & No & $\omega(n^{-2/r})$(C); $\omega(n^{-1})$(Ac) & Consistency\\\hline
			
			\citet{green2017bootstrapping} & Noisy & No & Mixed\footnotemark\label{footnote::greenshalizismoothness} & No & $R$ is Ac; or $\omega(n^{-1/(2r)})$(C)\footnotemark\label{footnote::greenshalizisparsity} &  Consistency\\\hline
			
			\citet{levin2019bootstrapping} & Noisy & No & 
			{Low-rank\footnotemark\label{footnote::levinlevinaassumptions}} & No &  $\omega(n^{-1}\cdot \log n)$ (Ac*)\footnotemark\label{footnote::levinalevinaresults} & Consistency\\\hline

			Bickel, Gotze and Zwet \cite{bickel1986edgeworth} & Noiseless & Yes & No & Yes & Not applicable& ${o}(n^{-1})$ {\bf (H)}\\\hline

			Bentkus, Gotze and Zwet \cite{bentkus1997edgeworth} & Noiseless & Yes & No & Yes & Not applicable& ${O}(n^{-1})$ {\bf (H)}\\\hline
			
			Putter and Zwet \cite{putter1998empirical} & Noiseless & No & No & Yes & Not applicable & $o_p(n^{-1/2})$ {\bf (H)}\\\hline
			
			\citet{bloznelis2003edgeworth} & Noiseless & No & No & Yes & Not applicable& $o_p(n^{-1/2})$ {\bf (H)}\\\hline
		\end{tabular}
	\end{adjustbox}
	\vspace{-1em}
\end{table}
\footnotetext[\getrefnumber{footnote::knownpopmoments}]{``Yes'' means need to know the population moments that appear in Edgeworth coefficients, i.e. $\xi_1$, $\ep[g_1^3(X_1)]$ and $\ep[g_1(X_1)g_1(X_2)g_2(X_1,X_2)]$.}
\footnotetext[\getrefnumber{footnote::sparsity-assumption}]{To compare $\rho_n$ assumptions, see our Remark \ref{remark::rho-n-lower-bound}}
\footnotetext[\getrefnumber{footnote::motifshape}]{(C): {\bf c}yclic $R$; (Ac): {\bf ac}yclic $R$.}
\footnotetext[\getrefnumber{footnote::higher-order-accuracy}]{Recall ${\cal M}(\rho_n,n; R)$ defined in \eqref{M} {and $\tOp$ defined in Section \ref{subsec::tOp}}.  {\bf (H)}: higher-order accuracy results. ``Consistency'': only convergence, no error rate.}
\footnotetext[\getrefnumber{footnote::knowing-true-rho}]{In \cite{bickel2011method, bhattacharyya2015subsampling, lunde2019subsampling}, $\hat{U}_n-\mu_n$ was rescaled by $\rho_n$ and $n$.  Whether assuming the knowledge of the true $\rho_n$ or not does not matter for their $o_p(1)$ error bound, but it would make a difference if an $o_p(n^{-1/2})$ or finer bound is desired.  See our Remark \ref{remark::knowing-rho_n}.}
\footnotetext[\getrefnumber{footnote::greenshalizismoothness}]{The bootstrap based on denoised $A$ requires smoothness.  See Theorem 2 of \cite{green2017bootstrapping}.}
\footnotetext[\getrefnumber{footnote::greenshalizisparsity}]{It seems their assumption for cyclic $R$ was a typo, and $\rho_n=\omega(n^{-2/r})$ should suffice.  Also, they used \cite{bhattacharyya2015subsampling} in their proof, which requires $\rho_n = \omega(n^{-1})$ for (Ac).}
\footnotetext[\getrefnumber{footnote::levinlevinaassumptions}]{\cite{levin2019bootstrapping} assumed the graphon rank is low and known.}
\footnotetext[\getrefnumber{footnote::levinalevinaresults}]{(Ac*): They require the motif to be either acyclic or an $r$-cycle, see their Theorem 4.  Their Theorem 3 requires condition (8) that only holds when $R$ is a clique.}

\section{Theoretical and methodological applications}
\label{sec::applications}
\subsection{Higher-order accuracy of node sub- and re-sampling network bootstraps}
\label{subsec::bootstrap-accuracy}


One important corollary of our results is first higher-order accuracy proof of some network bootstrap schemes.  For a network bootstrap scheme that produces an estimated $\hat{U}_{n^*}^b$ and its jackknife\footnote{Here, we use the jackknife estimator in the bootstrap for a better connection with existing literature in the proof.} variance estimator $\hat{S}_{n^*}^*$, define $\hat T_{n^*}^*=(\hat U_{n^*}^b-\hat U_n)/\hat{S}_{n^*}^*$.  We are going to justify the following two schemes.

\begin{enumerate}[(a).]
	\item Sub-sampling \citep{bhattacharyya2015subsampling}: randomly sample $n^*$ nodes from $\{1,\ldots,n\}$ \emph{without replacement}, and compute $\hat{T}^*_{n^*}$ from the induced sub-network of $A$.
	\item Re-sampling \citep{green2017bootstrapping}: random sample $n$ nodes from $\{1,\ldots,n\}$ \emph{with replacement}, and compute $\hat{T}^*_{n^*}$ from the induced sub-network of $A$.
\end{enumerate}
\begin{remark}
	Notice that \cite{green2017bootstrapping} did not study the studentized form, and \cite{bhattacharyya2015subsampling} proposed a different variance estimator (what they call ``$\hat{\sigma}_{B_i}$'').  Our justifications focus on the \emph{sampling schemes} combined with some natural formulation, not necessarily the same formulation as in these two papers.
\end{remark}
\begin{remark}
	As noted in \cite{green2017bootstrapping}, scheme (b) can be viewed as our data generation procedure described in Sections \ref{subsec::problem-setup-graphon-model} and \ref{subsec::problem-setup-network-moments} but with the graphon $f$ replaced by the adjacency-induced graphon $A(u,v) = A_{\lceil nu\rceil, \lceil nv\rceil}$, where $\lceil y\rceil :=\mathrm{Ceiling}(y)$.  {We discuss their scheme (a) in Section \ref{sec::discussion}.}  This may seem similar to $f$-based data generaion, but in fact they are distinct.  The graphon $A(\cdot,\cdot)$ inherits the binary nature of $A$ and will necessarily yield a lattice $g_1^*(X_1^*)$ regardless of the original graphon $f$ and the motif ${\cal R}$, rendering most classical Edgeworth analysis techniques inapplicable.  But the real obstacle is that the bootstrapped network data from $A(\cdot,\cdot)$ have no edge-wise observational errors (i.e. no counterpart to the randomness in $A|W$).  Consequently, $\hat{T}_{n^*}^*$ loses the self-smoothing feature that $\hat{T}_n$ enjoys.
\end{remark}

\begin{thm}\label{thm::bootstrap-accuracy}
	Assume $g_1(X_1)$ satisfies a Cramer's condition such that $\limsup_{t\to\infty}\left|\ep\left[ e^{\ii t g_1(X_1)\cdot \xi_1^{-1}} \right]\right|<1$. Under the conditions of Theorem \ref{thm::main-empirical}, we conclude for the following bootstrap schemes:
	\begin{enumerate}[(a).]
		\item Sub-sampling:  choosing $n^* \asymp n$ and $n-n^* \asymp n$, we have
			\begin{equation}
			\label{theorem::application::boot-sub-sampling}
				\left\|  F_{\hat{T}_{n^*}^*}(u) - F_{\hat{T}_{n^*(1-n^*/n)}}(u)  \right\|_\infty = o_p\left((n^*)^{-1/2}\right) = o_p(n^{-1/2}).
			\end{equation}
		\item Re-sampling:  choosing $n^* = n$, we have
			\begin{equation}
			\label{theorem::application::boot-re-sampling}
				\left\|  F_{\hat{T}_{n^*}^*}(u) - F_{\hat{T}_{n^*}}(u)  \right\|_\infty = o_p\left((n^*)^{-1/2}\right) = o_p(n^{-1/2}).
			\end{equation}
	\end{enumerate}
	
\end{thm}

\begin{remark}
	In the proof of Theorem \ref{thm::bootstrap-accuracy}, we combined our main results with the results of \cite{bloznelis2003edgeworth} for finite population U-statistics.  It is important to notice that all existing works under the finite populations did assume non-lattice with population size growing to infinity, see condition (1.13) in Theorem 1 of \cite{bloznelis2003edgeworth}.  Consequently, the higher-order accuracy of some network bootstraps is only proved under Cramer's condition by so far.
\end{remark}

Part (a) of Theorem \ref{theorem::application::boot-sub-sampling} quantifies the \emph{effective sample size} in the sub-sampling network bootstrap:  sampling $n^*$ out of $n$ nodes without replacement, the resulting bootstrap $\hat{T}_{n^*}^*$ approximates the distribution of $\hat{T}_m$ where $m = \left\{n^*/n\cdot (1-n^*/n)\right\}\times n$.  Consequently, in order to approach the sampling distribution of $\hat{T}_n$ with higher-order accuracy using sub-sampling \cite{bhattacharyya2015subsampling}, one must have an observed network of at least $4n$ nodes, from which she shall repeatedly sub-sample $2n$ nodes without replacement.

\subsection{One-sample t-test for network moments under general null graphon models}
\label{subsec::one-sample t-test}

In this and the next subsections, we showcase how our results immediately lead to useful inference procedures for network moments.  For a given motif $R$, we test on its population mean frequency $\mu_n$.  Since $\mu_n$ depends on $n$ through $\rho_n$, we formulate the hypotheses as follows
\begin{center}
	$H_0: \mu_n = c_n$, versus $H_a: \mu_n \neq c_n$.
\end{center}
where $c_n$ is a speculated value of $\mu_n = \ep[h(A_{1,\ldots,r})]$.  In practice, $c_n$ may come from a prior study on a similar data set or fitting a speculated model to the data (for concrete examples on $c_n$ guesses, see Section 6.1 of \cite{bhattacharyya2015subsampling}).

Here for simplicity we only discuss a two-sided alternative, and one-sided cases are exactly similar.  The p-value can be formulated using our empirical Edgeworth expansion $\hat{G}_n(\cdot)$ in \eqref{Edgeworth::empirical}:
\begin{equation}
\textrm{Estimated p-value} = 2\cdot \min\left\{ \hat{G}_n(t^{\mathrm{(obs)}}),  1-\hat{G}_n(t^{\mathrm{(obs)}}) \right\}.
\label{eqn::application::one-sample-test}
\end{equation}
where $t^{\mathrm{(obs)}} := (\hat{u}_n^{\mathrm{(obs)}} - {c_n})/\hat{s}_n^{\mathrm{(obs)}}$, and $\hat{u}_n^{\mathrm{(obs)}}$ and $\hat{s}_n^{\mathrm{(obs)}}$ are the observed $\hat{U}_n$ and $\hat{S}_n$, respectively.  We have the following explicit Type-II error rate.

\begin{thm}
	\label{theorem::application-one-sample-test}
	Under the conditions of Theorem \ref{thm::main-empirical}, we have the following results:
	\begin{enumerate}
		\item The Type-I error rate of test \eqref{eqn::application::one-sample-test} is $\alpha + {\Ohighorderbound}$.
		\item The Type-II error rate of this test is ${o}(1)$
		when $|c_n-{\mu_n}| = \omega\left(\rho_n^s\cdot n^{-1/2}\right)$.
	\end{enumerate}
\end{thm}

\begin{remark}
	The null model we study is complementary to the degenerate Erdos-Renyi null model in \citep{lei2016goodness, gao2017testinga, gao2017testingb}.  The scientific questions are also different: they test model goodness-of-fit whereas we test population moment values.  Notice that distinct network models may possibly share some common population moments.  These approaches also use very different methods and analysis techniques.
\end{remark}

\subsection{Cornish-Fisher confidence intervals for network moments}
\label{subsec::cornish-fisher-CI}

Noticing that $\hat{G}_n$ is almost never a valid CDF, in order to preserve the higher-order accuracy of $\hat{G}_n$, we use the Cornish-Fisher expansion \citep{cornish1938moments, fisher1960percentile} to approximate the quantiles of $F_{\hat{T}_n}$.   A Cornish-Fisher expansion is the inversion of an Edgeworth expansion, and its validity hinges on the validity of its corresponding Edgeworth expansion.  Using the technique of \cite{hall1983inverting}, we have
\begin{thm}
	\label{theorem::application-CI}
	
	{
    For any $\alpha\in(0,1)$, define the lower $\alpha$ quantile of the distribution of $\hat{T}_n$ by
    \begin{equation}
        q_{\hat{T}_n;\alpha} := \arg\inf_{q\in\mathbb{R}}\  F_{\hat T_n}(q)\geq \alpha
    \end{equation}
    }
    {
    and define the approximation
    \begin{align}
    	\hat{q}_{\hat{T}_n;\alpha}&:=z_{\alpha}-\frac{1}{\sqrt{n}\cdot \hat{\xi}_1^3} \cdot \Bigg\{ \frac{ 2z_{\alpha}^2 + 1}6\cdot \hat{\ep}[g_1^3(X_1)]\notag\\
    	& + \frac{r-1}2\cdot \left(  z_{\alpha}^2 + 1 \right)\hat{\ep}[g_1(X_1)g_1(X_2)g_2(X_1,X_2)]\Bigg\},
    	\label{formula::quantile-for-CI}
	\end{align}
	where $z_{\alpha} := \Phi^{-1}(\alpha)$.  Then under the conditions of Theorem \ref{thm::main-empirical}, we have the following results
	\begin{enumerate}[(a).]
	    \item The discrepancy between nominal and actual percentage-below for $q_{\hat T_n;\alpha}$ is bounded by
	    \begin{equation}
    	    |F_{\hat T_n}(q_{\hat T_n;\alpha}) - \alpha| = \Ohighorderbound
    	    \label{application::CI-quantile-discrepancy}
    	\end{equation}
    	
    	\item A ``horizontal'' error bound:
    	\begin{equation}
        	\left| \hat{q}_{\hat{T}_n;\alpha}  -  q_{\hat{T}_n;\alpha} \right| = \tildehighorderbound 
    	\label{application::CI-q-q-hat-error-bound}
    	\end{equation}
    	
	    \item A uniform ``vertical'' error bound
	    \begin{equation}
    		\pr(\hat{T}_n\leq \hat{q}_{\hat{T}_n;\alpha}) = \alpha + \Ohighorderbound.
    		\label{application::CI-coverage-error-bound}
    	\end{equation}
    	
	\end{enumerate}
	}
\end{thm}

{
The vertical error bound describes the approximation error between the nominal and actual coverage probabilities, whereas the horizontal error bound governs the approximation of quantiles.  Using the vertical error bound, a $1-\alpha$ two-sided symmetric Cornish-Fisher confidence interval for estimating $\mu_n$ can be easily constructed as follows}
\begin{equation}
    \left( {\hat{U}_n}-\hat{q}_{\hat{T}_n;1-\alpha/2}\cdot {\hat S_n}, {\hat{U}_n}-\hat{q}_{\hat{T}_n;\alpha/2}\cdot {\hat S_n}  \right)    
    \label{formula::EEE-confidence-interval}
\end{equation}
and by Theorem \ref{theorem::application-CI}, we know this CI has a $1-\alpha + {\Ohighorderbound}$ coverage probability.  One-sided confidence intervals can be constructed similarly.

\section{Simulations}
\label{sec::simulations}

\subsection{{Simulation 1:  Higher-order accuracy of empirical Edgeworth expansion}}
\label{subsec::simulation-1}

{In the first simulation, }our numerical studies focus on the CDF of $F_{\hat{T}_n}$.  In an illustrative example, we simulate with a lattice $g_1(X_1)$ and show the distinction between $F_{\hat{T}_n}$ and $F_{T_n}$ that clearly illustrates the self-smoothing effect in $\hat{T}_n$.  Then we systematically compare the performance of our empirical Edgeworth expansion to benchmarks that demonstrates the clear advantage of our method in both accuracy and computational efficiency.

We begin by describing the basic settings.  We range the network size $n$ in an exponentially spaced set $n\in\{10,20,40,80,160\}$.  Synthetic network data are generated from three graphons marked by their code-names in our figures:  (1). \texttt{"BlockModel":} This is an ordinary stochastic block model with $K=2$ equal-sized communities and the following edge probabilities $B = (0.6,0.2;0.2,0.2)$;  (2). \texttt{"SmoothGraphon":} Graphon 4 in \cite{zhang2017estimating}, i.e. $f(u,v):=(u^2+v^2)/3\cdot \cos(1/(u^2+v^2))+0.15$.  This graphon is smooth and full-rank \citep{zhang2017estimating};  (3). \texttt{"NonSmoothGraphon"}\citep{choi2017co}:  We set up a high-fluctuation area in a smooth $f$ to emulate the sampling behavior of a non-smooth graphon, as follows
$$
f(u,v):=0.5 \cos\left\{ 0.1/((u-1/2)^2+(v-1/2)^2)^{-1} + 0.01 \right\}\max\{u,v\}^{2/3}+0.4.
$$
We test the {four} simplest motifs: \emph{edge}, \emph{triangle}, \emph{V-shape}\footnote{A ``V-shape'' is the motif obtained by disconnecting one edge in a triangle.  In the language of \cite{bickel2011method}, it is a 2-star.}{, and a \emph{three-star} among 4 nodes with edge set $\{(1,2),(1,3),(1,4)\}$}.  The main computation bottleneck lies in the evaluation of $F_{\hat{T}_n}$.  Network bootstraps also becomes costly as $n$ increases.

The benchmarks are:  1. $N(0,1)$ (its computation time is deemed zero and not compared to others); 2. sub-sampling by \cite{bhattacharyya2015subsampling} with $n^*=n/2$; 3. re-sampling $A$ by \cite{green2017bootstrapping}; 4. latent space bootstrap called ``ASE plug-in'' defined in Theorem 2 of \cite{levin2019bootstrapping}.  Notice that we equipped \cite{levin2019bootstrapping} with an adaptive network rank estimation\footnote{Consequently, our enhanced version of this benchmark can decently denoise some smooth but high-rank graphons, in view of the remarks in \cite{zhang2017estimating} and the results of \cite{xu2018rates}.} by USVT \citep{chatterjee2015matrix}.

For each (graphon, motif, $n$) tuple,  we first evaluate the true sampling distribution of $\hat{T}_n$ by a Monte-Carlo approximation that samples $n_{\mathrm{MC}}:= 10^6$ networks from the graphon.  Next we start $30$ repeated experiments: in each iteration, we sample $A$ from the graphon and approximate $F_{\hat{T}_n}$ by all methods, in which we draw $n_{\mathrm{boot}}=2000$ bootstrap samples for each bootstrap method -- notice that this is 10 times that in \cite{levin2019bootstrapping}.  We compare
\begin{equation}
\mathrm{Error}(\hat{F}_{\hat{T}_n}) := \sup_{u\in[-2,2]; 10u\in\mathbb{Z}}\left|\hat{F}_{\hat{T}_n}(u) - F_{\hat{T}_n}(u)\right|.
\label{simu::loss-function}
\end{equation}
\begin{remark}
	We need many Monte-Carlo repetitions, because the uniform accuracy of the empirical CDF of an i.i.d. sample is only $O_p(n_{\mathrm{MC}}^{-1/2})$ \citep{dvoretzky1956asymptotic, kosorok2007introduction}, and for the noiseless and noisy U-statistic setting, the bound might be worse than the i.i.d. setting due to dependency\footnote{This is not to be confused with the Edgeworth approximation error bound.  In this Monte Carlo procedure, both the true and approximate $F_{\hat{T}_n}$ are oracle.}.  Therefore, we set $n_{\mathrm{MC}}\gg \max(n^2)= 160^2$ to prevent the errors of the compared methods being dominated by the error of the Monte-Carlo procedure; while keep our simulations reproducible with moderate computation cost.  We did find smaller $n_{\mathrm{MC}}$ such as $10^5$ to cloud the performance of our method.
\end{remark}

\begin{figure}[h!]
	\centering
	\includegraphics[width=0.95\textwidth]{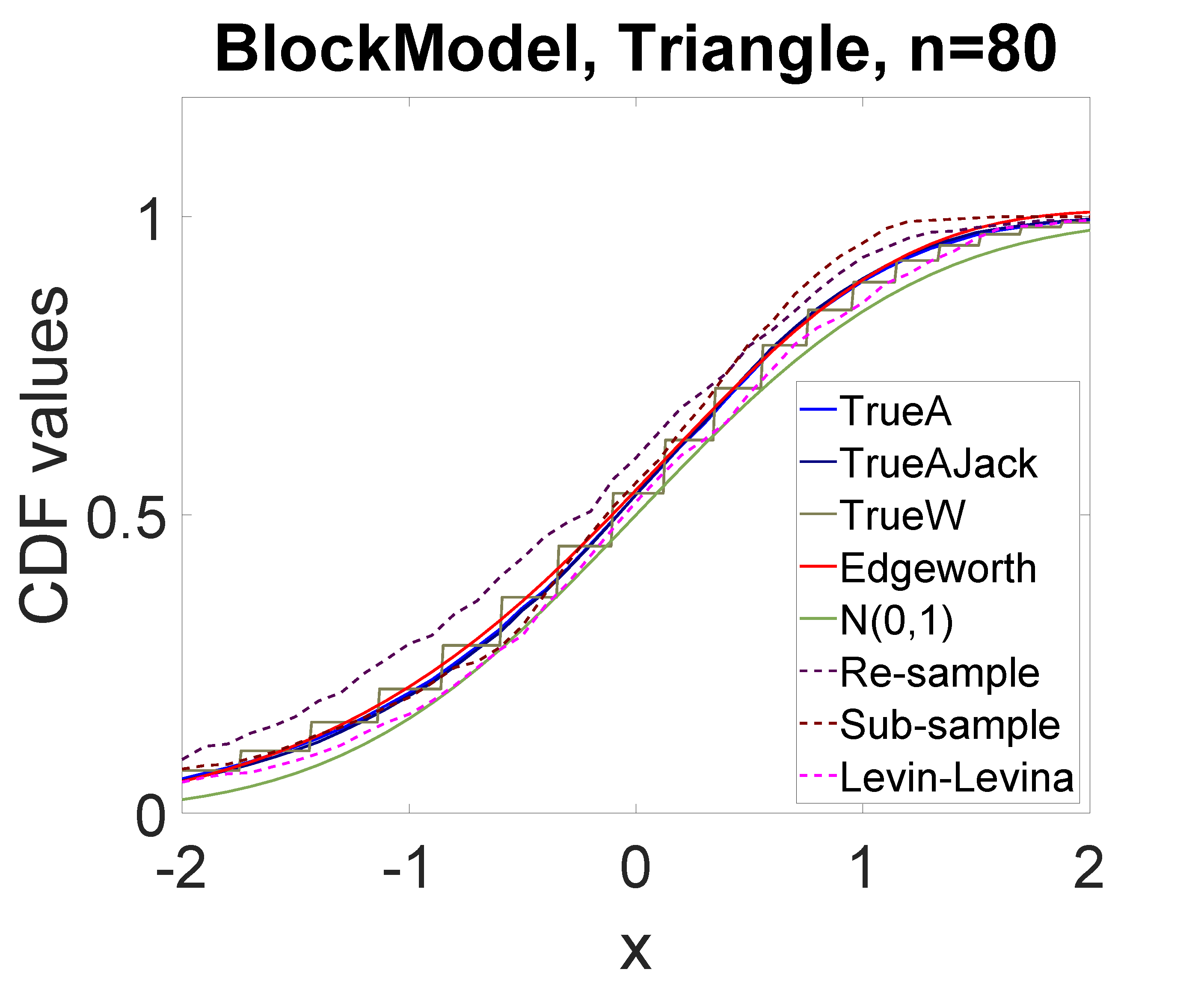}
	\caption{CDF curves of the studentization forms and approximations.  Network size $n=80$.  The graphon is the ``\texttt{BlockModel}'' we described earlier in this section, and the motif is triangular.  Each bootstrap method draws 500 random samples.  \texttt{TrueA} is $F_{\hat{T}_n}$; \texttt{TrueAJack} is $F_{\hat{T}_{n;\mathrm{jackknife}}}$; \texttt{TrueW} is $F_{T_n}$; \texttt{Edgeworth} is our empirical Edgeworth expansion; \texttt{Re-sample} is node re-sampling $A$ in \cite{green2017bootstrapping}; \texttt{Sub-sample} is node sub-sampling $A$ in \cite{bhattacharyya2015subsampling}; \texttt{Levin-Levina} is the ``ASE plug-in'' bootstrap in \cite{levin2019bootstrapping}.}
	\label{fig::illustration-1}
\end{figure}

Now we present the results. 
We first present the illustrative simulation for just one specific setting.  Figure \ref{fig::illustration-1} shows the distribution approximation curves under a block model graphon that yields a lattice $g_1(X_1)$.   Lines correspond to the population CDF of $\hat{T}_n$, its jackknife version and noiseless version, all evaluated by Monte-Carlo procedures;  our proposed empirical Edgeworth expansion; and benchmarks.  We make two main observations.  First, \texttt{TrueA} and \texttt{TrueAJack} are almost indistinguishable, echoing our Theorem \ref{thm::jackknife}; meanwhile, they are both smooth and rather different from the step-function \texttt{TrueW}.  This clearly demonstrates the self-smoothing feature of $\hat{T}_n$ in the lattice case.  If we change the graphon to a smooth one, these curves would all be smooth and close to each other.  Second, we observe the higher accuracy of our empirical Edgeworth expansion compared to competing methods.  In fact, repeating this experiment multiple times,  our method shows significantly stabler approximations than bootstraps.

\begin{figure}[htbp!]
	\begin{adjustwidth}{-\oddsidemargin-1.5in}{-\rightmargin-2in}
		\centering
	\includegraphics[width=0.4\textwidth]{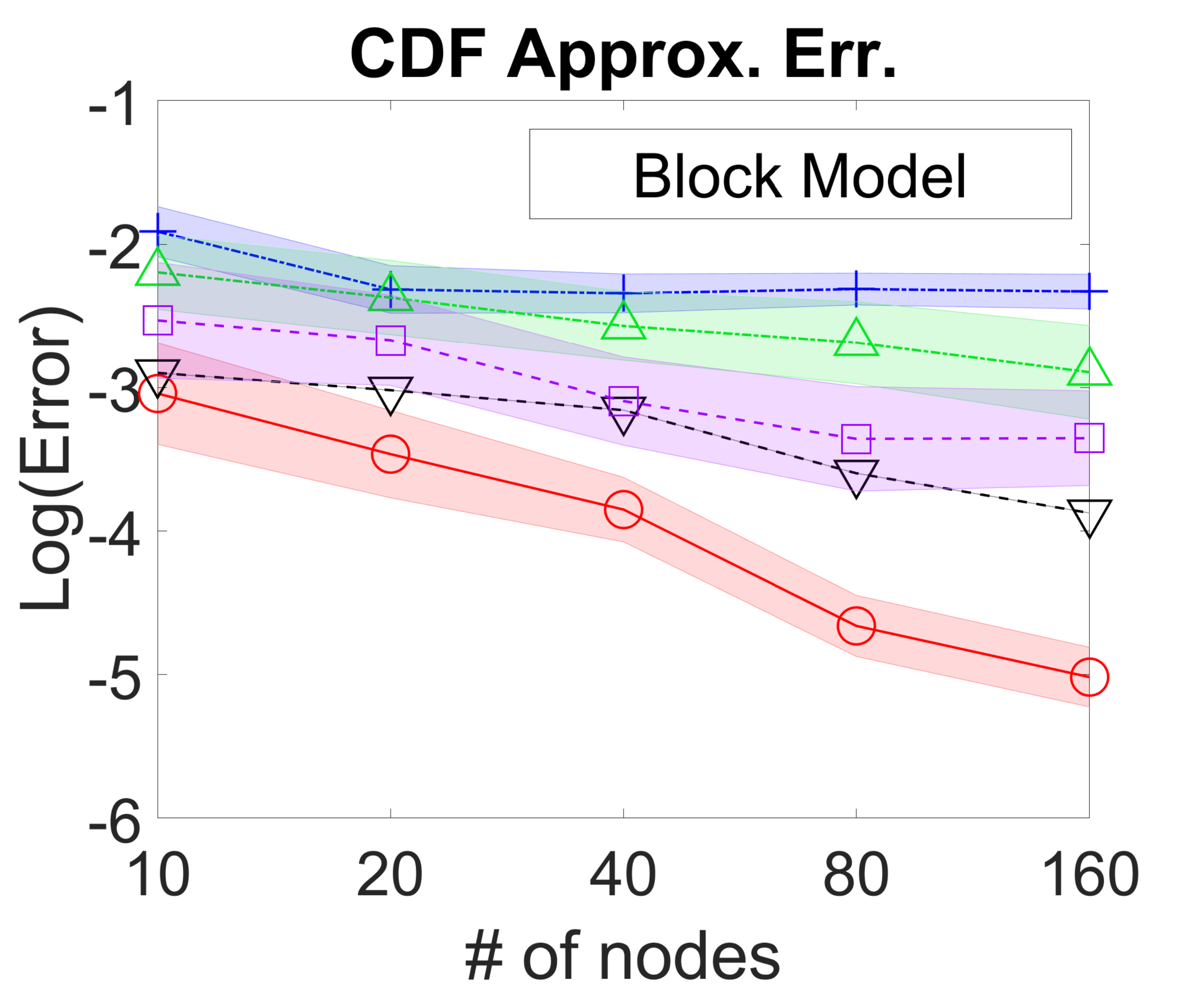}
	\includegraphics[width=0.4\textwidth]{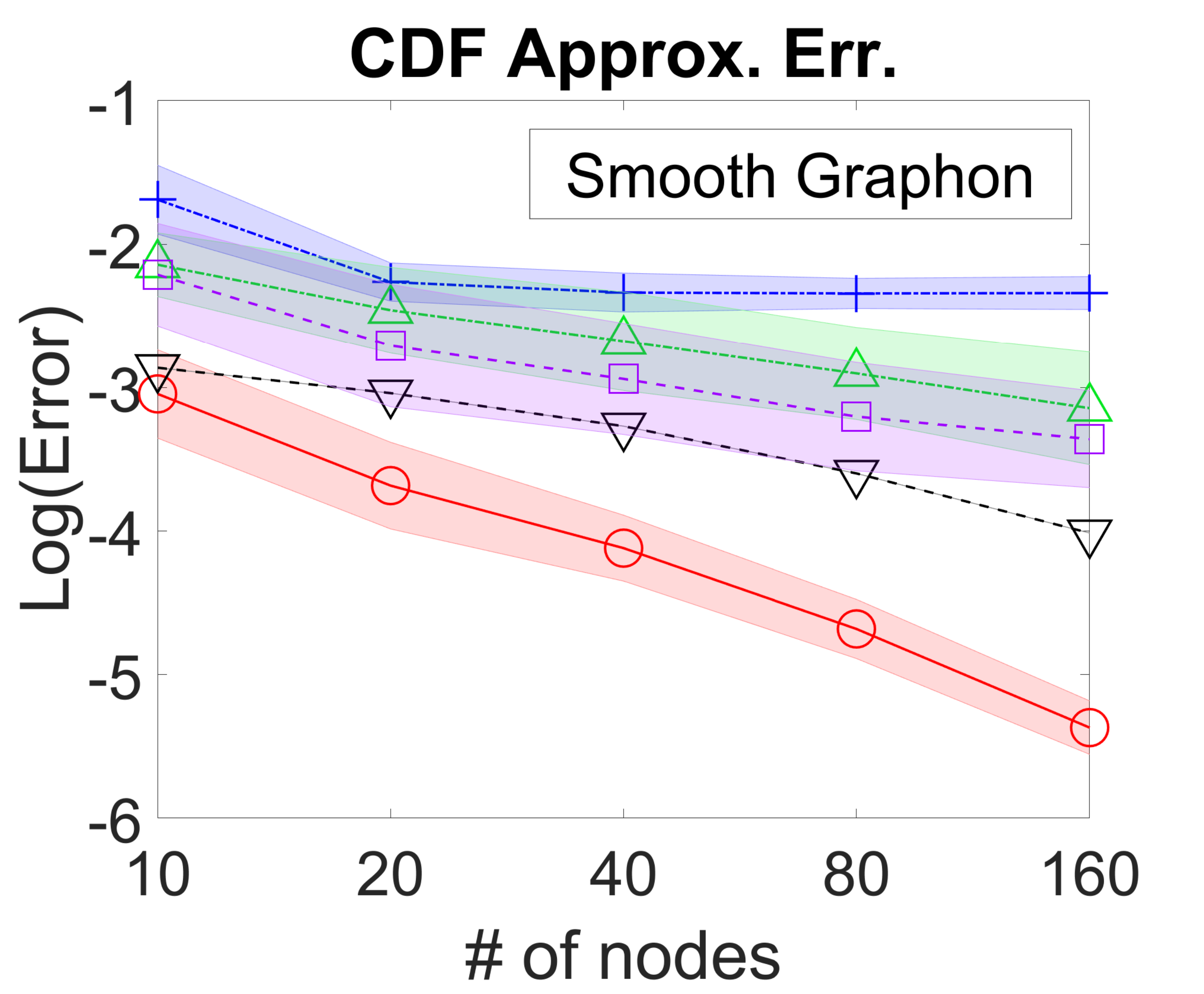}
	\includegraphics[width=0.4\textwidth]{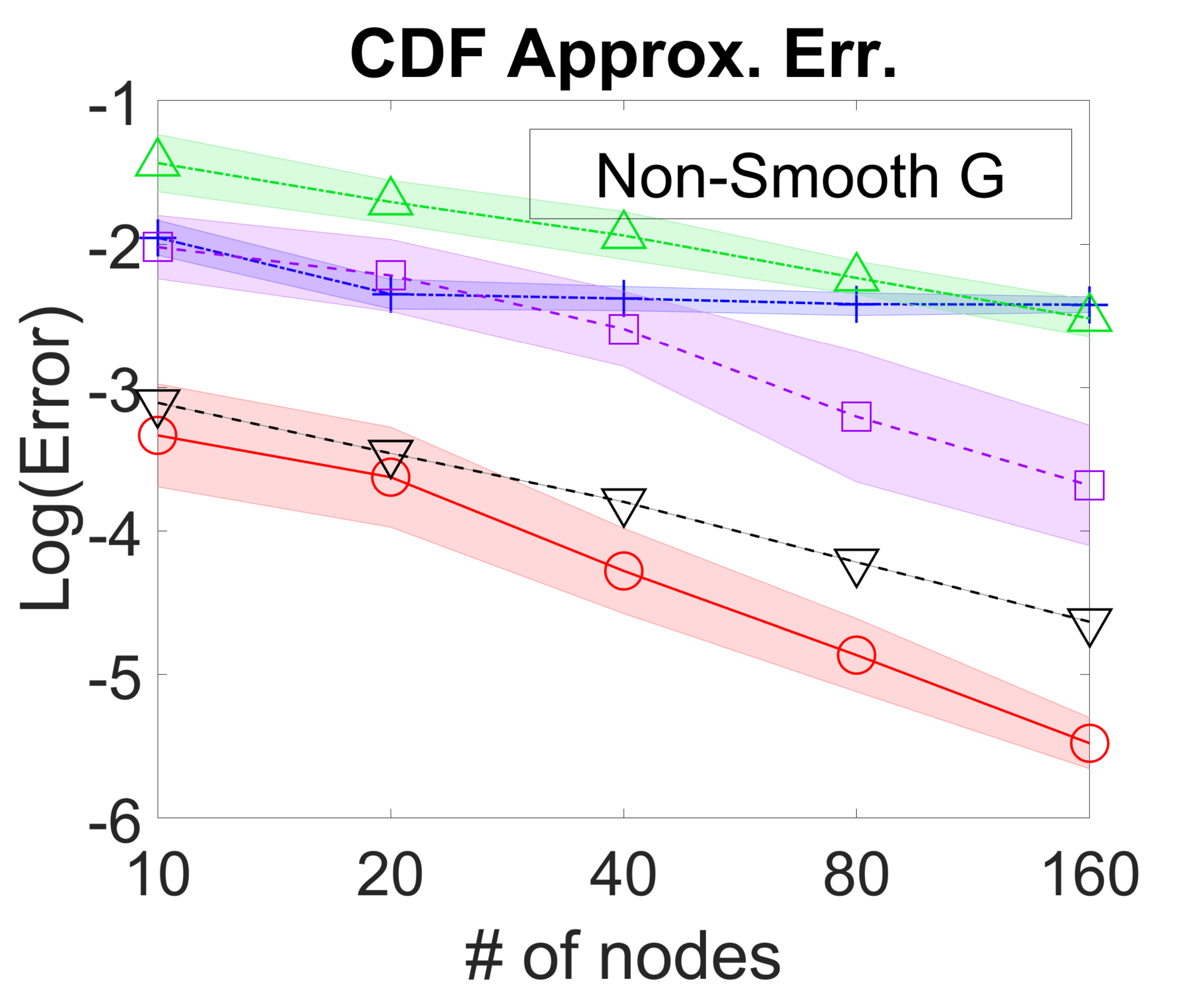}\\
	\includegraphics[width=0.4\textwidth]{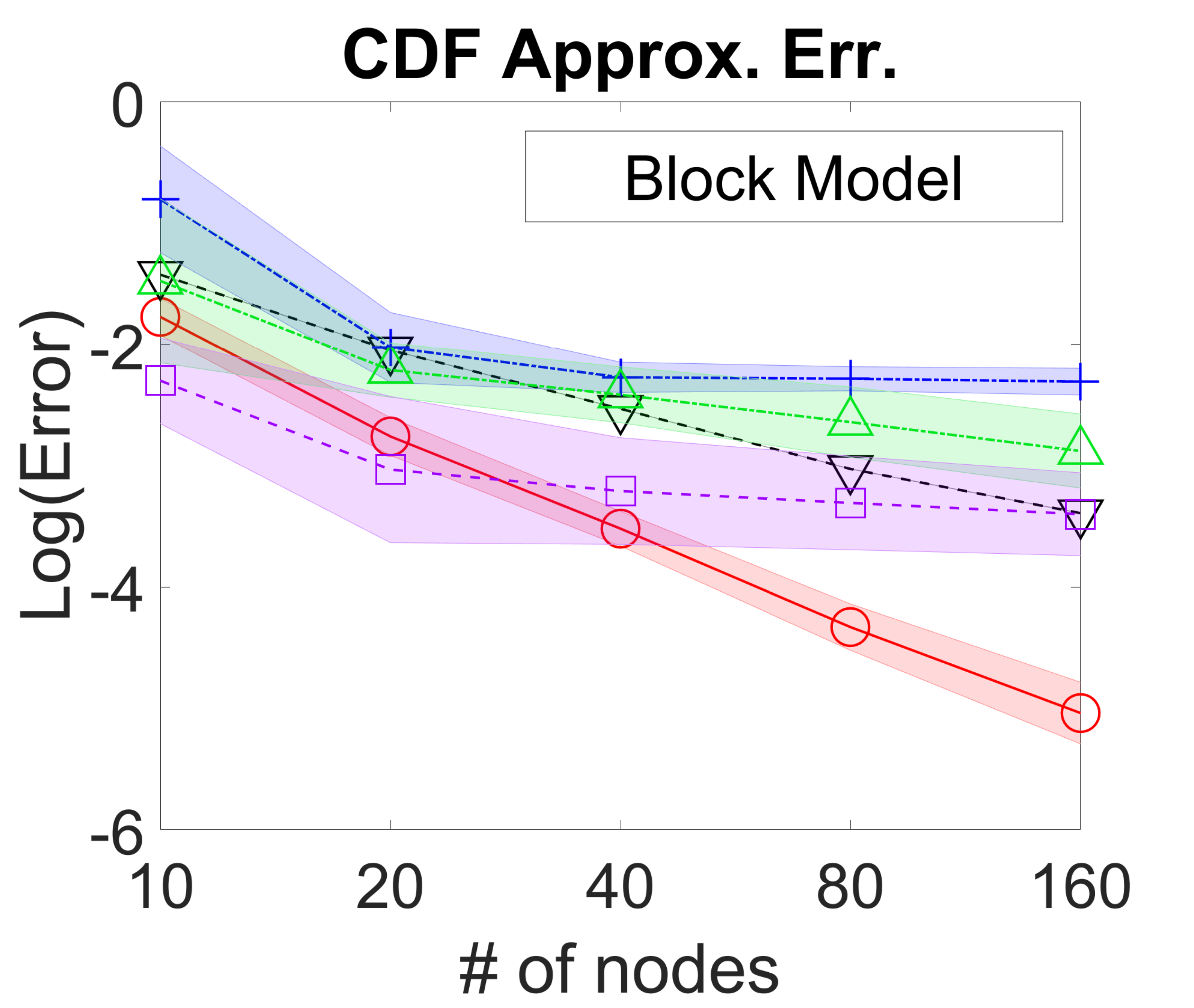}
	\includegraphics[width=0.4\textwidth]{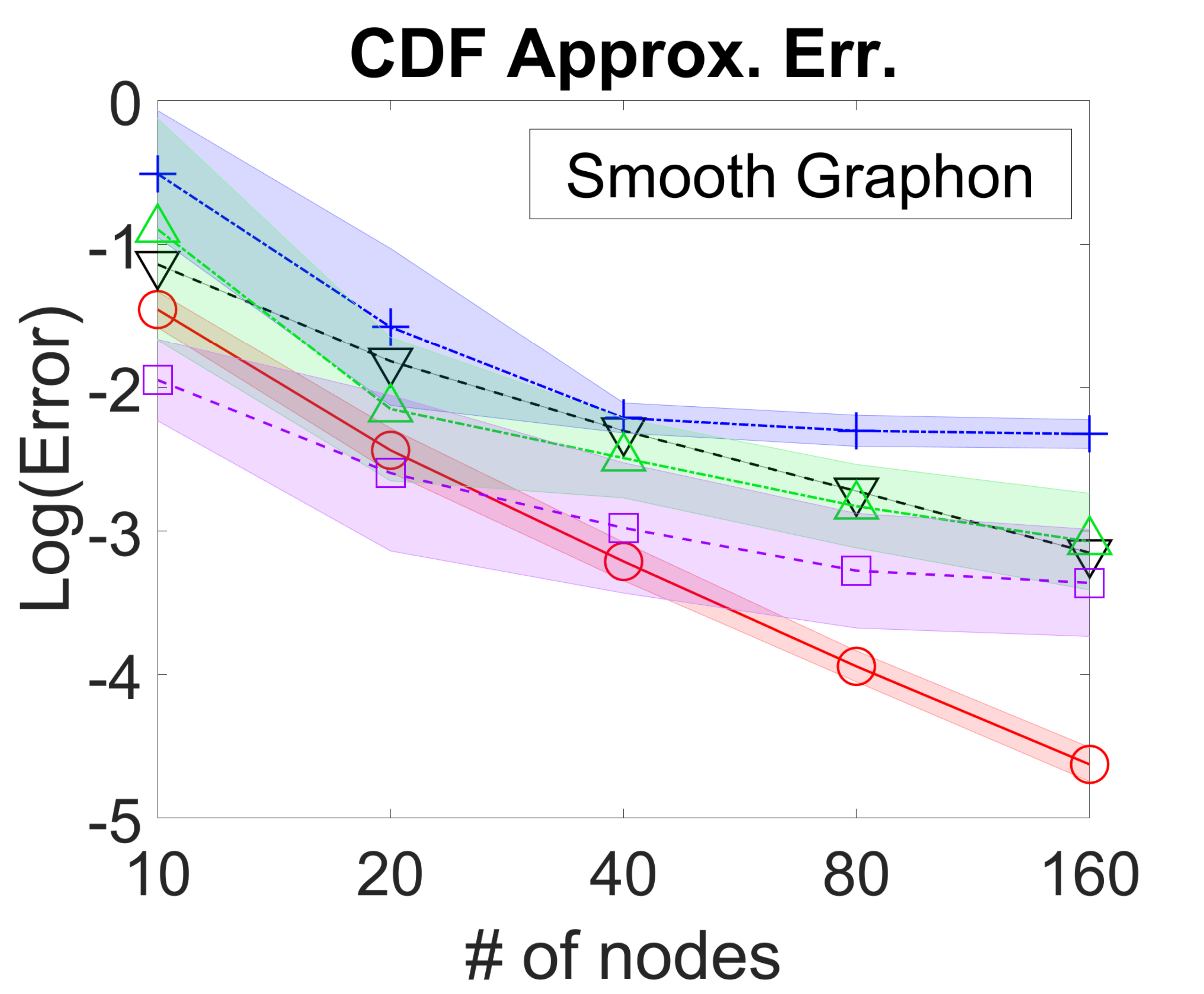}
	\includegraphics[width=0.4\textwidth]{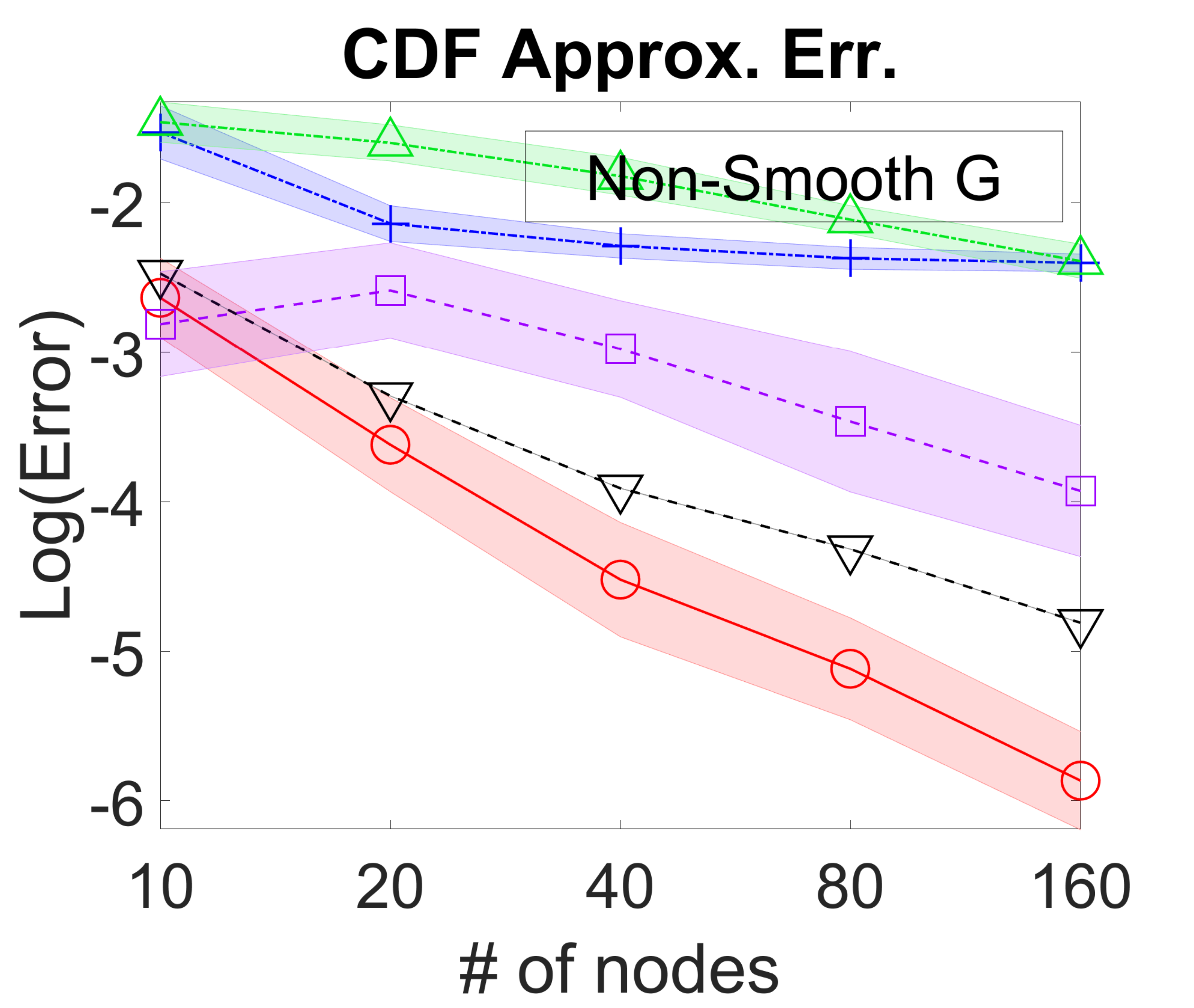}\\
	\includegraphics[width=0.4\textwidth]{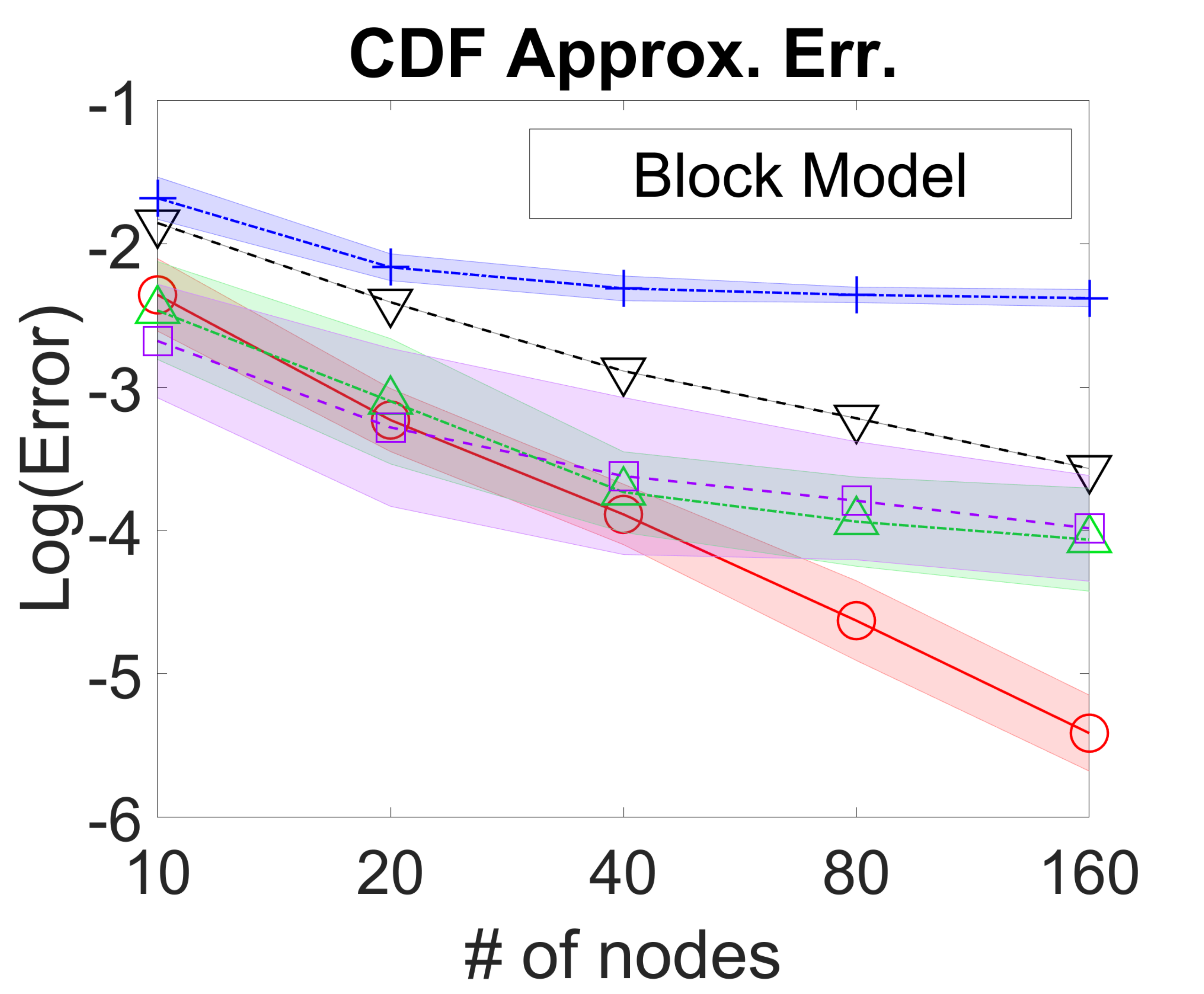}
	\includegraphics[width=0.4\textwidth]{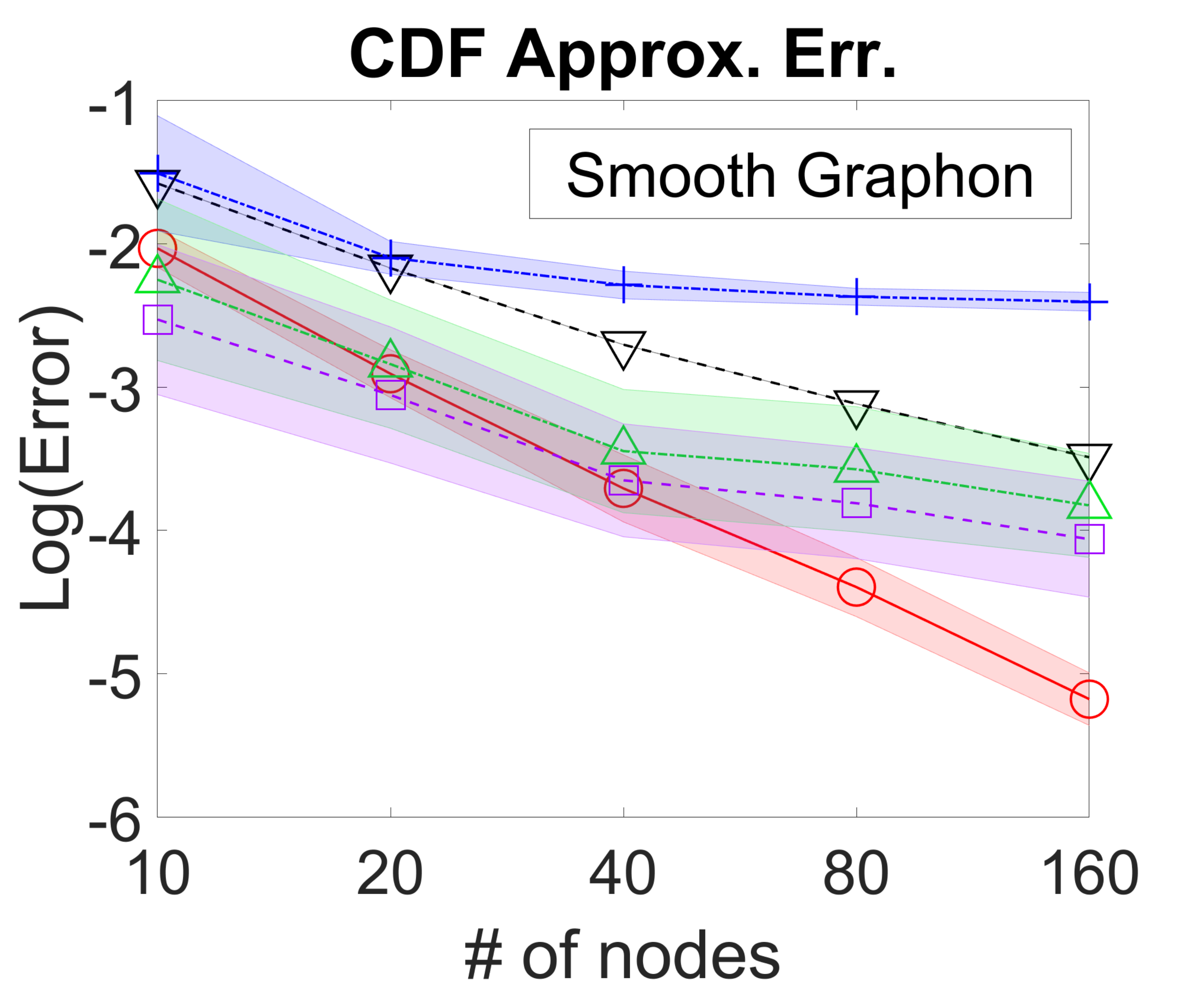}
	\includegraphics[width=0.4\textwidth]{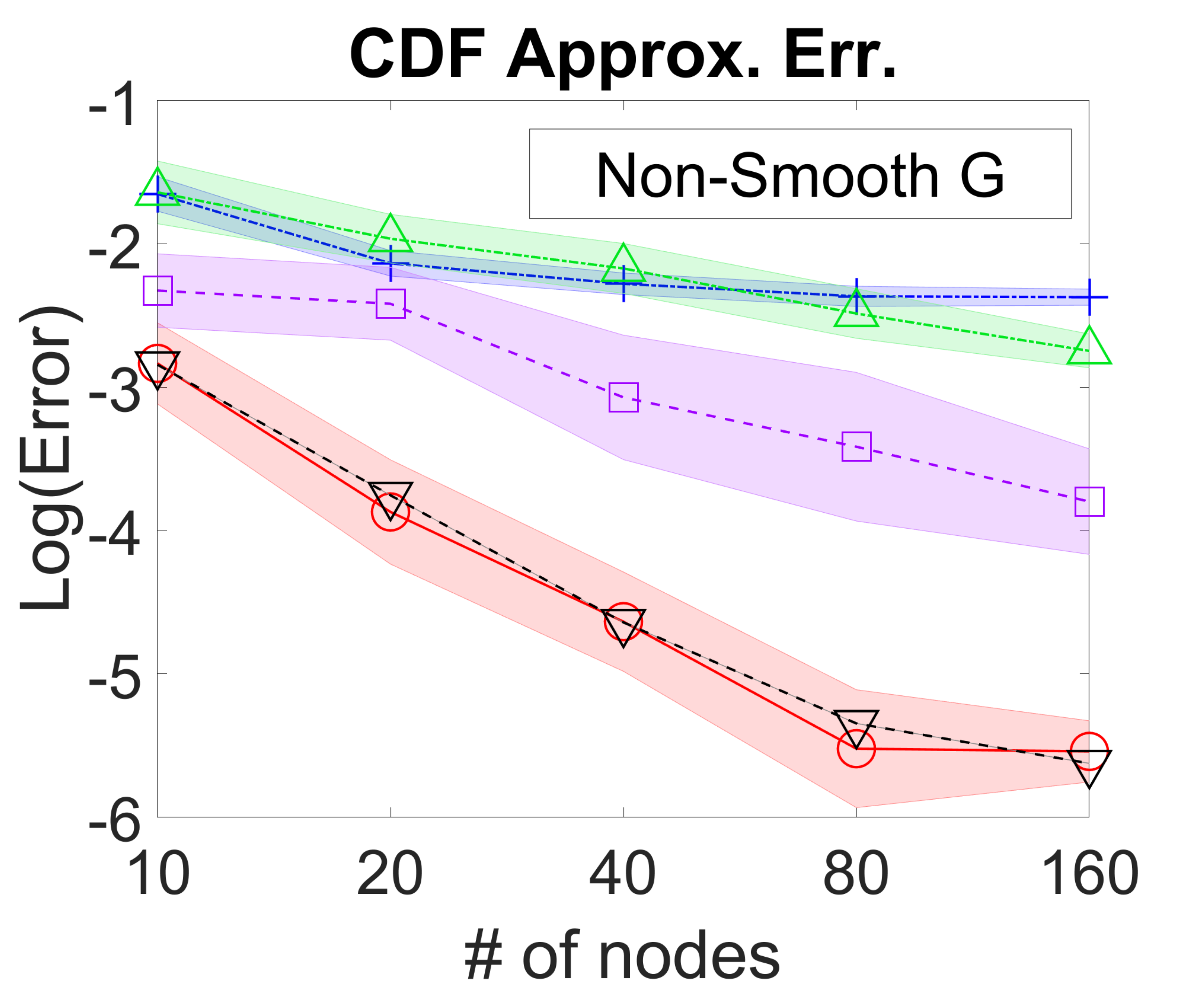}\\
	\includegraphics[width=0.4\textwidth]{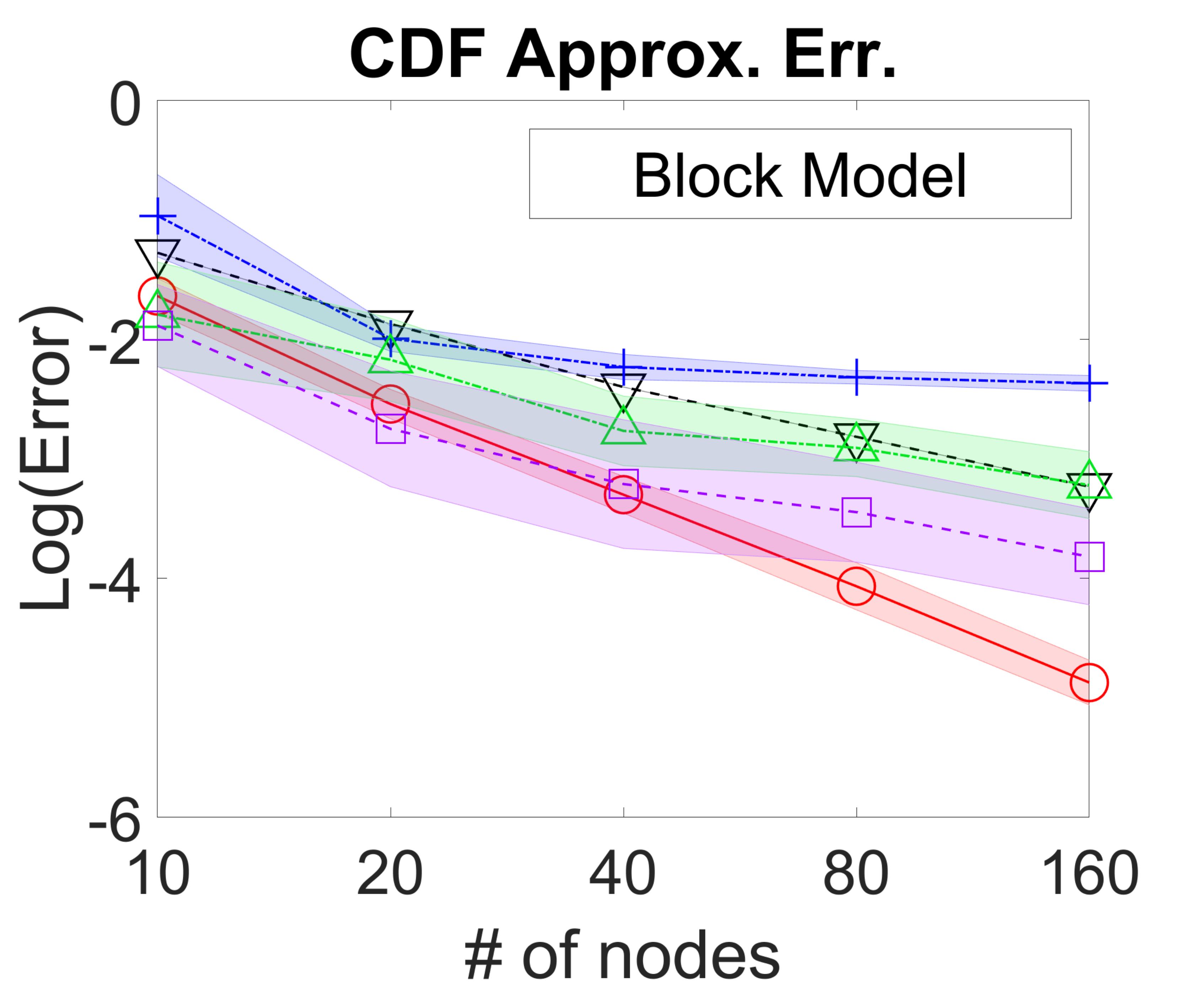}
	\includegraphics[width=0.4\textwidth]{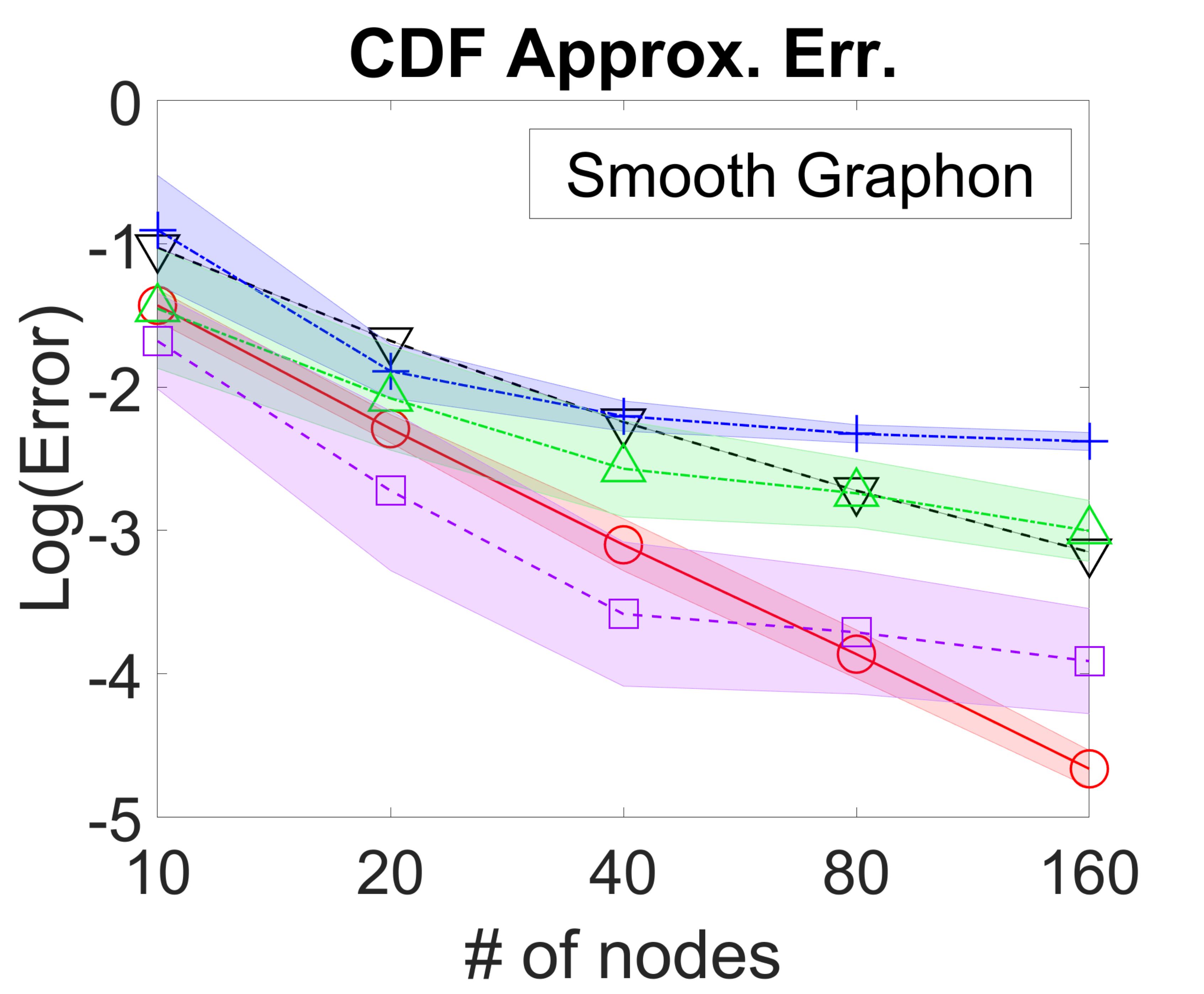}
	\includegraphics[width=0.4\textwidth]{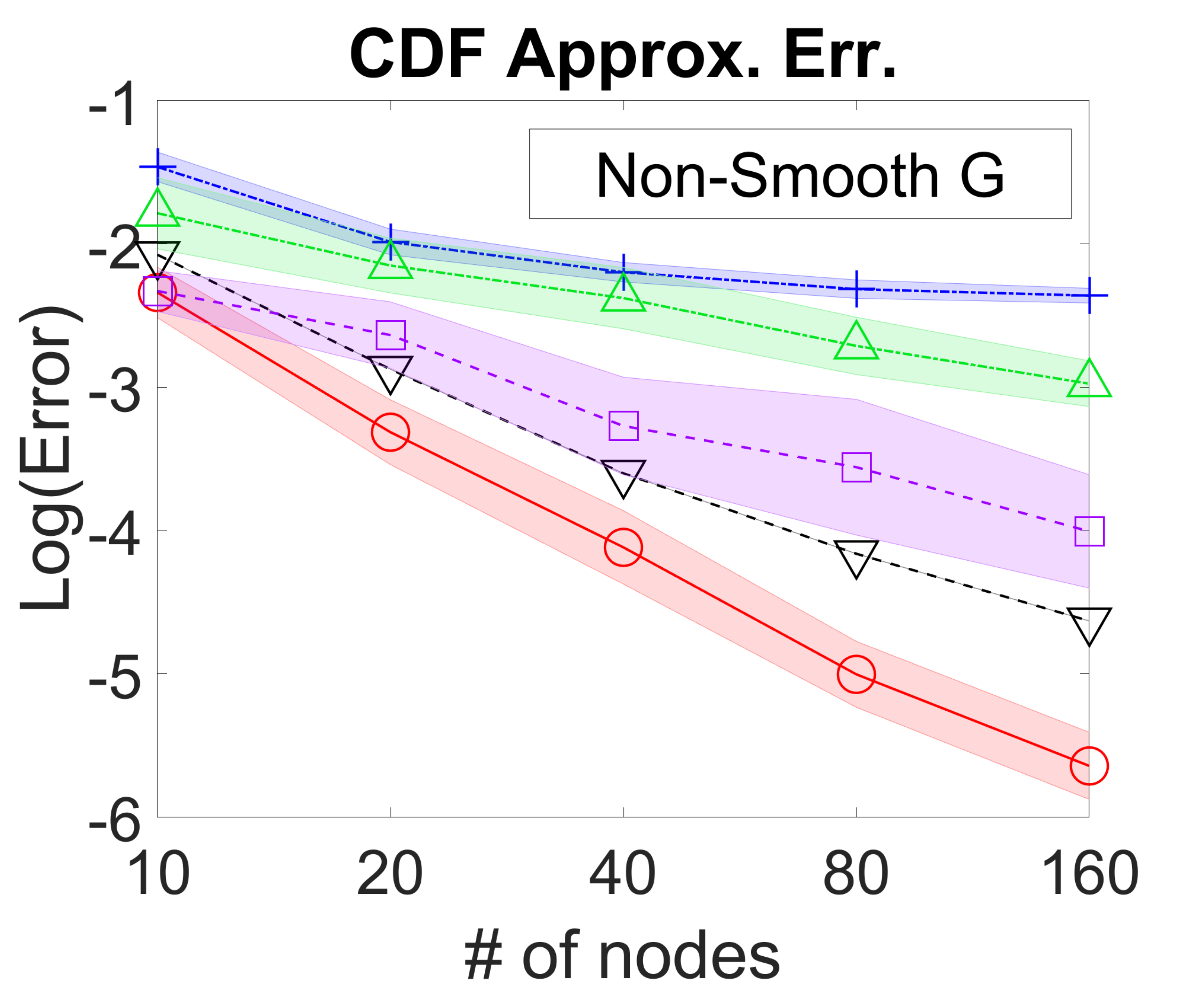}
	\end{adjustwidth}
	\caption{{CDF approximation errors.  Both axes are log(e)-scaled.} {\bf Motifs:} row 1: {\tt Edge}; row 2: {\tt Triangle}; row 3: {\tt Vshape}{; row 4: {\tt ThreeStar}}.    \tred{Red solid curve marked circle}: our method (empirical Edgeworth); black dashed curve marked down-triangle: $N(0,1)$ approximation; \textcolor{green}{green dashed curve marked up-triangle}: re-sampling of $A$ in \cite{green2017bootstrapping}; \tblue{blue dashed curve marked plus}: \cite{bhattacharyya2015subsampling} sub-sampling $\asymp n$ nodes; \textcolor{magenta}{magenta dashed line with square markers}: ASE plug-in bootstrap in \cite{levin2019bootstrapping}.}
	\label{fig::numerical-1}
\end{figure}

\begin{figure}[htbp!]
	\begin{adjustwidth}{-\oddsidemargin-1.5in}{-\rightmargin-2in}
		\centering
		\includegraphics[width=0.4\textwidth]{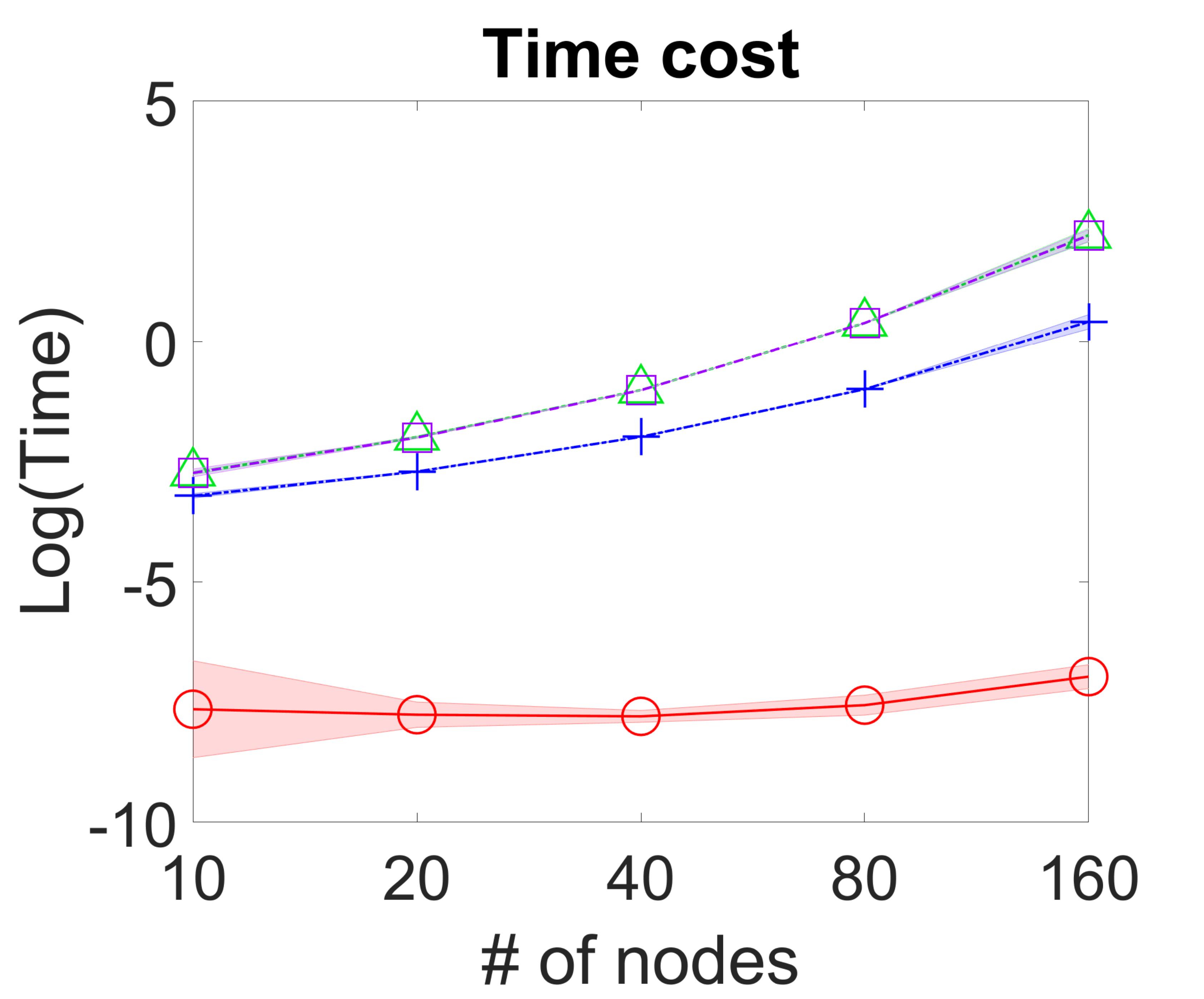}
		\includegraphics[width=0.4\textwidth]{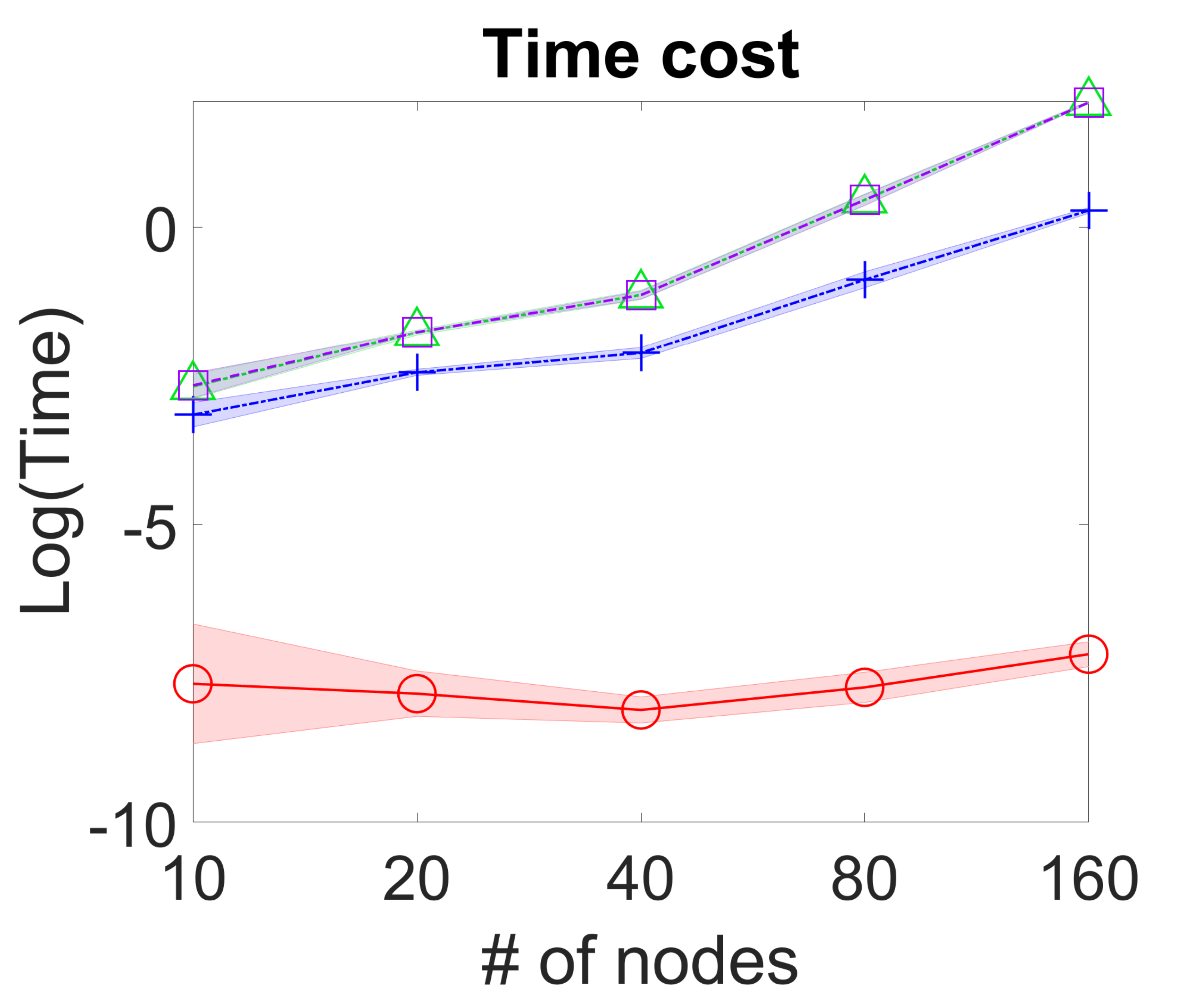}
		\includegraphics[width=0.4\textwidth]{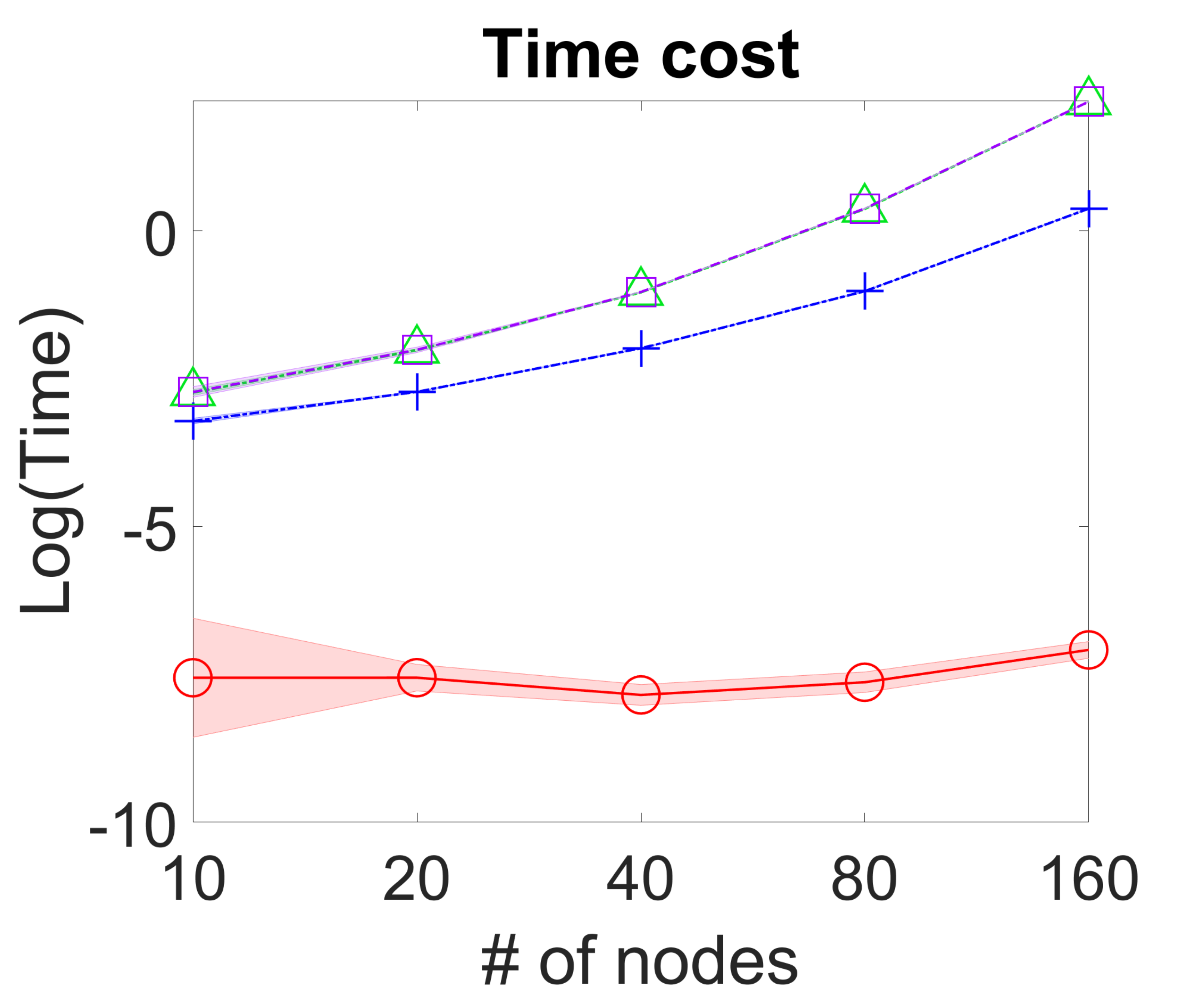}\\
		\includegraphics[width=0.4\textwidth]{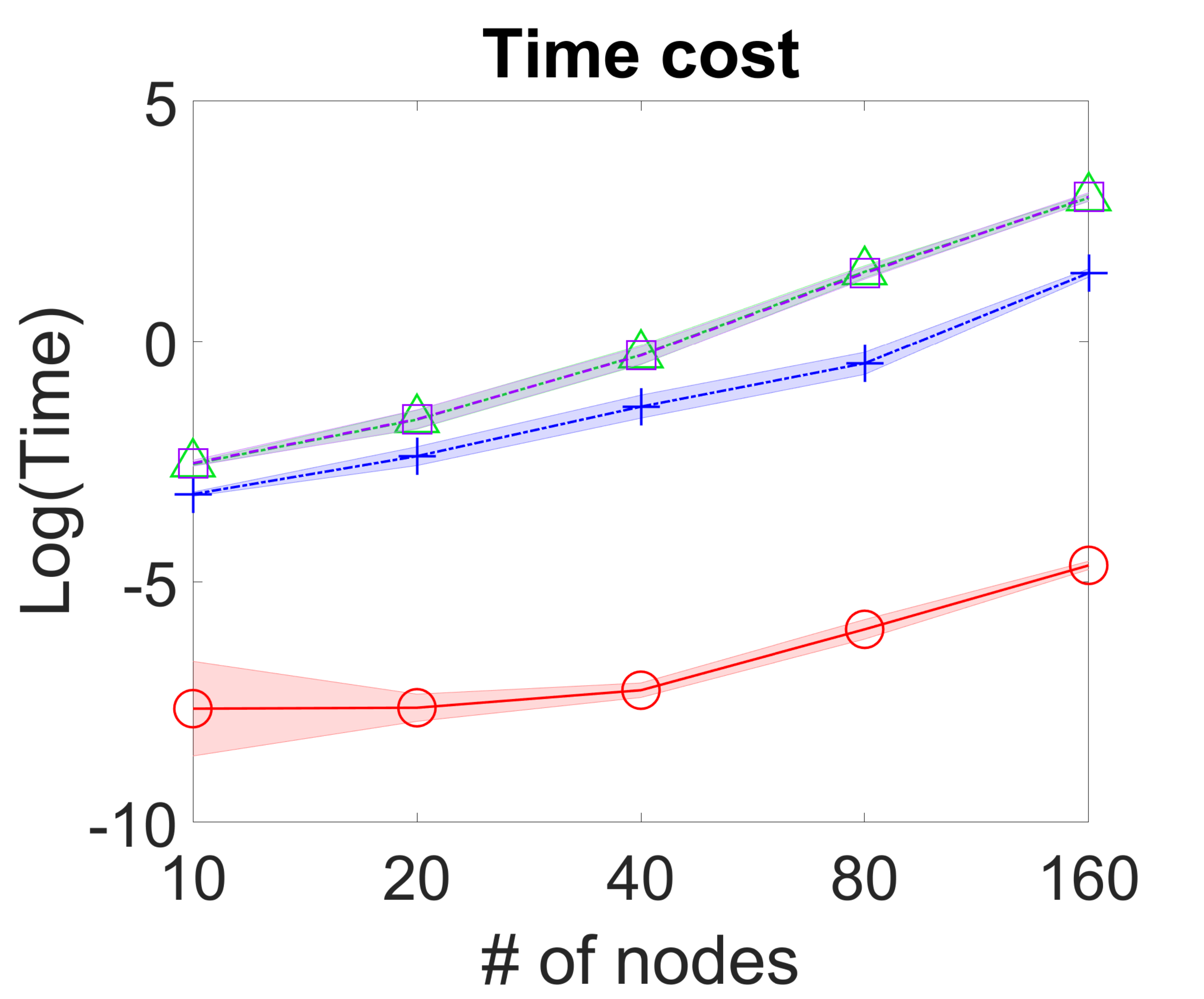}
		\includegraphics[width=0.4\textwidth]{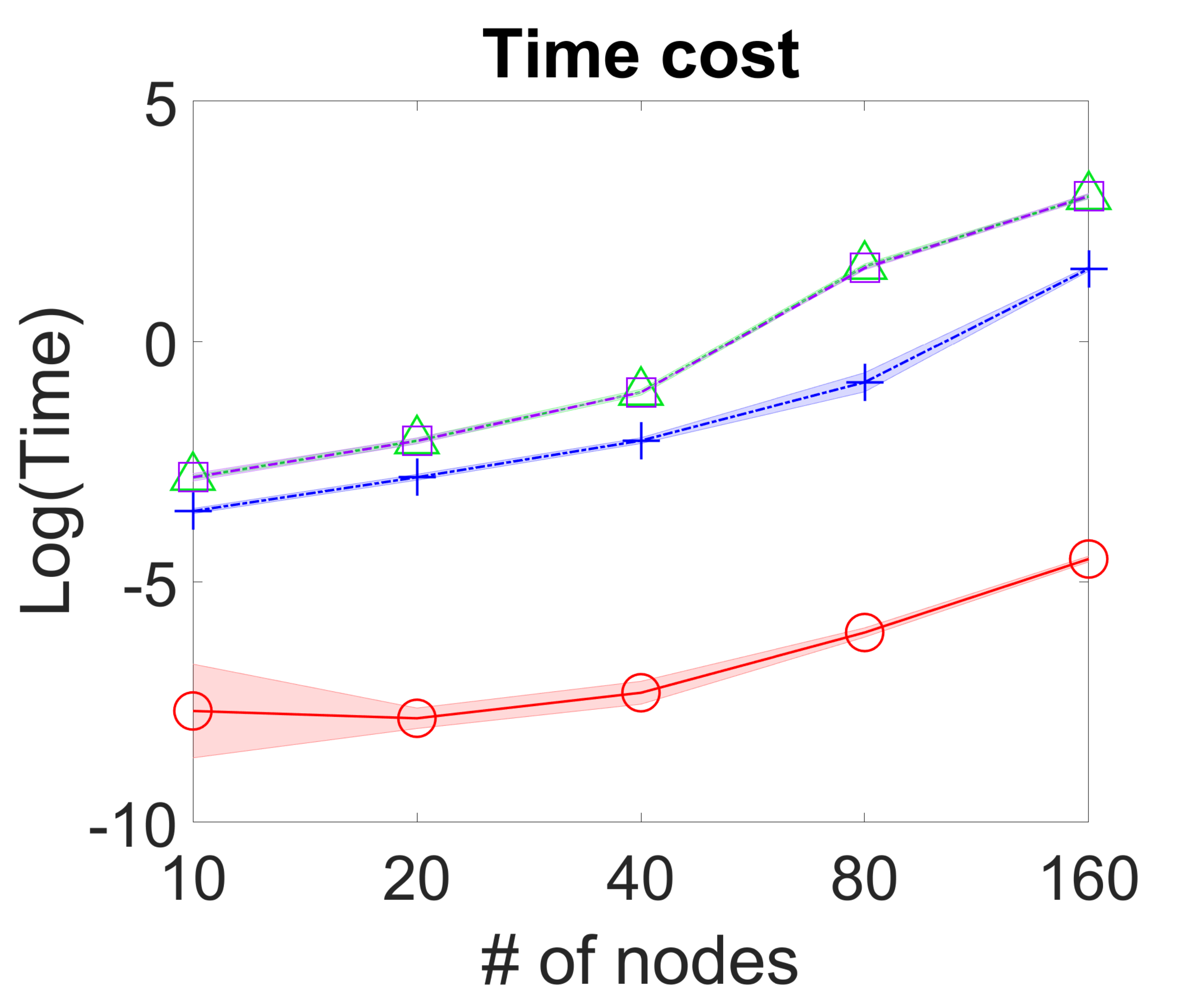}
		\includegraphics[width=0.4\textwidth]{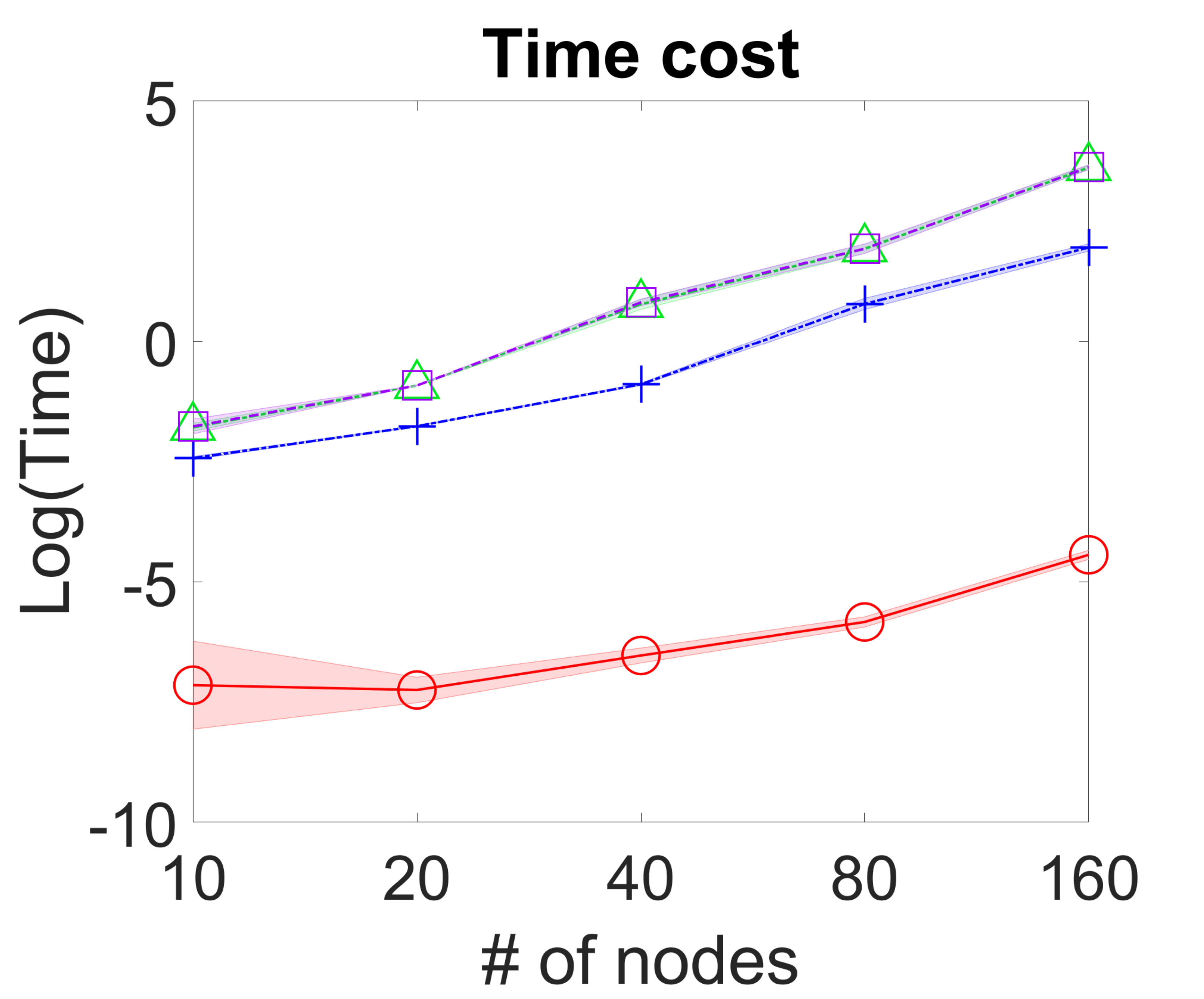}\\
		\includegraphics[width=0.4\textwidth]{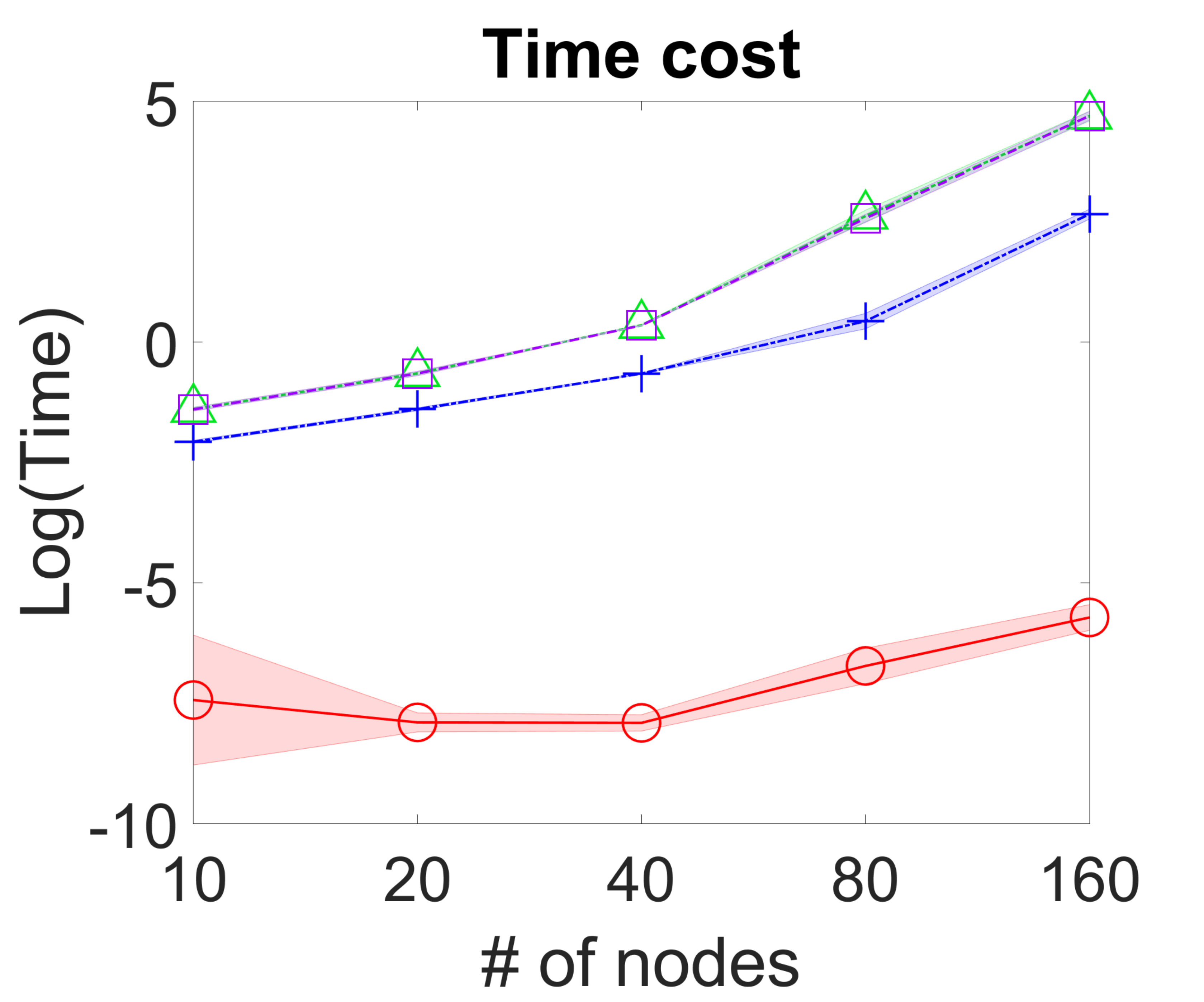}
		\includegraphics[width=0.4\textwidth]{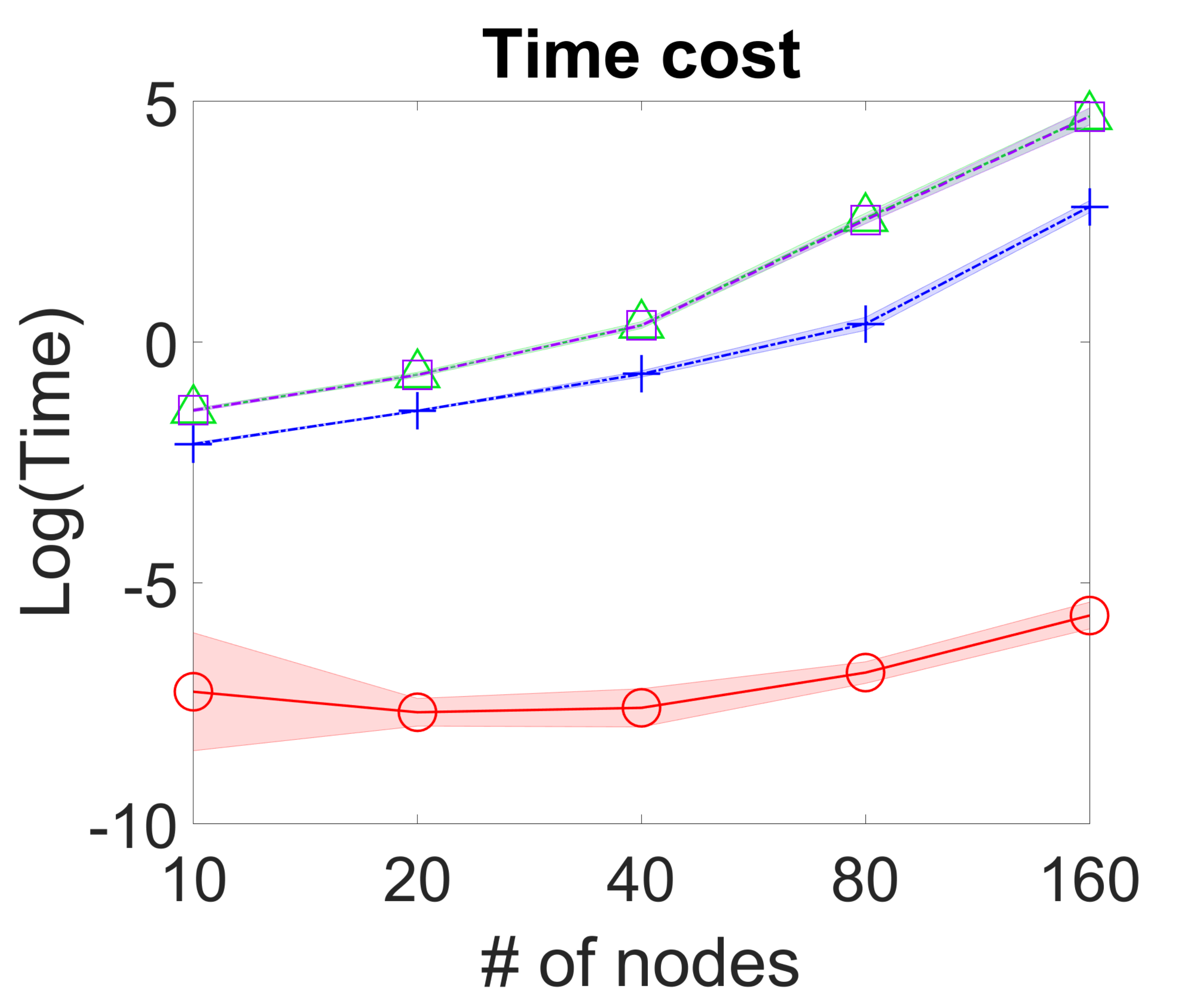}
		\includegraphics[width=0.4\textwidth]{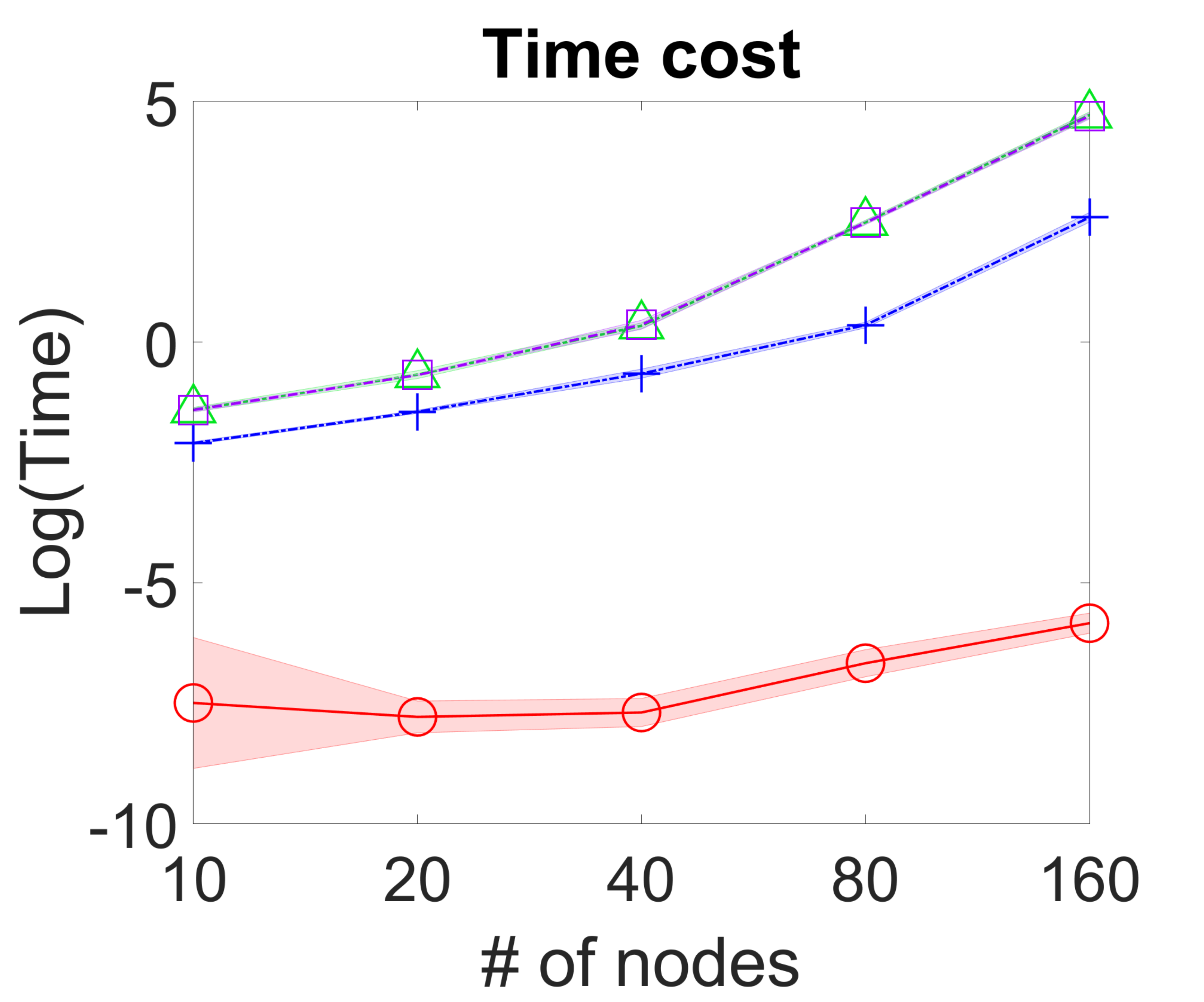}\\
		\includegraphics[width=0.4\textwidth]{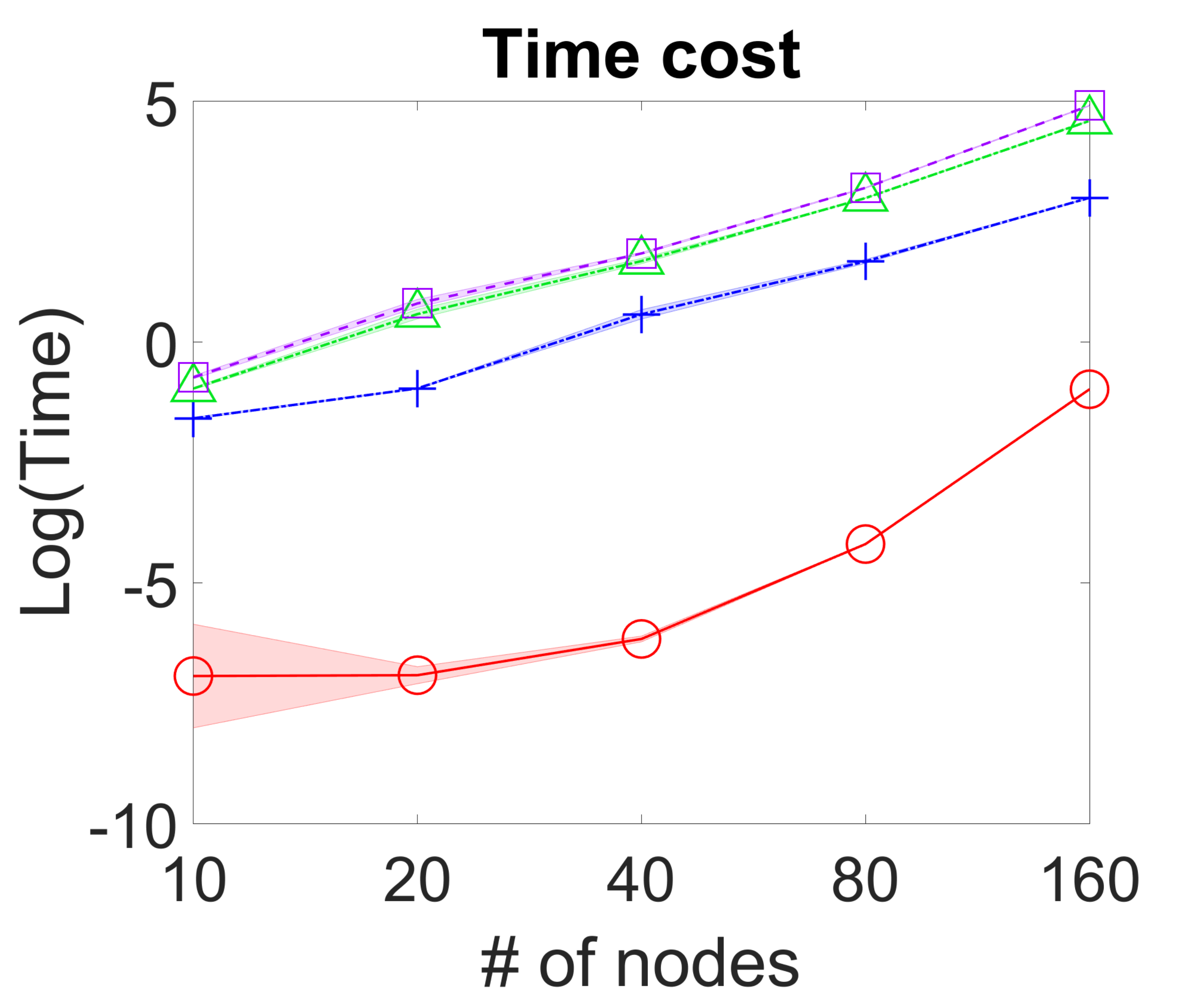}
	\includegraphics[width=0.4\textwidth]{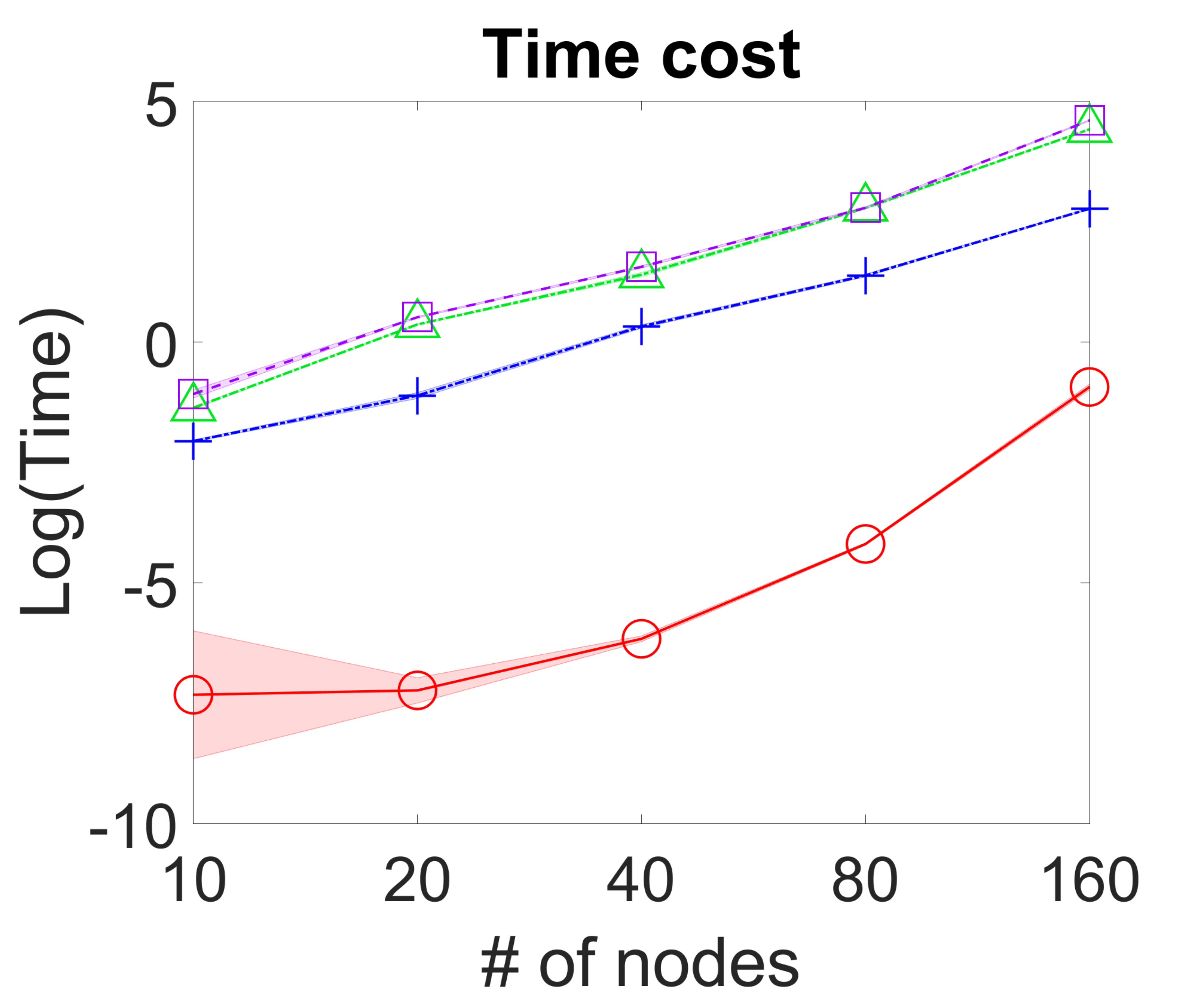}
	\includegraphics[width=0.4\textwidth]{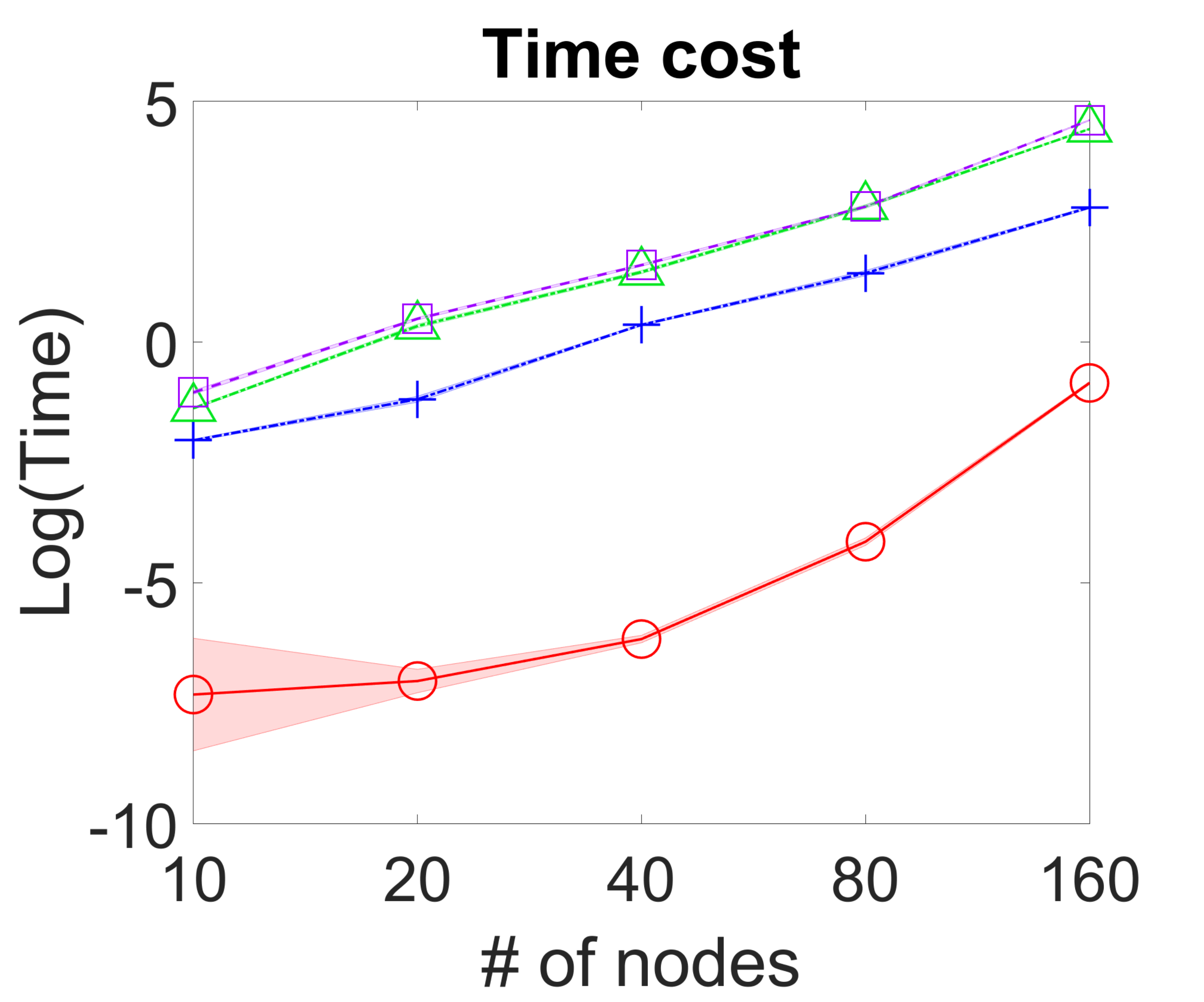}
	\end{adjustwidth}
	\caption{{Time costs (in seconds) of all methods.}  Both axes are log(e)-scaled.  {\bf Motifs:} row 1: {\tt Edge}; row 2: {\tt Triangle}; row 3: {\tt Vshape}{; row 4: {\tt ThreeStar}}.      \tred{Red solid curve marked circle}: our method (empirical Edgeworth); \textcolor{green}{green dashed curve marked up-triangle}: re-sampling of $A$ in \cite{green2017bootstrapping}; \tblue{blue dashed curve marked plus}: \cite{bhattacharyya2015subsampling} sub-sampling $\asymp n$ nodes; \textcolor{magenta}{magenta dashed line with square markers}: ASE plug-in bootstrap in \cite{levin2019bootstrapping}.  We regarded $N(0,1)$ as zero time cost so does not appear in the time cost plot.}
	\label{fig::numerical-2}
\end{figure}

Next, we conduct a systematic comparison of the performances of all methods across many settings.  We mainly varied three factors: graphon type, motif type and network size, over the previously described ranges.  Our experiment results under different network sparsity levels would have to sink to Supplemental Material due to page limit, and here we keep $\rho_n=1$.  Results are shown in Figure \ref{fig::numerical-1} (error) and Figure \ref{fig::numerical-2} (time cost), where error bars show standard deviations.

In all experiments, our empirical Edgeworth expansion approach exhibited clear advantages over benchmark methods in all aspects: the absolute values of errors, the diminishing rates of errors, and computational efficiency.  Our method shows a higher-order accuracy by slopes steeper than $-1/2$ and much steeper than other methods.  On computation efficiency, our method is the second cheapest after the simple $N(0,1)$ approximation (that does not need computation) and much faster than network bootstraps.  It typically costs about $e^{-5}\approx 1/150$ the time of sub-sampling and about $e^{-7}\approx 1/1000$ the time of re-sampling.  Our method only needs one run and does not require repeated sampling.

Notice that there is no simple rule to judge the difficulty of different scenarios, which jointly depends on the graphon and the motif through implicit and complex relationship.  In our experience, triangle may be more difficult than V-shape under some graphons, but easier under some others, and this comparison may vary from method to method.  Answering this question requires calculation of the population Edgeworth expansion up to $o(n^{-1})$ remainder, and the leading term in the remainder of the one-term Edgeworth expansion would then quantify the real difficulty.  But the calculation is very complicated and outside the scope of this paper.

We did not observe the higher-order accuracy of bootstrap methods as our results predicted.  One likely reason is the numerical accuracy limited by the $n_{\mathrm{boot}}$ that our computing servers can afford.  We did see an observable improvement in the performances of network bootstraps as we increased $n_{\mathrm{boot}}$ from 200 suggested by \cite{levin2019bootstrapping} to the current 2000.  But further increasing $n_{\mathrm{boot}}$ will also increase their time costs and potentially memory usage.  We ran each experiment on 36 parallel Intel(R) Xeon(R) X5650 CPU cores at \@ 2.67GHz with 12M cache and 2GB RAM.  It took roughly 3$\sim$8 hours to run each experiment that produces one individual plot in Figures \ref{fig::numerical-1} and \ref{fig::numerical-2}.

{

\subsection{Simulation 2:  Finite-sample performance of Cornish-Fisher confidence interval}
\label{subsec::simulation-2::CI}

In this simulation, we numerically assess the performance of our Cornish-Fisher confidence interval, compared to benchmark methods.  Throughout this subsection, we set $\alpha=0.2$ and focus on symmetric two-sided confidence intervals.
We inherit most simulation settings from Section \ref{subsec::simulation-1} with some modifications we now clarify.  The main difference is that in this simulation, we must conduct many repeated experiments in order to accurately evaluate the coverage probability (each iteration produces a binary outcome of whether the CI contains the population parameter).  We repeated the experiment 10000 times for our method and normal approximation, and 500 times for the much slower bootstrap methods.  Due to the computer limitations, while we can keep the same number of Monte Carlo evaluations, in order to repeat the entire experiment 500 times to accurately evaluate the actual CI coverage rates of bootstraps, we have to reduce their numbers of bootstrap samples to 500 (still more than the 200 in \cite{levin2019bootstrapping}).  
We evaluate three performance measures: {\tt coverage}: \emph{actual coverage probability}; {\tt length}: \emph{confidence intereval length}; and {\tt time}: \emph{computation time in seconds}.

\begin{table}[htbp!]
\centering
\caption{Performance measures of 95\% confidence intervals\\$n=80$, $\rho_n\asymp 1$, graphon: block model}
\label{table::n=80::sparse-dense::graphon-block model}
\begin{adjustbox}{center}
\begin{tabular}{c|cccc}\hline
Method & Edge & Triangle & V-shape & Three star \\\hline
Our method & \bigcell{c}{Coverage~$=0.957(0.202)$\\Length~$=0.097(0.010)$\\LogTime~$=-8.092(0.141)$} & \bigcell{c}{$0.953(0.211)$\\$0.040(0.008)$\\$-7.292(0.094)$} & \bigcell{c}{$0.956(0.205)$\\$0.200(0.033)$\\$-7.393(0.135)$} & \bigcell{c}{$0.952(0.213)$\\$0.145(0.033)$\\$-4.184(0.023)$}\\\hline
Norm. Approx. & \bigcell{c}{$0.950(0.218)$\\$0.097(0.010)$\\No time cost} & \bigcell{c}{$0.934(0.248)$\\$0.040(0.008)$\\No time cost} & \bigcell{c}{$0.942(0.235)$\\$0.200(0.033)$\\No time cost} & \bigcell{c}{$0.932(0.251)$\\$0.145(0.033)$\\No time cost}\\\hline
\citet{bhattacharyya2015subsampling} & \bigcell{c}{$0.814(0.389)$\\$0.069(0.009)$\\$-1.444(0.045)$} & \bigcell{c}{$0.832(0.374)$\\$0.032(0.006)$\\$-1.129(0.004)$} & \bigcell{c}{$0.826(0.379)$\\$0.148(0.025)$\\$-1.026(0.117)$} & \bigcell{c}{$0.828(0.378)$\\$0.114(0.025)$\\$0.348(0.016)$}\\\hline
\citet{green2017bootstrapping} & \bigcell{c}{$0.946(0.226)$\\$0.097(0.012)$\\$-0.052(0.027)$} & \bigcell{c}{$0.950(0.218)$\\$0.045(0.009)$\\$1.306(0.096)$} & \bigcell{c}{$0.944(0.230)$\\$0.206(0.036)$\\$1.119(0.071)$} & \bigcell{c}{$0.950(0.218)$\\$0.153(0.036)$\\$1.543(0.011)$}\\\hline
\citet{levin2019bootstrapping} & \bigcell{c}{$0.960(0.196)$\\$0.100(0.012)$\\$-0.052(0.028)$} & \bigcell{c}{$0.958(0.201)$\\$0.044(0.009)$\\$1.273(0.099)$} & \bigcell{c}{$0.960(0.196)$\\$0.211(0.036)$\\$1.115(0.069)$} & \bigcell{c}{$0.956(0.205)$\\$0.157(0.036)$\\$1.728(0.010)$}\\\hline
\end{tabular}
\end{adjustbox}
\end{table}

\begin{table}[htbp!]
\centering
\caption{Performance measures of 95\% confidence intervals\\$n=80$, $\rho_n\asymp 1$, graphon: smooth graphon}
\label{table::n=80::sparse-dense::graphon-smooth graphon}
\begin{adjustbox}{center}
\begin{tabular}{c|cccc}\hline
Method & Edge & Triangle & V-shape & Three star \\\hline
Our method & \bigcell{c}{Coverage~$=0.958(0.201)$\\Length~$=0.092(0.009)$\\LogTime~$=-8.054(0.052)$} & \bigcell{c}{$0.940(0.238)$\\$0.021(0.005)$\\$-7.235(0.110)$} & \bigcell{c}{$0.951(0.216)$\\$0.141(0.025)$\\$-7.466(0.106)$} & \bigcell{c}{$0.942(0.235)$\\$0.083(0.021)$\\$-4.137(0.022)$}\\\hline
Norm. Approx. & \bigcell{c}{$0.951(0.216)$\\$0.092(0.009)$\\No time cost} & \bigcell{c}{$0.920(0.271)$\\$0.021(0.005)$\\No time cost} & \bigcell{c}{$0.938(0.242)$\\$0.141(0.025)$\\No time cost} & \bigcell{c}{$0.923(0.266)$\\$0.083(0.021)$\\No time cost}\\\hline
\citet{bhattacharyya2015subsampling} & \bigcell{c}{$0.870(0.337)$\\$0.066(0.008)$\\$-1.436(0.007)$} & \bigcell{c}{$0.868(0.339)$\\$0.018(0.005)$\\$-0.690(0.007)$} & \bigcell{c}{$0.876(0.330)$\\$0.111(0.020)$\\$-1.004(0.036)$} & \bigcell{c}{$0.870(0.337)$\\$0.073(0.017)$\\$0.062(0.136)$}\\\hline
\citet{green2017bootstrapping} & \bigcell{c}{$0.954(0.210)$\\$0.092(0.011)$\\$0.004(0.006)$} & \bigcell{c}{$0.952(0.214)$\\$0.025(0.006)$\\$1.206(0.014)$} & \bigcell{c}{$0.956(0.205)$\\$0.148(0.027)$\\$1.188(0.069)$} & \bigcell{c}{$0.956(0.205)$\\$0.090(0.023)$\\$1.432(0.087)$}\\\hline
\citet{levin2019bootstrapping} & \bigcell{c}{$0.962(0.191)$\\$0.096(0.010)$\\$-0.021(0.005)$} & \bigcell{c}{$0.962(0.191)$\\$0.024(0.006)$\\$1.213(0.011)$} & \bigcell{c}{$0.962(0.191)$\\$0.154(0.027)$\\$1.192(0.092)$} & \bigcell{c}{$0.962(0.191)$\\$0.095(0.023)$\\$1.527(0.079)$}\\\hline
\end{tabular}
\end{adjustbox}
\end{table}

\begin{table}[htbp!]
\centering
\caption{Performance measures of 95\% confidence intervals\\$n=80$, $\rho_n\asymp 1$, graphon: non-smooth graphon}
\label{table::n=80::sparse-dense::graphon-non-smooth graphon}
\begin{adjustbox}{center}
\begin{tabular}{c|cccc}\hline
Method & Edge & Triangle & V-shape & Three star \\\hline
Our method & \bigcell{c}{Coverage~$=0.956(0.205)$\\Length~$=0.116(0.009)$\\LogTime~$=-8.082(0.066)$} & \bigcell{c}{$0.957(0.203)$\\$0.135(0.010)$\\$-7.342(0.159)$} & \bigcell{c}{$0.957(0.202)$\\$0.422(0.027)$\\$-7.383(0.145)$} & \bigcell{c}{$0.957(0.203)$\\$0.531(0.040)$\\$-4.276(0.042)$}\\\hline
Norm. Approx. & \bigcell{c}{$0.952(0.215)$\\$0.116(0.009)$\\No time cost} & \bigcell{c}{$0.949(0.220)$\\$0.135(0.010)$\\No time cost} & \bigcell{c}{$0.951(0.215)$\\$0.422(0.027)$\\No time cost} & \bigcell{c}{$0.950(0.218)$\\$0.531(0.040)$\\No time cost}\\\hline
\citet{bhattacharyya2015subsampling} & \bigcell{c}{$0.870(0.337)$\\$0.081(0.008)$\\$-1.328(0.022)$} & \bigcell{c}{$0.876(0.330)$\\$0.096(0.009)$\\$0.445(0.006)$} & \bigcell{c}{$0.876(0.330)$\\$0.297(0.024)$\\$-1.027(0.081)$} & \bigcell{c}{$0.878(0.328)$\\$0.377(0.033)$\\$0.022(0.024)$}\\\hline
\citet{green2017bootstrapping} & \bigcell{c}{$0.946(0.226)$\\$0.113(0.011)$\\$0.145(0.018)$} & \bigcell{c}{$0.940(0.238)$\\$0.136(0.012)$\\$1.595(0.010)$} & \bigcell{c}{$0.946(0.226)$\\$0.417(0.035)$\\$1.188(0.104)$} & \bigcell{c}{$0.950(0.218)$\\$0.531(0.048)$\\$1.435(0.021)$}\\\hline
\citet{levin2019bootstrapping} & \bigcell{c}{$0.960(0.196)$\\$0.116(0.011)$\\$0.146(0.017)$} & \bigcell{c}{$0.960(0.196)$\\$0.138(0.013)$\\$1.602(0.008)$} & \bigcell{c}{$0.958(0.201)$\\$0.427(0.035)$\\$1.173(0.058)$} & \bigcell{c}{$0.958(0.201)$\\$0.544(0.049)$\\$1.408(0.015)$}\\\hline
\end{tabular}
\end{adjustbox}
\end{table}

Due to page limit, here we present the results for the setting $n=80$ and $\rho_n=1$ in Tables \ref{table::n=80::sparse-dense::graphon-block model} (block model), \ref{table::n=80::sparse-dense::graphon-smooth graphon} (smooth graphon) and \ref{table::n=80::sparse-dense::graphon-non-smooth graphon} (non-smooth graphon).  Each entry is formatted ``mean(standard deviation)''.  We sink the remaining results to Supplemental Materials.  Our method exhibits very accurate actual coverage probabilities consistently close to the nominal confidence level.  Our method is the only method that can always achieve a $\leq 0.010$ coverage error across all settings.  It also produces competitively short confidence interval lengths, again, reflecting the high accuracy of the method.  The comparison of computational efficiency between different methods echoes the qualitative results in Figure \ref{fig::numerical-2} despite slightly different settings and confirms our method's huge speed advantage over all bootstrap methods.  We remark that for the three star motif, we used a double for-loop in evaluating $\hat\ep[g_1(X_1)g_1(X_2)g_2(X_1,X_2)]$ for easiness of derivation at the price of speed.

It is interesting to observe that under the setting of this simulation, our empirical Edgeworth expansion method always produces the same interval length as the normal approximation.  This is not a coincidence in view of \eqref{formula::quantile-for-CI}, \eqref{formula::EEE-confidence-interval} and that $z_{\alpha/2}^2=z_{1-\alpha/2}^2$.  In other words, our two-sided Edgeworth confidence interval is a bias-corrected version (by mean-shift) of the corresponding ordinary CLT confidence interval.

}

\subsection{Simulation 3:  Numerical evaluation of the finite-sample impact of sparsity}
\label{subsec::simulation-3::sparse-networks}

In this part, we conduct numerical studies to evaluate the finite sample performances of our method compared to benchmarks as the network grows sparser under fixed $n$.  Despite in Simulation \ref{subsec::simulation-1}, we tested different network sparsity settings (see Supplemental Material), it would still be interesting to more directly illustrate the impact of $\rho_n$ for each fixed network size.  The simulation set up carries over the same set of graphon models, motif shapes and compared methods from Simulation \ref{subsec::simulation-1}.  Here, for simplicity, we only tested $n = 80, 160$ and varied $\rho_n$ in $\{1\textrm{ (``dense'')},n^{-1/4},n^{-1/2},n^{-1}\}$.

\begin{figure}[htbp!]
	\begin{adjustwidth}{-\oddsidemargin-1.5in}{-\rightmargin-2in}
		\centering
		\includegraphics[width=0.4\textwidth]{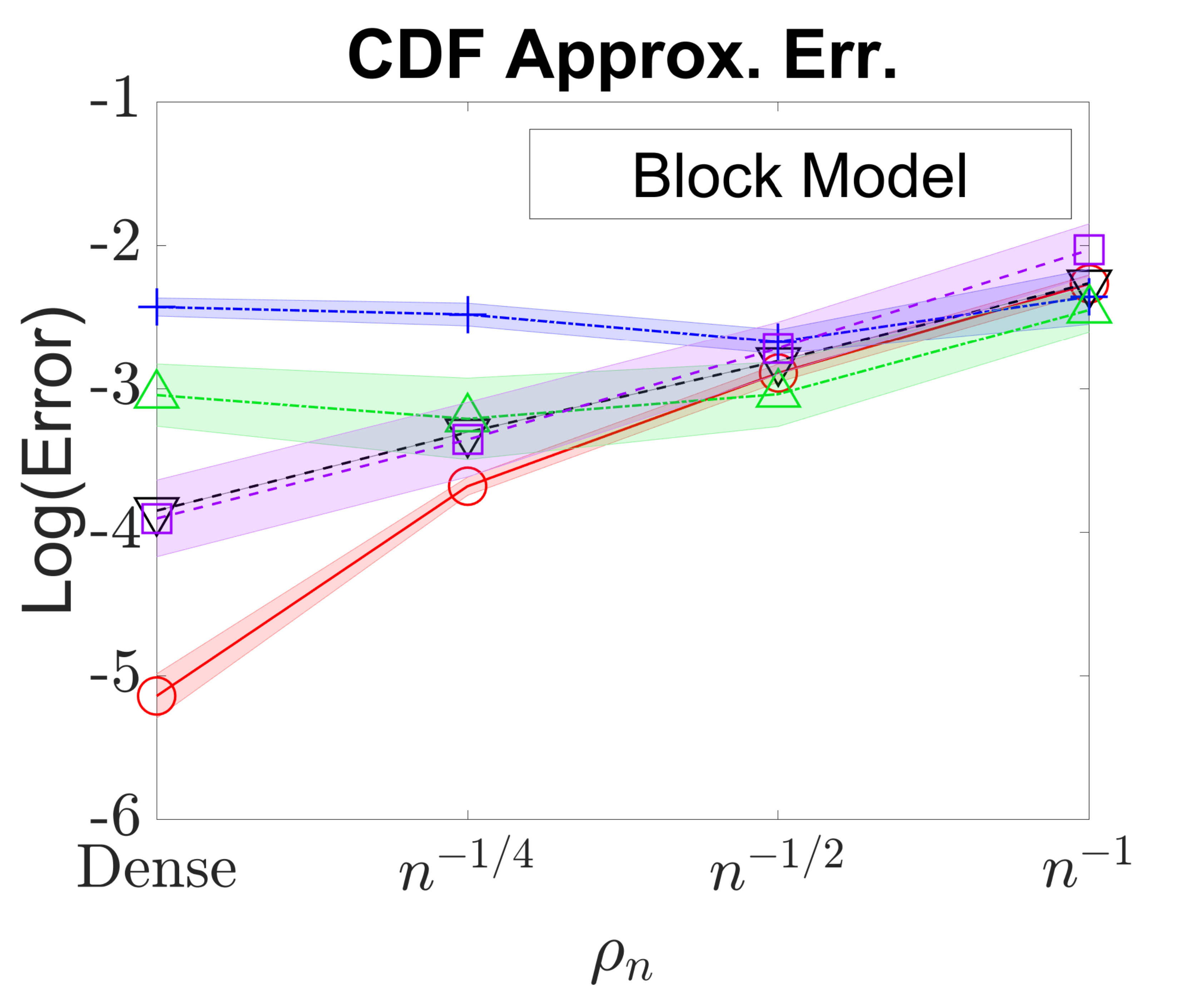}
    	\includegraphics[width=0.4\textwidth]{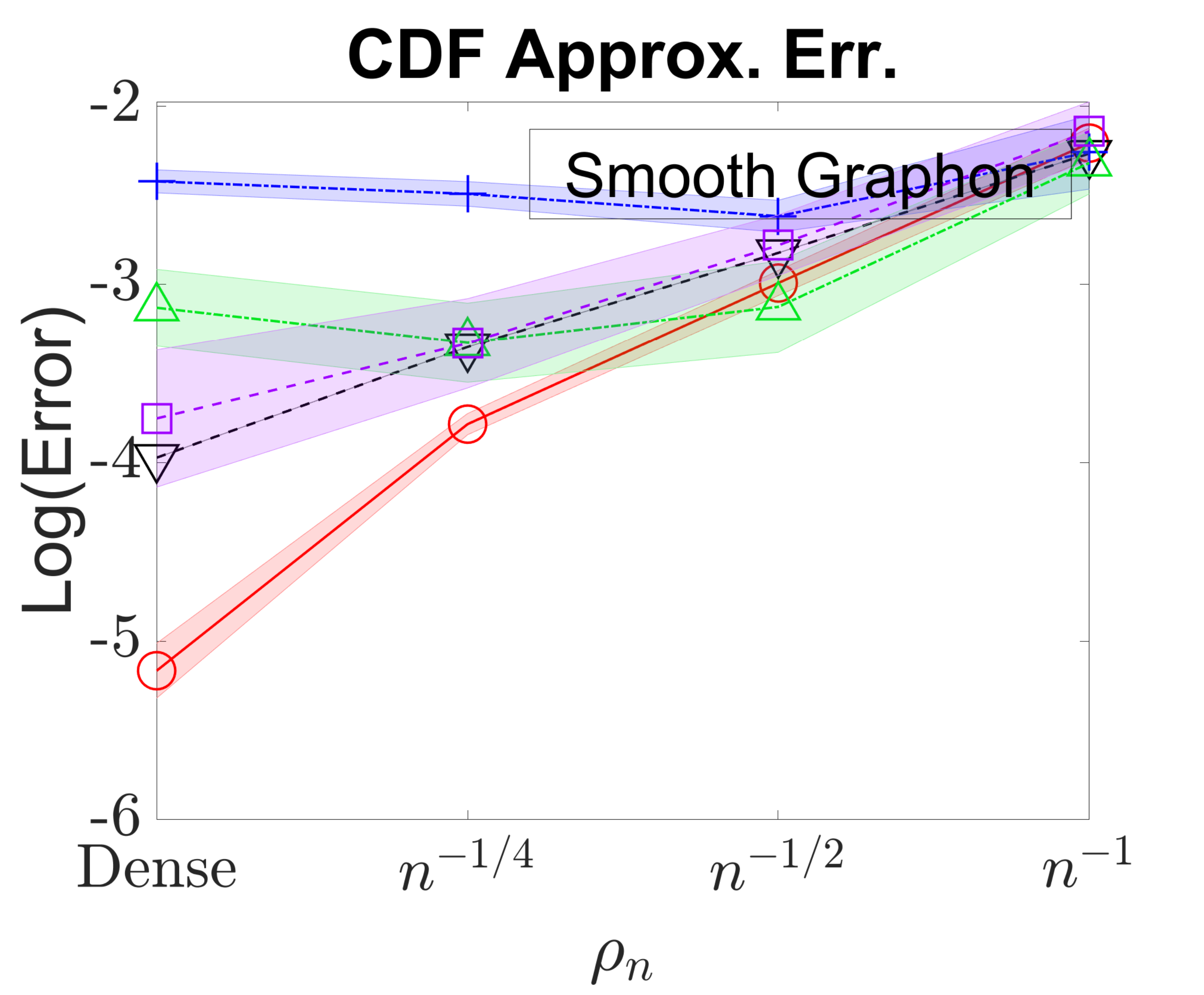}
    	\includegraphics[width=0.4\textwidth]{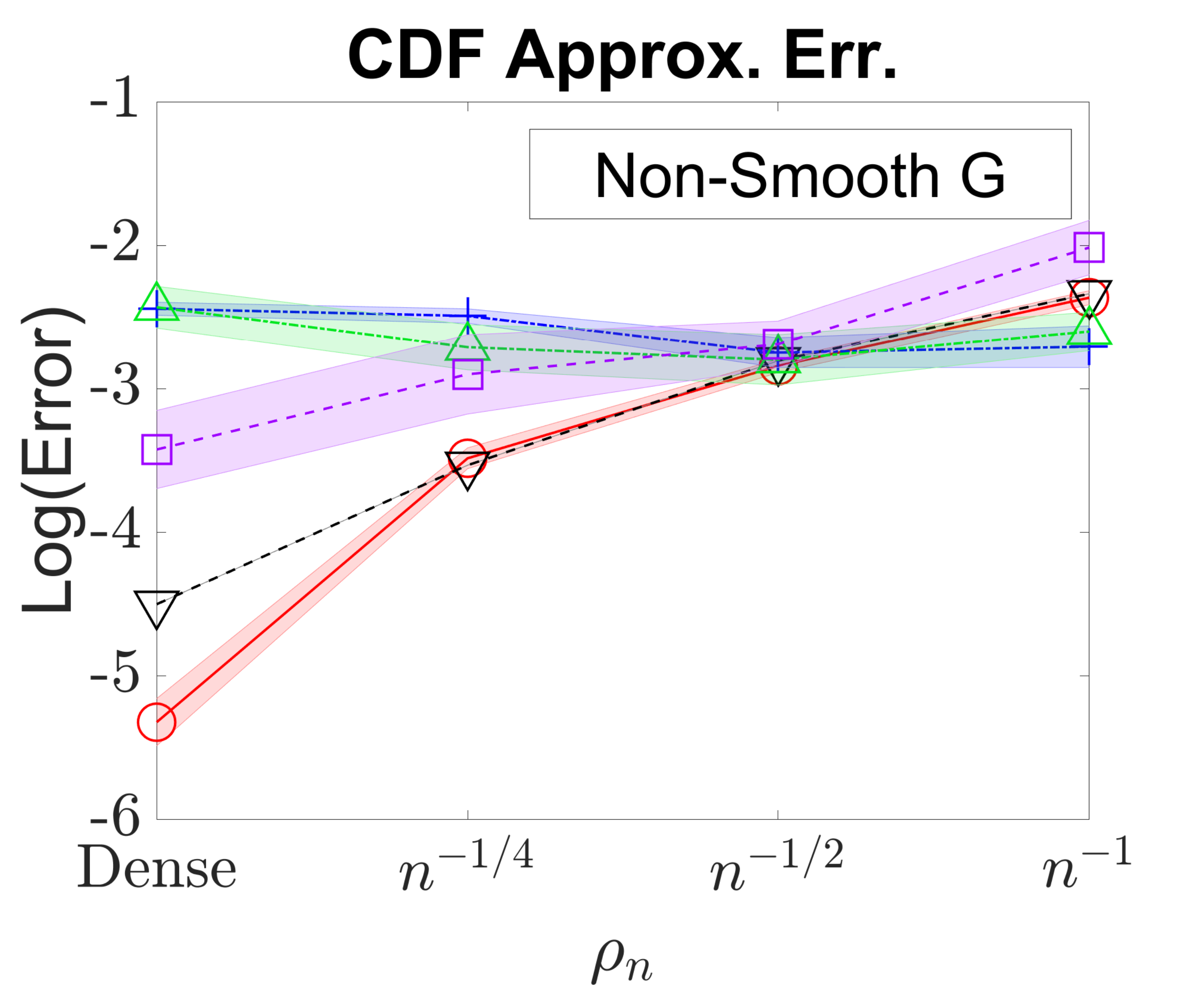}\\
    	\includegraphics[width=0.4\textwidth]{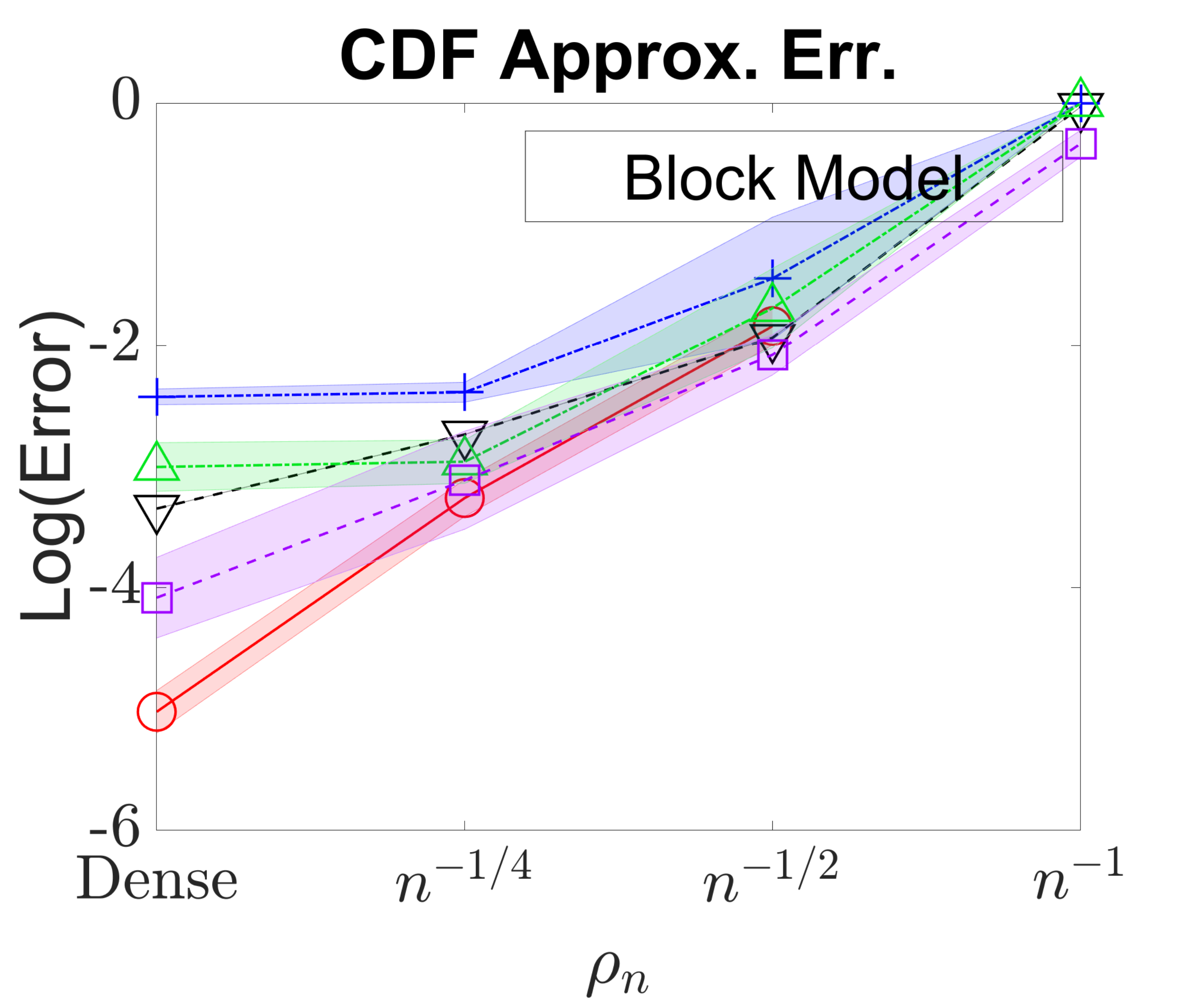}
    	\includegraphics[width=0.4\textwidth]{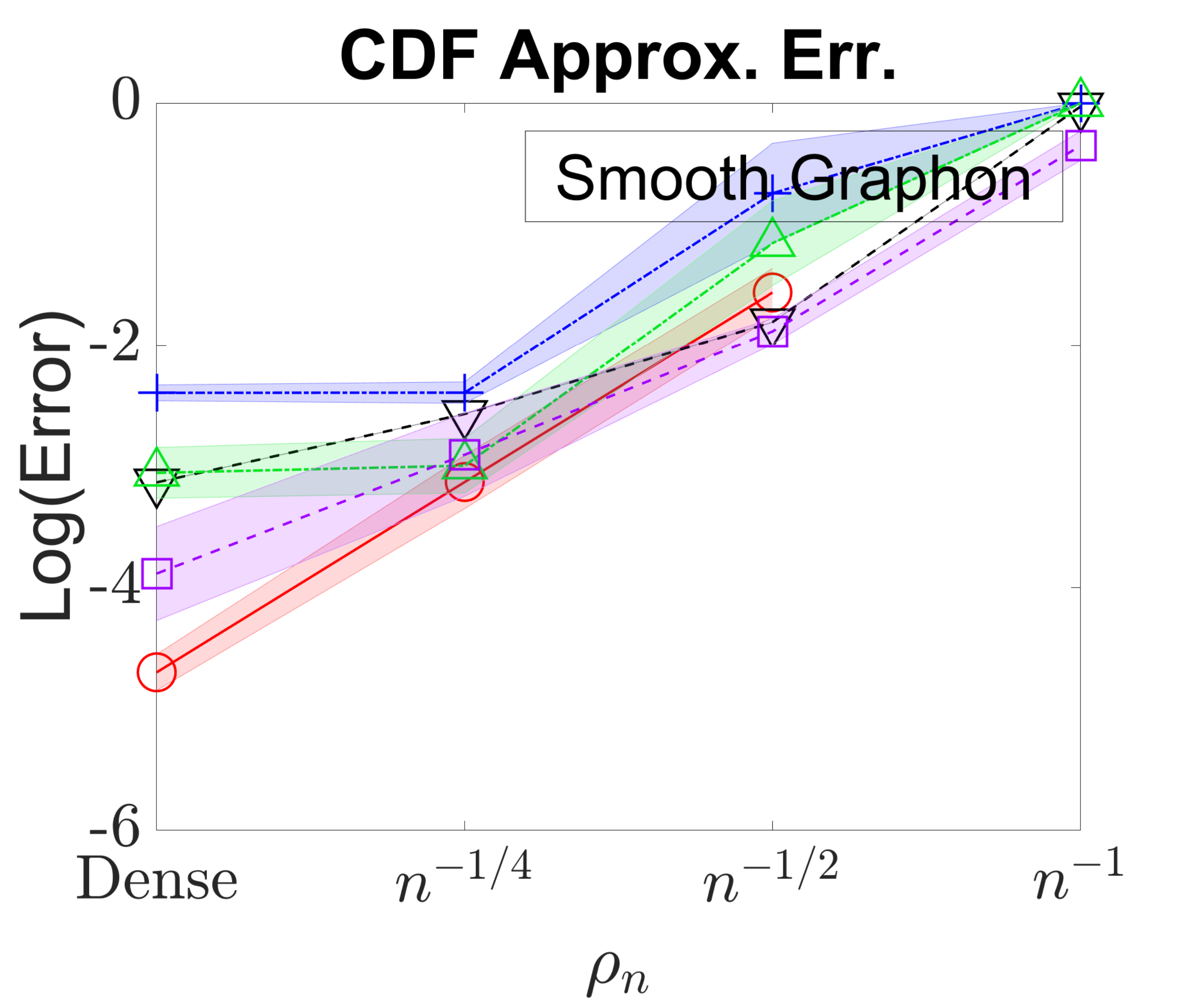}
    	\includegraphics[width=0.4\textwidth]{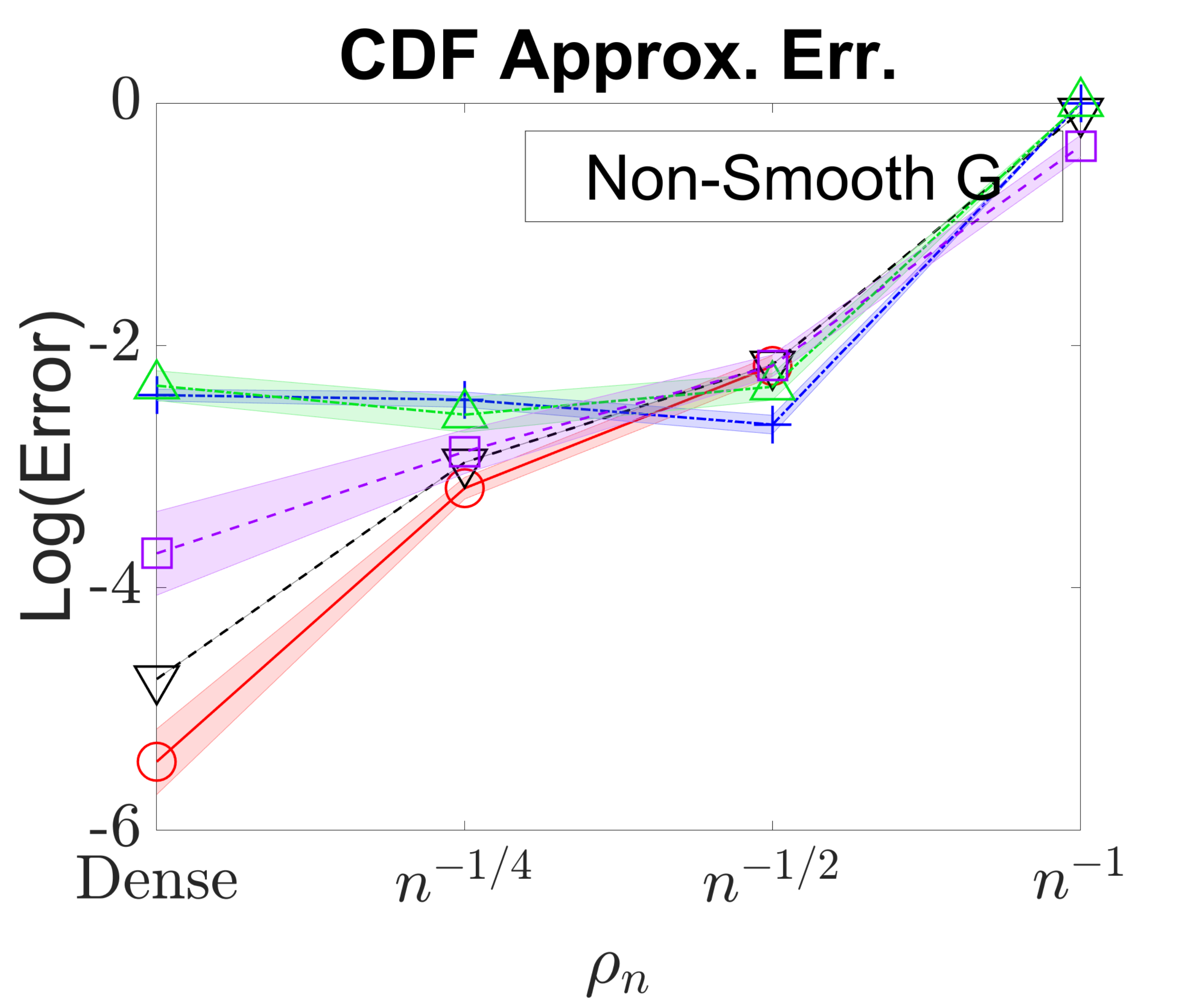}\\
    	\includegraphics[width=0.4\textwidth]{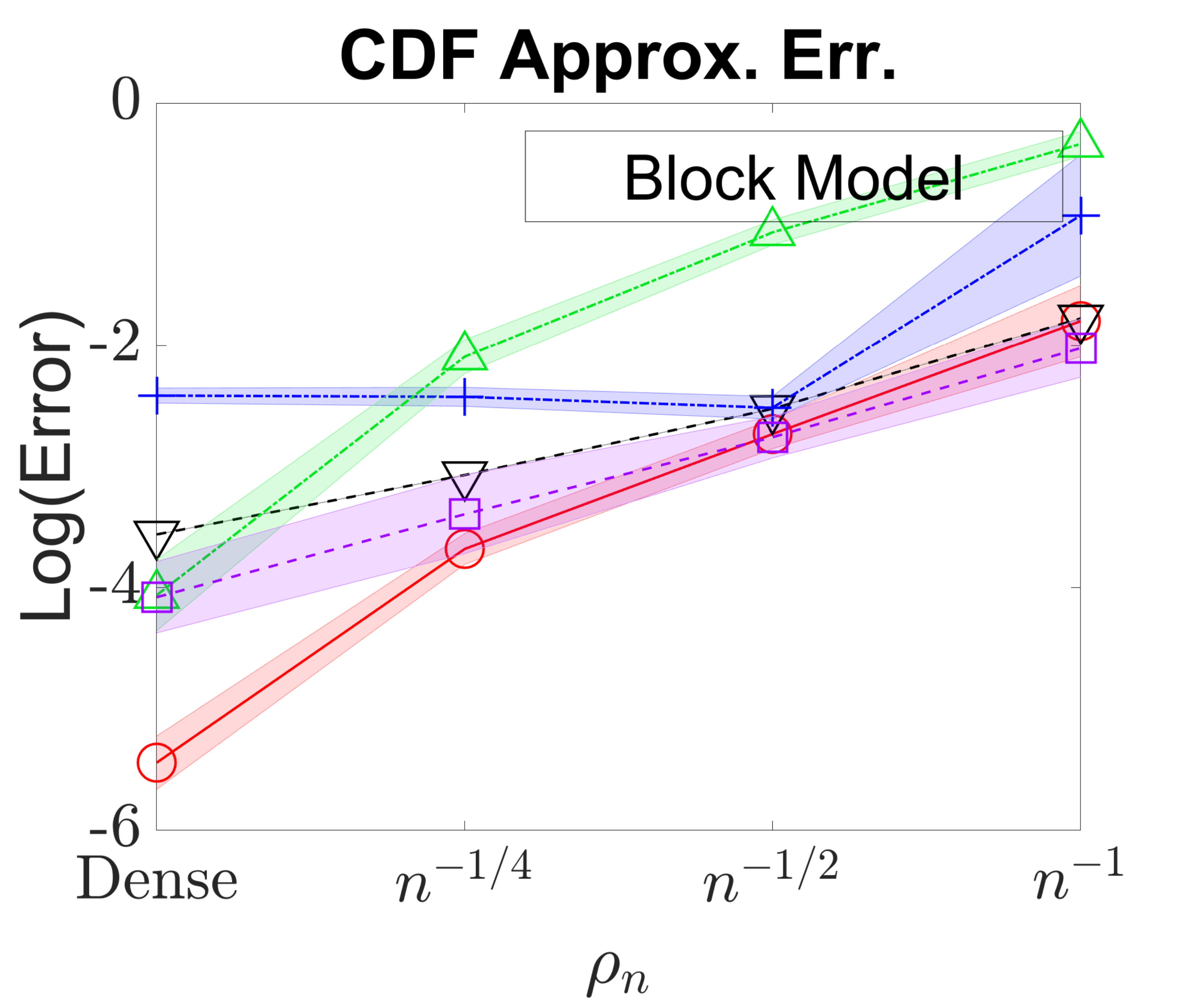}
    	\includegraphics[width=0.4\textwidth]{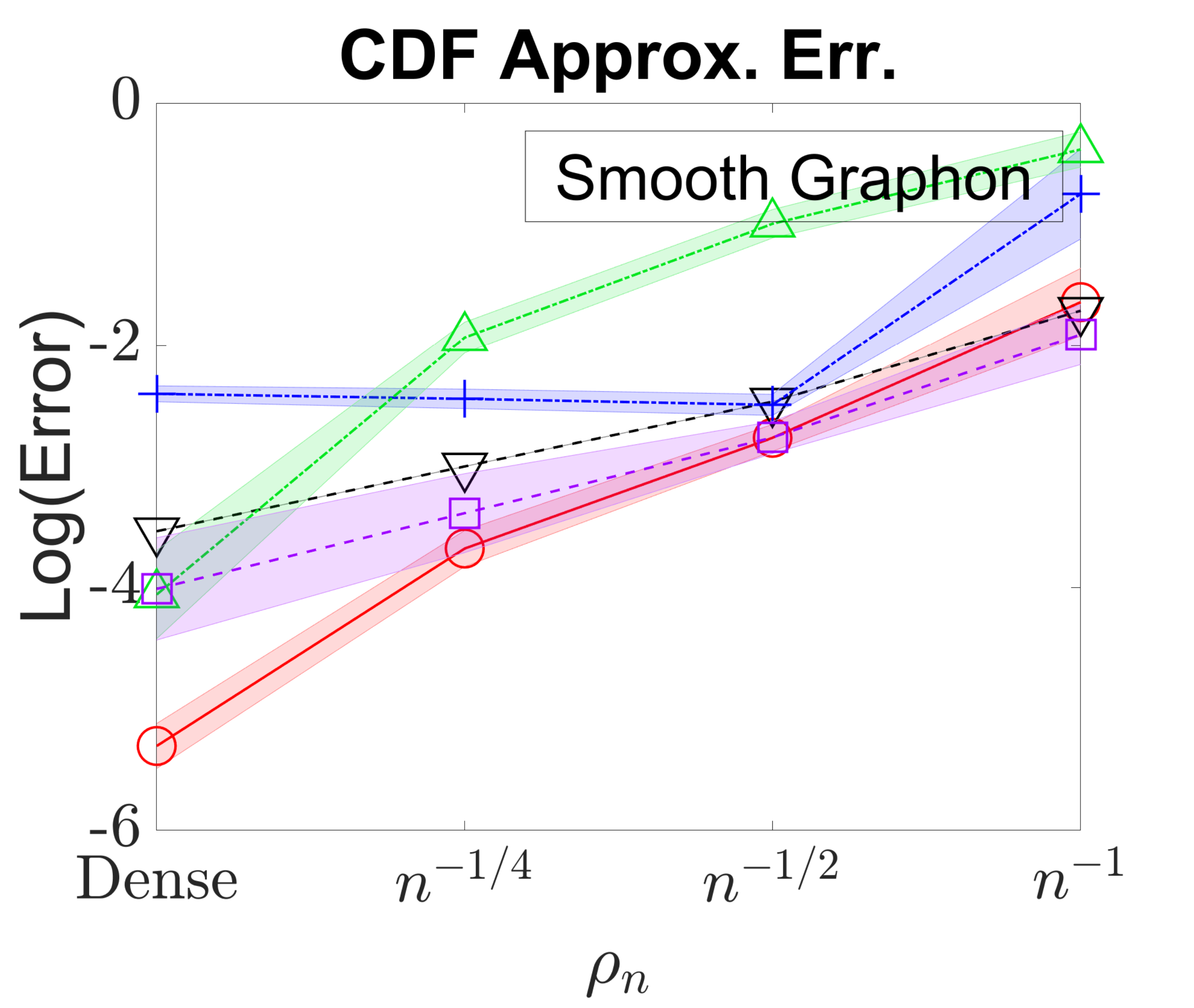}
    	\includegraphics[width=0.4\textwidth]{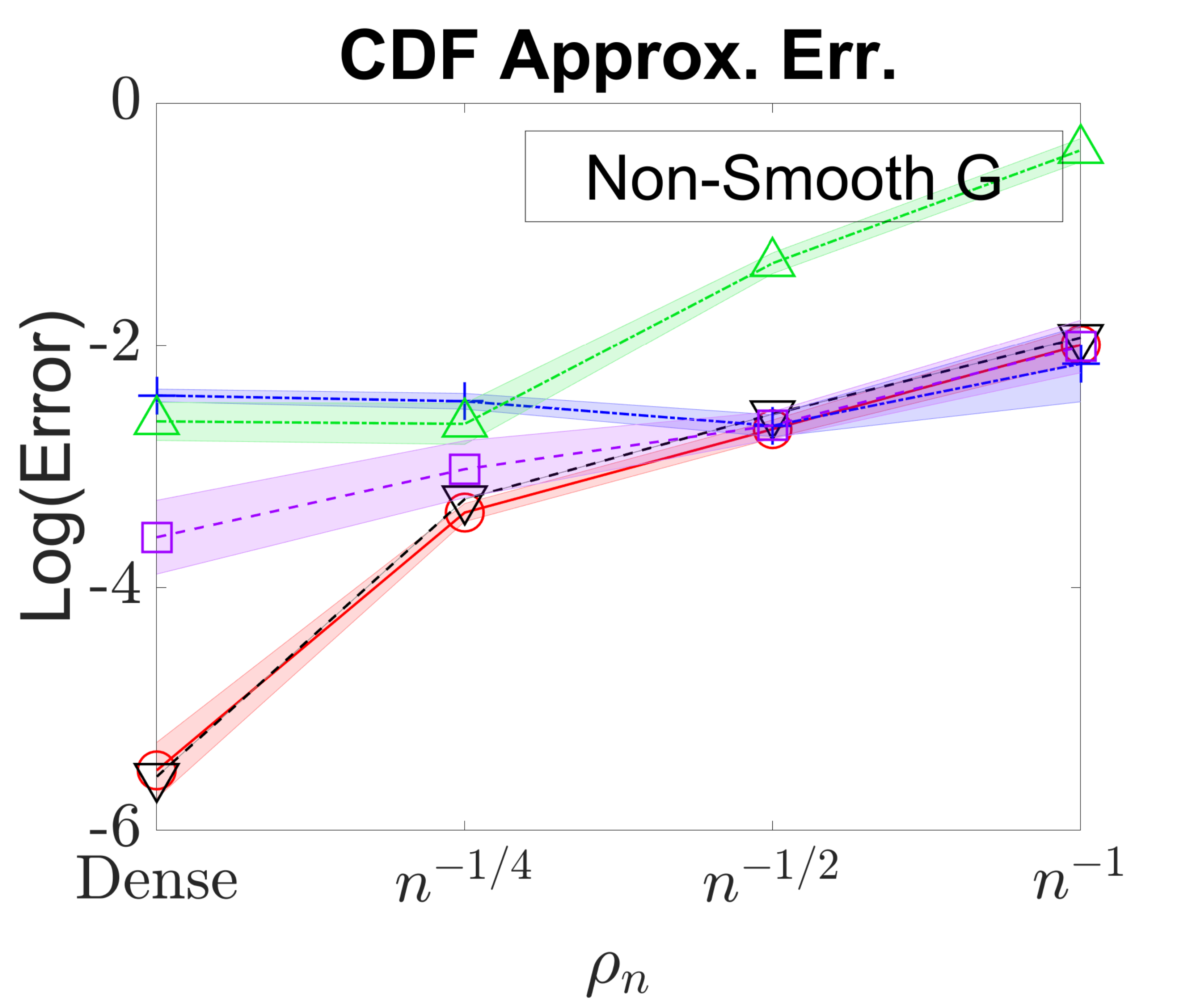}\\
    	\includegraphics[width=0.4\textwidth]{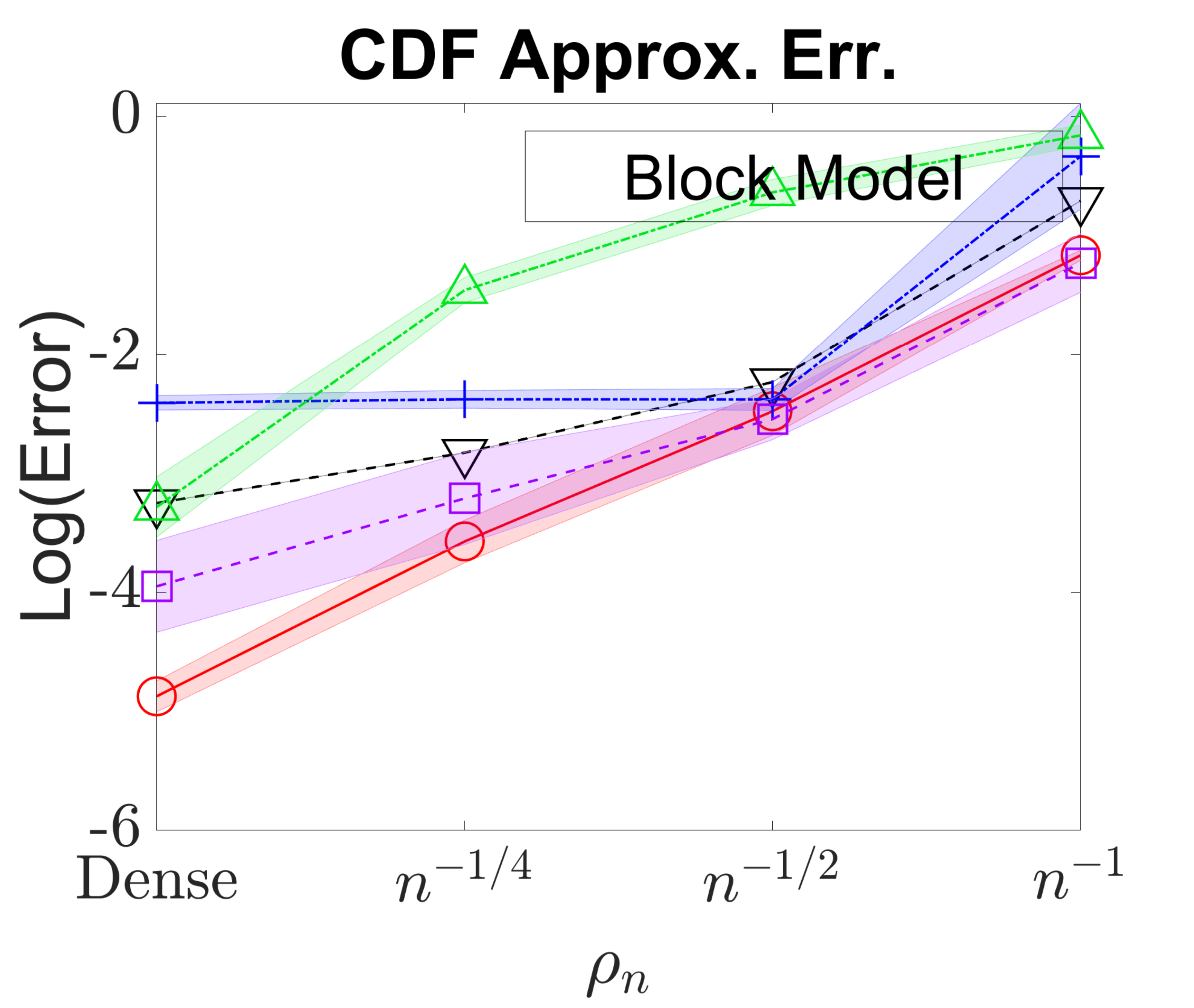}
    	\includegraphics[width=0.4\textwidth]{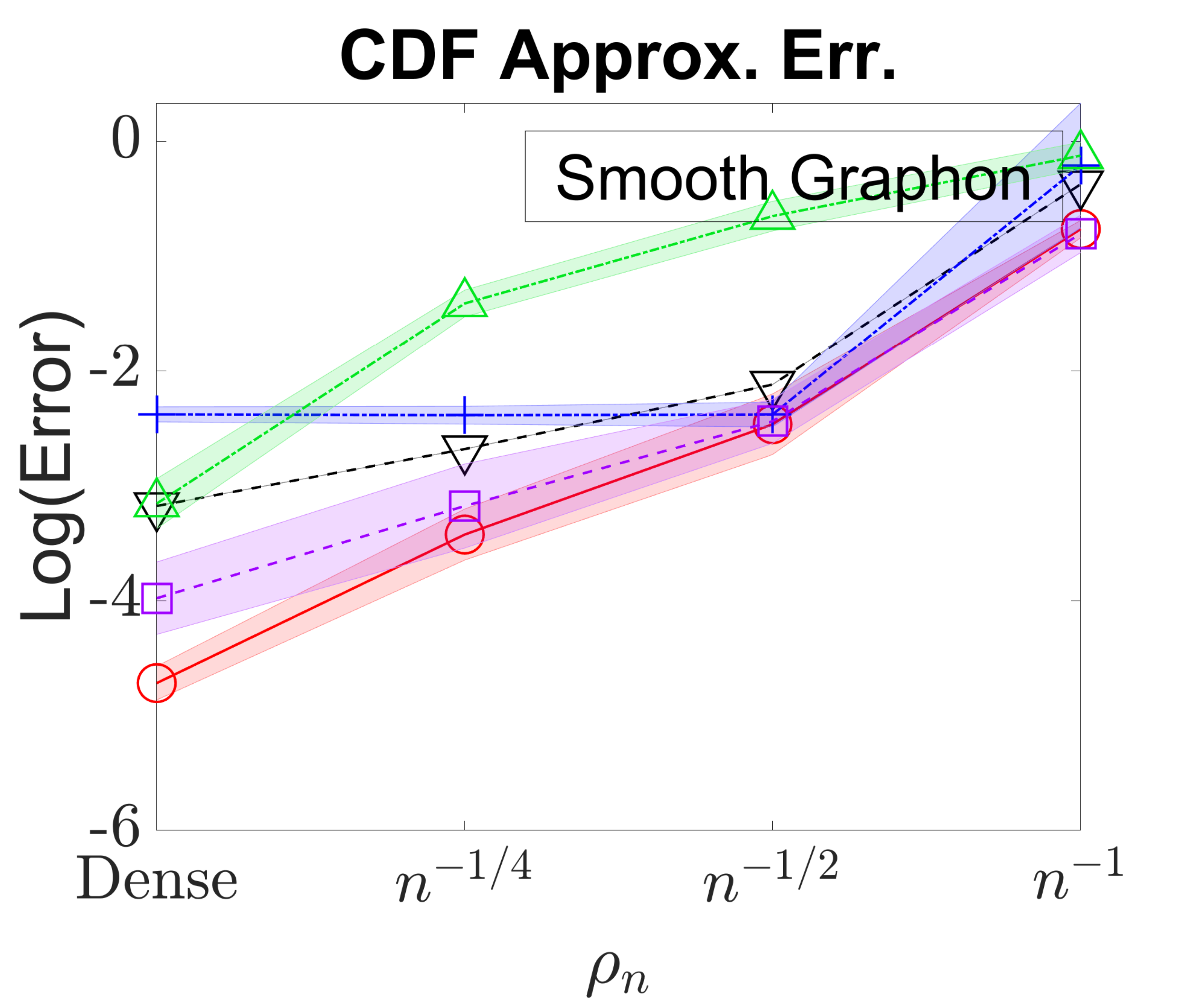}
    	\includegraphics[width=0.4\textwidth]{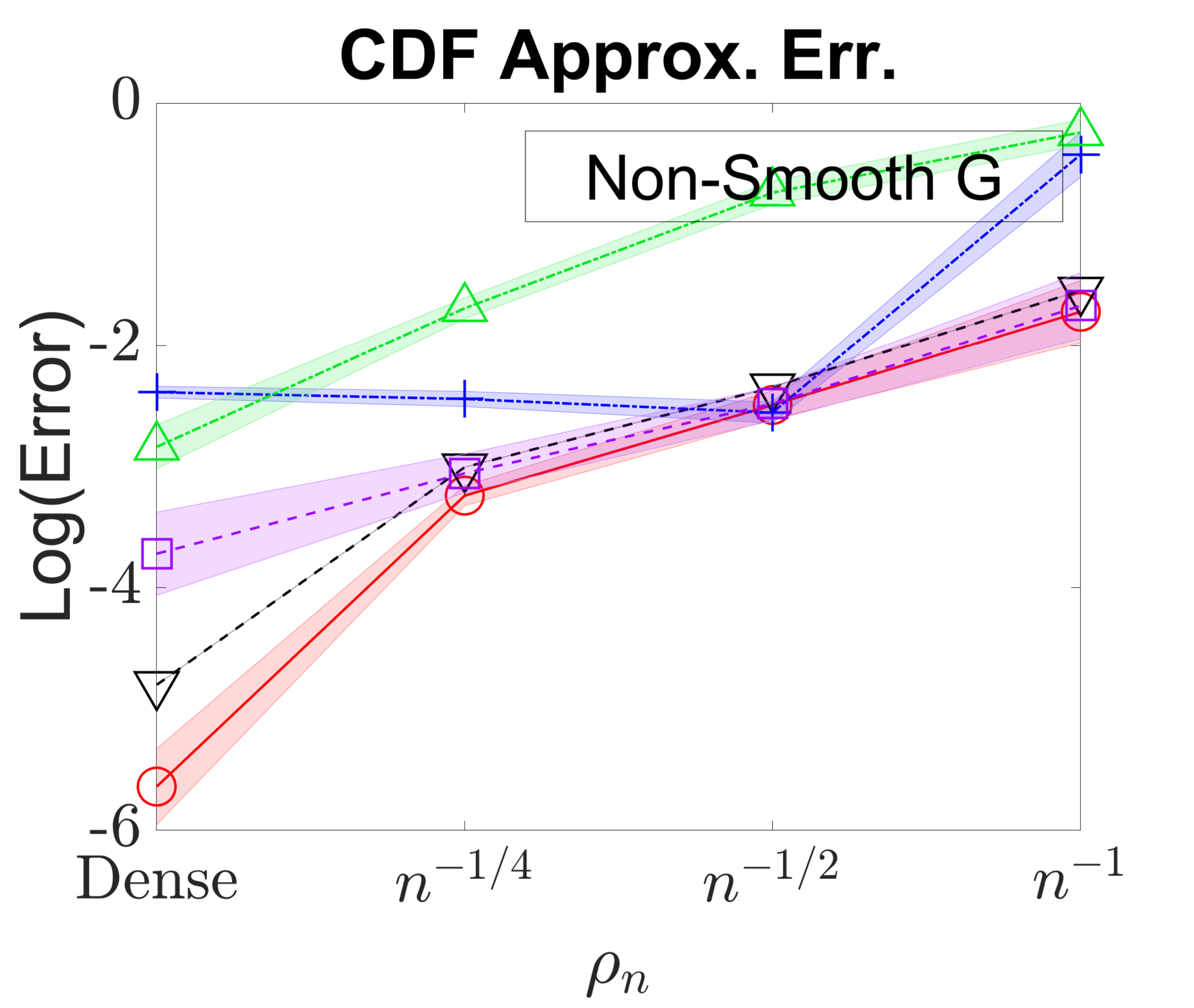}
	\end{adjustwidth}
	\caption{ {Impact of sparsity on approximation errors, $n=160$.} Both axes are log(e)-scaled. {\bf Motifs:} row 1: {\tt Edge}; row 2: {\tt Triangle}; row 3: {\tt Vshape}{; row 4: {\tt ThreeStar}}.      \tred{Red solid curve marked circle}: our method (empirical Edgeworth);
	black dashed curve marked down-triangle: $N(0,1)$ approximation;
	\textcolor{green}{green dashed curve marked up-triangle}: re-sampling of $A$ in \cite{green2017bootstrapping}; \tblue{blue dashed curve marked plus}: \cite{bhattacharyya2015subsampling} sub-sampling $\asymp n$ nodes; \textcolor{magenta}{magenta dashed line with square markers}: ASE plug-in bootstrap in \cite{levin2019bootstrapping}.}
	\label{fig::numerical-sparsity-main}
\end{figure}

Figure \ref{fig::numerical-sparsity-main} shows the CDF approximation errors under different $\rho_n$ settings for $n=160$.  Aligned with our theory's prediction, we observe that as the network grows sparser, our method's performance depreciates and gradually regresses to the performance of normal approximation.  Due to page limit, we sink the approximate error plots for $n=80$ and the time cost plots for both $n$ settings to Supplemental Materials.

\section{Discussion}
\label{sec::discussion}
{
Our results do not cover the case where $g_1(X_1)$ is lattice and $\rho_n\asymp 1$.  An ad-hoc remedy is to simply introduce artificial missing links by sparsifying $A$:
$$
    \tilde A_{ij} := \tilde A_{ji} := 
    \begin{cases}
        A_{ij}=A_{ji}, \textrm{ with probability }1-\tilde\rho_n\\
        0, \textrm{ with probability }\tilde\rho_n
    \end{cases}
$$
where we set $\tilde\rho_n\asymp  (\log n)^{-1}$.  One can then make inferences about the population network moment $\tilde\rho_n\cdot \mu_n$ using $\tilde A$ as the input data (notice $\tilde\rho_n$ is known).  This reinstates the $(\log n)^{-1}$ sparsification that we need to overcome the latticeness at the price of a very minor information loss.
}

The Edgeworth expansion we derived for Bernoulli $A_{ij}|W_{ij}$ distributions can be readily extended to general weighted networks, where the conditionally independent $A_{ij}|W_{ij}$ distributions may either depend on $W_{ij}$ or not.  A distinct feature of our setting is that the edge-wise observational errors are a contributing component of $\hat{T}_n$ that smooths the distribution.  In contrast to matrix estimation problems, where such noise is to be suppressed \citep{candes2010matrix,xia2019statistical}, a moderate amount of tailedness can strengthen the smoothing effect in $A|W$ and might improve finite sample performances.  Notice that similar to \cite{bhattacharyya2015subsampling,green2017bootstrapping,levin2019bootstrapping,lunde2019subsampling}, in our main theorems, we omitted finite-moment assumptions on $h(\cdot)$ since it is naturally bounded in network settings.

{
This paper focuses on studying the marginal randomness in $A$ jointly contributed by the randomness in $W$ and $A|W$.  In this study, we take the sparsity-scaled graphon $\rho_n\cdot f$ as the population and the graphon feature $\mu_n$ as the ultimate inference goal.  Our approach is nonparametric and directly approximates $F_{\hat T_n}$ \emph{without requiring} a graphon estimation $\hat W$.  If one regards $W$ as the population and wants to make inference for $U_n$, she would need a CDF approximation to $(\hat U_n-U_n)|X_1,\ldots,X_n$.  This distribution is asymptotically normal as has been described by \eqref{eqn::Gaussian-term-convergence} in our Lemma \ref{lemma::term-approx}-\eqref{lemma::term-approx-Delta-hat-conditional-normal}.  However, estimating the normal variance typically \emph{requires} a graphon estimation $\hat W$\footnote{The expression of $\sigma_w^2$ contains $W_{ij}^2$ terms originated from ``$W_{ij}(1-W_{ij})$'' terms, which could not be estimated without a graphon estimation.}.  On the other hand, the bootstrapping of $\hat T_n|X_1,\ldots,X_n$ would also (seemingly unavoidably) be a parametric bootstrap.  
In view of Lemma \ref{lemma::term-approx}-\eqref{lemma::term-approx-Delta-hat-conditional-normal}, asympotically 
\begin{equation}
    (\rho_n\cdot n)^{1/2}\cdot \frac{\hat U_n-U_n}{\sigma_n} \stackrel{d}\approx N(0,\sigma_w^2)
    \label{eqn::discussion::Un-hat-minus-Un}
\end{equation}
given $X_1,\ldots,X_n$, where recall that $\sigma_w\asymp 1$.
However, the minimax rate for sparse graphon estimation (see \cite{gao2016optimal, klopp2017oracle}) is
$$
    \textrm{Rescaled MSE: }(\rho_n\cdot n)^{-2}\cdot\|\hat W - W\|_F^2 \asymp (\rho_n\cdot n)^{-1}\cdot \log n
$$
If we use this error bound to control the estimation error of $\sigma_w^2$, then this yields an error of $|\hat\sigma_w^2-\sigma_w^2|\asymp n^{-1}\|\hat W-W\|_F\asymp \rho_n^{1/2}\cdot n^{-1/2}\cdot \log^{1/2} n$. 
This error may dominate the $n^{-1/2}$ correction term in an Edgeworth expansion for dense networks (e.g., Cramer's condition holds and $\rho_n\asymp 1$).  Moreover, the minimax graphon estimation rate has not yet been achieved by any polynomial-time algorithm (see \cite{zhang2017estimating, xu2018rates} for comments), and using a practically feasible $\hat W$ would cause an error $\gg n^{-1/2}$, ignoring $\rho_n$ and log.
Therefore, it might be challenging to accurately approximate the distribution of the LHS of \eqref{eqn::discussion::Un-hat-minus-Un} beyond asymptotic normality.
Our observation here echos the common practice in network bootstrap literature \citep{bhattacharyya2015subsampling, green2017bootstrapping, levin2019bootstrapping} that they unanimously focus on the marginal distribution of $\hat U_n$, rather than $(\hat U_n-U_n)|X_1,\ldots,X_n$\footnote{For example, in \citet{levin2019bootstrapping}, the authors used a low-rank decomposition of $A$, which directly leads to an estimated $\hat W$.  But they also solely focused the marginal distribution of $U_n$ (in our notation system).}.
}

A retrospection on our simulation setting provides an interesting insight.  In fact, the population Edgeworth expansion provides a much more efficient Monte Carlo procedure for simulating the true distribution $F_{\hat{T}_n}$.  Indeed, estimating $\xi_1$, $\ep[g_1^3(X_1)]$ and $\ep[g_1(X_1)g_1(X_2)g_2(X_1,X_2)]$ with $n_{\textrm{MC}}\asymp n$ Monte Carlo samples yields a CDF approximation rate of $\Ohighorderbound=o(n^{-1/2})$ when $\rho_n$ satisfies the conditions of Theorem \ref{thm::main-theorem}.  This is much more efficient than the empirical CDF which requires $n_{\mathrm{MC}}\succeq n^2$ to achieve the same accuracy order. 

In the application of our results, we focused on node sampling network bootstraps.  It is an interesting future work to investigate the higher-order accuracy properties of other schemes, such as sub-graph sampling \citep{bhattacharyya2015subsampling} and (artificially) weighted bootstrap \citep{levin2019bootstrapping}.  Also comprehensive numerical comparisons of different schemes under various settings would certainly be interesting for practitioners.
{
As a closely related point, this paper studies the \emph{complete} noisy U-statistic, ``complete'' in the sense that $(i_1,\ldots,i_r)$ ranges over all $\binom{n}r$ possible $r$-tuples.  As one of the anonymous referees pointed out, evaluating the moment corresponding to an $r$-node motif would cost $O(n^r)$, which is still expensive for large $n$.  Papers \cite{blom1976some, nowicki1988subgraph, chen2019randomized, kong2020design, song2019approximating} investigated this topic for the conventional noiseless U-statistic setting, where the same challenge exists, such as formulating the Edgeworth expansion for incomplete U-statistics.  They study noiseless incomplete U-statistics, and \cite{bhattacharyya2015subsampling} proposed a ``subgraph subsampling'' scheme (their scheme (a)) that computes noisy incomplete U-statistics, which we denote by $\hat U_n^{\rm (Incomplete)}$ for the network setting.  Formally, we have
$$
    \hat U_n^{{\rm (Incomplete)}}
    :=
    \frac{\sum_{1\leq i_1<\cdots<i_r\leq n} I_{i_1,\ldots,i_r}\cdot h(A_{i_1,\ldots,i_r})}{\sum_{1\leq i_1<\cdots<i_r\leq n} I_{i_1,\ldots,i_r}}
$$
where $I_{i_1,\ldots,i_r}$'s are random variables independent of the network data.  These $I_{i_1,\ldots,i_r}$'s can be i.i.d. Bernoulli, multinomial (if a given proportion of sub-sampled motifs is desired), or other reasonable sampling scheme distributions.
It would be an interesting future research to carefully explore and quantify the self-smoothing effect for $\hat U_n^{\rm (Incomplete)}$.
}


%

\section*{Acknowledgments}
{We thank the Editor, the Associate Editor and three referees for their insightful and constructive comments that led to significant improvements in many aspects of this paper.} 
We thank Professor Bing-Yi Jing, Professor Junhui Wang and Professor Ji Zhu for helpful comments; 
Vincent Q. Vu for helpful feedback; 
Linjun Zhang and Yunpeng Zhao for constructive suggestions; 
and Keith Levin and Tianxi Li for discussions on network bootstraps.
Dong Xia's research was partially supported by Hong Kong RGC Grant ECS 26302019 and Adobe Research Gift Award.

%
\begin{supplement}
\sname{Supplement for}\label{suppA}
\stitle{``Edgeworth expansions for network moments''}
\slink[url]{URL to be added}
\sdescription{The supplementary material contains: (1). Definition of $\sigma_w$ in Lemma \ref{lemma::term-approx}-\eqref{lemma::term-approx-Delta-hat-conditional-normal}; (2). All proofs; and (3). Additional simulation results and accompanying interpretations.}
\end{supplement}

\bibliographystyle{imsart-nameyear}
\bibliography{Edgeworth-ref}

\end{document}